\documentclass[a4paper,10pt]{article}

\def\zibreport{0}
\def\guidelines{0}
\def\expensiveFigures{1}
\def\longtitle{The~SCIP~Optimization~Suite~10.0}
\def\shortfunding{
The work for this article was partially supported through the \emph{Research Campus MODAL} funded by the German Federal Ministry of Research,
Technology, and Space (fund numbers \emph{05M14ZAM}, \emph{05M20ZBM}) and the \emph{Deutsche Forschungsgemeinschaft (DFG)} through the \emph{DFG Cluster of Excellence MATH+}.
It has also been partly supported the project \emph{Making Mixed-Integer Programming Solvers Smarter and Faster using Network Matrices} (with
project number \emph{OCENW.M20.151}) of the research programme \emph{NWO Open Competition Domain Science -- M}
which is (partly) financed by the Dutch Research Council (NWO).
}

\makeatletter
\DeclareRobustCommand*{\escapeus}[1]{%
  \begingroup\@activeus\scantokens{#1 }\endgroup}
\begingroup\lccode`\~=`\_\relax
   \lowercase{\endgroup\def\@activeus{\catcode`\_=\active \let~\_}}
\makeatother

\usepackage[utf8]{inputenc}
\usepackage[font=small, labelfont=bf, margin=1cm]{caption}
\usepackage{amsmath,amsfonts,amssymb}
\usepackage{tabularx,booktabs}
\usepackage{url}
\usepackage{rotating}
\usepackage{microtype}

\usepackage[textsize=small,textwidth=3.3cm]{todonotes}

\usepackage{xspace}
\usepackage{etex}
\usepackage[inline]{enumitem}
\usepackage{dsfont}
\usepackage{mathtools}
\usepackage[linesnumbered,ruled,vlined]{algorithm2e}

\SetCommentSty{AlgoCommentStyle}
\usepackage{changes}
\usepackage{longtable}
\usepackage{tabularx}
\usepackage{varwidth}
\usepackage{listings}
\usepackage{diagbox}
\usepackage{makecell}

\usepackage{geometry}
\usepackage{amsthm}
\usepackage{multirow}

\usepackage[breaklinks,colorlinks=true,citecolor=black,linkcolor=black,urlcolor=blue,hyperfootnotes=false,
  pdftitle={\longtitle},
  pdfsubject={\longtitle}]{hyperref}

\usepackage[capitalize]{cleveref}
\usepackage[numbers]{natbib}
\usepackage{doi}

\usepackage{tikz}
\usepackage{tikz-3dplot}
\usepackage{tikzscale}
\usetikzlibrary{arrows}
\usetikzlibrary{shapes}
\usetikzlibrary{snakes}
\usetikzlibrary{patterns}
\usetikzlibrary{backgrounds,topaths}
\usetikzlibrary{matrix,chains,scopes,positioning,arrows,fit}

\usepackage{pgfplots}
\pgfplotsset{compat=1.15}
\usepackage{subcaption}
\usepackage{textcomp}
\usepgfplotslibrary{groupplots}
\usepgfplotslibrary{units}

\usepackage{graphicx}

\makeatletter
\newcommand*\rel@kern[1]{\kern#1\dimexpr\macc@kerna}
\newcommand*\widebar[1]{%
  \begingroup
  \def\mathaccent##1##2{%
    \rel@kern{0.8}%
    \overline{\rel@kern{-0.8}\macc@nucleus\rel@kern{0.2}}%
    \rel@kern{-0.2}%
  }%
  \macc@depth\@ne
  \let\math@bgroup\@empty \let\math@egroup\macc@set@skewchar
  \mathsurround\z@ \frozen@everymath{\mathgroup\macc@group\relax}%
  \macc@set@skewchar\relax
  \let\mathaccentV\macc@nested@a
  \macc@nested@a\relax111{#1}%
  \endgroup
}
\makeatother

\newcommand{\LP}{{LP}\xspace}
\newcommand{\LPs}{{LPs}\xspace}

\newcommand{\CIPs}{{CIPs}\xspace}
\newcommand{\MIP}{{MIP}\xspace}
\newcommand{\MIPs}{{MIPs}\xspace}
\newcommand{\MILP}{{MILP}\xspace}
\newcommand{\MILPs}{{MILPs}\xspace}
\newcommand{\NLP}{{NLP}\xspace}

\newcommand{\MINLP}{{MINLP}\xspace}
\newcommand{\MINLPs}{{MINLPs}\xspace}

\newcommand{\T}{\top}

\newcommand{\defi}{\coloneqq}

\DeclareMathOperator{\conv}{conv}

\newcommand{\linobj}{c}
\newcommand{\nonlinobj}{f}
\newcommand{\linmatrix}{A}
\newcommand{\nonlincons}{g}
\newcommand{\rhs}{b}
\newcommand{\lb}{\ell}
\newcommand{\ub}{u}

\newcommand{\consindex}{\mathcal{M}}
\newcommand{\varindex}{\mathcal{N}}
\newcommand{\intvarindex}{\mathcal{I}}

\newcommand{\N}{\mathds{N}}
\newcommand{\R}{\mathds{R}}
\newcommand{\Q}{\mathds{Q}}
\newcommand{\Z}{\mathds{Z}}
\newcommand{\Rinf}{\ensuremath{\widebar{\mathds{R}}}\xspace}

\usepackage[binary-units=true]{siunitx}

\newcommand{\cleaninst}{all}
\newcommand{\affected}{affected}
\newcommand{\alloptimal}{{both-solved}\xspace}
\newcommand{\difftimeouts}{{diff-timeouts}\xspace}
\newcommand{\bracket}[2]{[#1,#2]}

\usepackage{xparse}
\ExplSyntaxOn
\def\myround#1{\num{\fp_eval:n {round(#1, 2)}}}
\ExplSyntaxOff

\definecolor{c1}{HTML}{000060}
\definecolor{c2}{HTML}{0000FF}
\definecolor{c3}{HTML}{36648B}
\definecolor{c4}{HTML}{4682B4}
\definecolor{c5}{HTML}{5CACEE}
\definecolor{c6}{HTML}{00FFFF}
\definecolor{c7}{HTML}{008888}
\definecolor{c8}{HTML}{00DD99}
\definecolor{c9}{HTML}{527B10}
\definecolor{c10}{HTML}{7BC618}
\definecolor{c11}{HTML}{33AA00}

\definecolor{scipoldcol}{HTML}{36648B}
\definecolor{scipnewcol}{HTML}{7BC618}

\newcommand{\solver}[1]{\textsc{#1}\xspace}
\newcommand{\scipopt}{\scip Optimization Suite\xspace}
\newcommand{\scipprevversion}{9.0\xspace}
\newcommand{\scipversion}{10.0\xspace}

\newcommand{\scipoptv}{\scipopt~\scipversion{}\xspace}
\newcommand{\scip}{\solver{SCIP}}
\newcommand{\scipv}{\solver{SCIP}~\scipversion{}\xspace}
\newcommand{\scippv}{\solver{SCIP}~\scipprevversion{}\xspace}
\newcommand{\soplex}{\solver{SoPlex}}
\newcommand{\soplexversion}{8.0}
\newcommand{\soplexv}{\solver{SoPlex}~\soplexversion{}\xspace}
\newcommand{\papilo}{\solver{PaPILO}}
\newcommand{\papiloversion}{3.0}
\newcommand{\papilov}{\solver{PaPILO}~\papiloversion{}\xspace}
\newcommand{\mipdd}{\solver{MIP-DD}}
\newcommand{\zimpl}{\solver{Zimpl}} 

\newcommand{\ug}{\solver{UG}}
\newcommand{\presollib}{\solver{PaPILO}}

\newcommand{\gcg}{\solver{GCG}}
\newcommand{\gcgversion}{4.0}
\newcommand{\gcgv}{\gcg~\gcgversion{}\xspace}

\newcommand{\scipsdp}{\solver{SCIP-SDP}}

\newcommand{\qsoptex}{\solver{QSopt\_ex}}

\newcommand{\param}[1]{\texttt{#1}\xspace}
\newcommand{\method}[1]{\texttt{#1}\xspace}

\newcommand{\plugin}[1]{\texttt{\escapeus{#1}}\xspace}

\newcommand{\cplex}{\solver{CPLEX}}

\newcommand{\nbsc}[1]{\mbox{#1}\xspace}
\newcommand{\miplib}{\nbsc{MIPLIB}}

\newcommand{\coral}{\nbsc{COR@L}}

\newcommand{\minlplibtwo}{\nbsc{MINLPLib}}

\definecolor{darkgreen}{HTML}{008800}

\newcommand{\floor}[1]{\lfloor #1\rfloor}

\newcommand{\fa}{\text{ for all }}

\newcommand{\defin}{\ensuremath{\coloneqq}}

\newcommand{\bliss}{\solver{bliss}}

\newcommand{\nauty}{\solver{nauty}}
\newcommand{\dejavu}{\solver{dejavu}}
\newcommand{\sassy}{\solver{sassy}}

\theoremstyle{plain}

\setlist[itemize]{leftmargin=3.45ex}
\setlist[itemize,1]{label=$-$,itemsep=0ex,topsep=0.9ex}
\setlist[itemize,2]{label=$\cdot$,topsep=0.5ex,leftmargin=2.75ex}
\setlist[enumerate]{leftmargin=3ex,itemsep=0.1ex,parsep=1ex,topsep=0.9ex}

\definecolor{tabcolor}{HTML}{6666AA}
\definecolor{f1}{HTML}{000060}
\definecolor{f2}{HTML}{0000FF}
\definecolor{f3}{HTML}{36648B}
\definecolor{f4}{HTML}{4682B4}
\definecolor{f5}{HTML}{5CACEE}
\definecolor{f6}{HTML}{00FFFF}
\definecolor{f7}{HTML}{00DD99}
\definecolor{f8}{HTML}{008888}
\definecolor{f9}{HTML}{000000}

\usepackage{float}
\newfloat{program}{thp}{lop}
\floatname{program}{Program}
\crefname{program}{program}{programs}

\newcommand{\inputExpensiveFigure}[1]{
\ifthenelse{\expensiveFigures = 1}{\input{#1}}{}
}

\usepackage{titling}
\pretitle{\begin{center}\linespread{1.05}\Large}
\posttitle{\par\end{center}\vskip 0.5em}
\preauthor{\begin{center}\linespread{1.15}\large \lineskip 0.6em\begin{tabular}[t]{c}}
\postauthor{\end{tabular}\par\end{center}\vskip 0.6em}
\predate{\begin{center}}
\postdate{\par\end{center}}
\usepackage[]{titlesec}
\usepackage{etoolbox}
\makeatletter
\patchcmd{\ttlh@hang}{\parindent\z@}{\parindent\z@\leavevmode}{}{}
\patchcmd{\ttlh@hang}{\noindent}{}{}{}
\makeatother
\titleformat{\section}
{\normalfont\large\bfseries}{\thesection}{0.9ex}{}
\titleformat{\subsection}
{\normalfont\normalsize\bfseries}{\thesubsection}{0.9ex}{}
\titleformat{\subsubsection}
{\normalfont\normalsize\upshape}{\thesubsubsection}{0.9ex}{}
\titleformat{\paragraph}[runin]
{\normalfont\normalsize\itshape}{\theparagraph}{1em}{}
\titleformat{\subparagraph}[runin]
{\normalfont\normalsize\itshape}{\theparagraph}{1em}{}
\titlespacing*{\section}     {0pt}{21dd plus 8pt minus 4pt}{10.5dd}
\titlespacing*{\subsection}   {0pt}{21dd plus 8pt minus 4pt}{10.5dd}
\titlespacing*{\subsubsection}{0pt}{19dd plus 8pt minus 4pt}{10.5dd}
\titlespacing*{\paragraph}   {0pt}{13pt plus 8pt minus 4pt}{1em}
\titlespacing*{\subparagraph}   {0pt}{13pt plus 8pt minus 4pt}{1em}

\usepackage[page]{appendix} 

\DeclareRobustCommand{\VAN}[3]{#2} 

\newcommand{\myorcidlink}[1]{\,\href{https://orcid.org/#1}{\raisebox{-0.45ex}{\includegraphics[width=1.8ex]{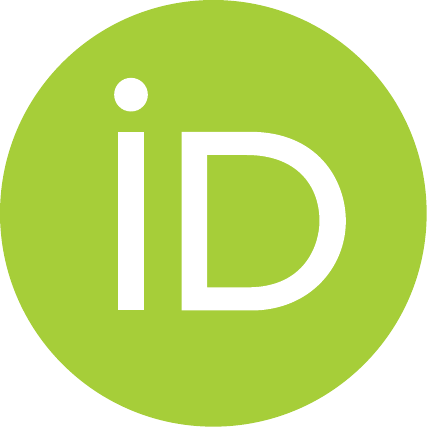}}}}

\newcommand{\scipoptauthors}{%
  Christopher Hojny\protect\myorcidlink{0000-0002-5324-8996},
  Mathieu Besançon\protect\myorcidlink{0000-0002-6284-3033},
  Ksenia Bestuzheva\protect\myorcidlink{0000-0002-7018-7099},
  Sander Borst\protect\myorcidlink{0000-0003-4001-6675},
  Antonia~Chmiela\protect\myorcidlink{0000-0002-4809-2958},
  Jo{\~a}o~Dion{\'i}sio\protect\myorcidlink{0009-0005-5160-0203},
  Johannes Ehls\protect\myorcidlink{0009-0005-1130-6683},
  Leon Eifler\protect{\myorcidlink{0000-0003-0245-9344}},
  Mohammed Ghannam\protect\myorcidlink{0000-0001-9422-7916},
  Ambros Gleixner\protect\myorcidlink{0000-0003-0391-5903},
  Adrian~G{\"o}{\ss}\protect\myorcidlink{0009-0002-7144-8657},
  Alexander Hoen\protect\myorcidlink{0000-0003-1065-1651},
  Jacob von Holly-Ponientzietz\protect\myorcidlink{0009-0002-2601-3689},
  Rolf van der Hulst\protect\myorcidlink{0000-0002-5941-3016},
  Dominik Kamp\protect\myorcidlink{0009-0005-5577-9992},
  Thorsten Koch\myorcidlink{0000-0002-1967-0077},
  Kevin Kofler,
  Jurgen Lentz\protect\myorcidlink{0009-0000-0531-412X},
  Marco~L\"ubbecke\protect\myorcidlink{0000-0002-2635-0522},
  Stephen J. Maher\protect\myorcidlink{0000-0003-3773-6882},
  Paul Matti Meinhold\protect\myorcidlink{0009-0003-5477-9152},
  Gioni Mexi \protect\myorcidlink{0000-0003-0964-9802},
  Til Mohr\protect\myorcidlink{0009-0001-9842-210X},
  Erik M\"uhmer\protect\myorcidlink{0000-0003-1114-3800},
  Krunal Kishor Patel\protect\myorcidlink{0000-0001-7414-5040},
  Marc E. Pfetsch\protect\myorcidlink{0000-0002-0947-7193},
  Sebastian Pokutta\protect\myorcidlink{0000-0001-7365-3000},
  Chantal Reinartz Groba\protect\myorcidlink{0009-0001-1820-3864},
  Felipe Serrano\protect\myorcidlink{0000-0002-7892-3951},
  Yuji Shinano\protect\myorcidlink{0000-0002-2902-882X},
  Mark Turner\protect\myorcidlink{0000-0001-7270-1496},
  Stefan Vigerske\protect\myorcidlink{0009-0001-2262-0601},
  Matthias Walter\protect\myorcidlink{0000-0002-6615-5983},
  Dieter Weninger\protect\myorcidlink{0000-0002-1333-8591},
  Liding Xu\protect\myorcidlink{0000-0002-0286-1109}
}

\ifthenelse{\zibreport = 1}{
  \let\pdfoutorg\pdfoutput
  \let\pdfoutput\undefined
  \usepackage{zibtitlepage}
  \let\pdfoutput\pdfoutorg

  \ZTPTitle{\longtitle}
  \ZTPAuthor{\scipoptauthors}
  \ZTPPreprint
  \ZTPNumber{21-41}
  \ZTPMonth{12}
  \ZTPYear{2021}
  \ZTPInfo{}
}{}

\newcommand{\myand}{$\cdot$\xspace}

\begin{document}

\ifthenelse{\guidelines = 1}{
\section*{Guidelines}

\subsection*{Managing Todos}

Please \textit{\bfseries always} use the todonotes package in order to mark todos.

\begin{itemize}
\item If you are adding a todo, prepend your name abbreviation, e.g.
  \todo[inline]{AG: don't forget about this}
\item If you are adding a todo or comment for someone, address using the ``@'' symbol, e.g.
  \todo[inline]{AG@GG: note that CHANGELOG says ``conflict'' graph}
\item \textit{\bfseries Do not delete todos yourself.}  The person who put it there should delete it when they think it is done, in exceptional cases also the lead author.  If you fixed a todo, simply add your comments, e.g.
  \todo[inline]{AG@SV: Please adjust to central notation above.  SV@AG: OK like this?}
  You can also use this to suggest a different solution or argue against the todo.
\end{itemize}

\subsection*{Editing Comments}

\begin{itemize}
\item \textit{\bfseries Be mindful of sentence length.}  For instance, If you give an example
  or an additional clarification, it can be useful to do so in a separate
  sentence.
\item Do not use personal pronouns. Write in the third person.\todo[inline]{GG: In some cases, I prefer ``we'', e.g., if it is really something that \textbf{we} observed in our experiments. But I totally agree that something implemented in the code should not be referenced by ``we''. AG@GG: I agree that there are some (few) cases where ``we'' is OK.  The one you mention is OK, also anything where we as development team made an individual choice (e.g., defining benchmark set) .  When encountered during proofreading, I may leave these ``we''s as they are or rewrite them.  Still, avoiding it where possible helps to solve the challenge of creating a unified style and sounds slightly more scientific, which is good because it counteracts the software-oriented perception of the report.}
\item Avoid the use of abbreviations, such as w.l.o.g, e.g., and i.e. Spell out what you are saying.
\item Use acronyms for regularly repeated terms.
\item The articles for common acronyms are: an LP, an MILP and an MINLP.
\item ``Code'' can be used in singular only.
\item Use author names when referencing.
\item Use parentheses sparingly. If something is important, integrate it properly into the sentence structure. If not, leave it out.
\end{itemize}

\subsection*{General Structure of a Feature Description}

\begin{itemize}
\item Length: the standard length should be half a page to a full page of text, not counting
  graphics and tables; some technical improvements can be shorter, some extensive features
  longer. \textit{\bfseries As compact as possible, as long as necessary.}
\item General structure:
  \begin{enumerate}
  \item Motivation and background: What is the problem, the issue that the
    features tries to address?  \textit{\bfseries Add references to the
      literature, not only to our own papers.}
  \item Description of the feature: Try to keep this as general (read:
    as non-\scip-specific as possible). \textit{\bfseries Where possible, use MILP/MINLP notation
      from Section~\ref{sect:introduction}} with the predefined macros.
  \item Implementation details: \scip specifics, involved plugins, performance
    critical implementation tricks, most important parameters, default
    parameters (activated?).  However: Avoid mentioning details about our
    development process like merge requests.
  \item Computational results (where applicable)
    Often it will suffice to separate the four parts only by paragraph, without additional section headers or \verb|\paragraph{}|.
  \end{enumerate}
\item Wherever possible refer to existing research articles for further details.
\end{itemize}
All these rules can be broken if there is a good reason, but generally they help
to create a unified form.

\subsection*{Reporting Performance Results}

\begin{itemize}
\item For MILP and MINLP, prefer reporting performance results on mipdev-solvable and
  minlpdev-solvable if it makes sense.
\item Use the same time limits as in the comparison to the last release.
\item You can report intermediate results, e.g., ``At time of activating this
  feature, it gave a performance improvement of \ldots''
\item Be careful when using the word \textit{\bfseries significant}, e.g., in
  ``no significant performance improvement''.  Use it exclusively in the strict
  statistical sense.  If you mean ``large'', then use words like considerable,
  noticable, large, ... Generally: be concrete and use numbers when reporting
  performance results wherever possible.
\item Be careful when reporting results on additional instances solved: ``X more
  instances solved'' is unclear on whether only the number increased or
  additional instances were solved, none were lost.  So either ``The number of
  solved instances increased by X'' or ``X instances were solved in addition,
  while only Y instances could not be solved anymore within the time limit.''
  \textit{\bfseries Only report number of solved instances if it is large and/or
    robust (over several seeds/permutations or clearly dominating)}, because
  small variations here may easily come from performance variability.
  \textit{\bfseries If not clear from a table, mention size of the testset.}
\item If you are taking the numbers from runs available on Rubberband,
  \textit{\bfseries add Rubberband link(s) in a todo environment}
  \verb|\todo[caption={}]{\url{...}}| such that proof-readers can make more
  informed comments.
\item Use tables only where necessary and useful.
\end{itemize}

\subsection*{\LaTeX Guidelines}

\begin{itemize}
\item Always check whether there is already a \textit{\bfseries macro} for a mathematical
  expression, a solver, testset, etc. name and use it; introduce new ones as
  necessary.  Examples: \verb|\R|, \verb|\miplib|, \verb|\coral|, \ldots
\item Do not use long lines in the source, but break them at around 80-120
  chars.  Otherwise the risk of conflicts increase and git diff does not show
  small changes clearly.
\item No spaces after \verb|$|, not \verb|$ a + b $|, but \verb|$a + b$|.
\end{itemize}

\subsection*{Unified Style and Notation}

\begin{itemize}
\item Use American English.
\item Use title case where applicable (section and paragraph headers, lemma and definition names).
\item No hyphen for ``LP relaxation'', ``MILP solver'', ``LP solver'', \ldots.
\item Hyphen in adjectives ``mixed-integer'', ``large-scale'', ``well-known'', ``problem-specific'', ``right-hand'', ``left-hand'', ``general-purpose'', ``user-defined'', \ldots.
\item No hyphen in ``nonzero'', ``nontrivial'', ``nonconvex''.
\item No hyphen in ``plugin''.
\item No space in ``testset''.
\item Use ``zeros'' and ``nonzeros''.
\item No period inside a \verb|\paragraph{}|.
\item All algorithms should use the \texttt{algorithm} environment.
\item The assignment operator in algorithm should be $i \leftarrow j$.
\item Use $\{i : x_{i} > 0\}$ for sets.
\item For single line problems $\{c^\T x : Ax = b, x \in \Z^{n}\}$. For multiple line problems
  \begin{equation}
    \begin{aligned}
      \min \quad& \linobj^\T x \\
      \text{s.t.} \quad& \linmatrix x \geq \rhs, \\
      &\lb_{i} \leq x_{i} \leq \ub_{i} && \fa i \in \varindex, \\
      &x_{i} \in \Z && \fa i \in \intvarindex,
    \end{aligned}
  \end{equation}
\item Define variable sets as $x \in \Z^{n}$, if mixed-integer then $(x,y)
  \in \Z^{p} \times \R^{n - p}$.
\item Use $\T$ for transpose (\verb|\T|).
\item Use $\R_+$, instead of $\R_{\geq 0}$.
\item Use $\bar{x}$ (\verb|\bar{x}|) instead of $\overline{x}$.
\item Use $\R$, $\Q$, $\N$, $\Z$ for the reals, rationals, naturals and integers.
\item The empty set should be given by $\emptyset$ (\verb|\emptyset|).
\item Equation environments:
  \begin{itemize}
  \item Single line equations: equation
  \item Multiple line equations with each line numbered: align
  \item Multiple line equations with a single number: equation wrapping aligned
  \end{itemize}
\item Problem definitions
  \begin{itemize}
  \item All problems should be numbered and referenced by their number.
  \item All multiple line problems should be given a single number, unless the individual lines need to be referenced.
  \end{itemize}
\item Captions are full sentences and end with a period.
\end{itemize}

\subsection*{References}

\begin{itemize}
\item only DOI if available, no additional URL, URN, ISBN, eprint, ...
\item The DOIs \emph{do not} include the web address
  \texttt{http://dx.doi.org/}. Only supply DOIs no additional web address
  to the DOIs -- this is redundant.
\item preprints with URL
\item page numbers with \verb|--| (not \verb|-|)
\item check capitalization of titles
\item ZIB-Reports should not be cited using the bibtex information that
  comes out of the system. It should rather be:
  \begin{itemize}
  \item The type of the bibtex entry should be \texttt{techreport}.
  \item The \texttt{institution} should be \texttt{Zuse Institute Berlin} (the
    abbreviation ZIB is not known to everyone).
  \item The \texttt{type} should be \texttt{ZIB-Report}.
  \item The number of the report appears in \texttt{numbeR} (not \texttt{volume}).
  \item The \texttt{address} and \texttt{URN} should not be given. It is
    not important for citing.
  \end{itemize}
  Here is an example:
\begin{verbatim}
@techreport{Gleixner2012,
  author =	 {Ambros M. Gleixner},
  title =	 {Factorization and update of a reduced basis matrix for the revised simplex method},
  institution =	 {Zuse Institute Berlin},
  type =	 {ZIB-Report},
  number =	 {12-36},
  year =	 {2012}
}
\end{verbatim}
\end{itemize}

Please extend the list if you notice any inconsistencies during proofreading.
}{}

\title{\longtitle}

\author{%
  Christopher Hojny\protect\myorcidlink{0000-0002-5324-8996} \myand
  Mathieu Besançon\protect\myorcidlink{0000-0002-6284-3033} \myand
  Ksenia Bestuzheva\protect\myorcidlink{0000-0002-7018-7099} \\
  Sander Borst\protect\myorcidlink{0000-0003-4001-6675} \myand
  Antonia~Chmiela\protect\myorcidlink{0000-0002-4809-2958},
  Jo{\~a}o~Dion{\'i}sio\protect\myorcidlink{0009-0005-5160-0203} \myand
  Johannes Ehls\protect\myorcidlink{0009-0005-1130-6683} \\
  Leon Eifler\protect{\myorcidlink{0000-0003-0245-9344}} \myand
  Mohammed Ghannam\protect\myorcidlink{0000-0001-9422-7916} \myand
  Ambros Gleixner\protect\myorcidlink{0000-0003-0391-5903} \\
  Adrian~G{\"o}{\ss}\protect\myorcidlink{0009-0002-7144-8657} \myand
  Alexander Hoen\protect\myorcidlink{0000-0003-1065-1651} \myand
  Jacob von Holly-Ponientzietz\protect\myorcidlink{0009-0002-2601-3689} \\
  Rolf van der Hulst\protect\myorcidlink{0000-0002-5941-3016} \myand
  Dominik Kamp\protect\myorcidlink{0009-0005-5577-9992} \myand
  Thorsten Koch\myorcidlink{0000-0002-1967-0077} \\
  Kevin Kofler \myand
  Jurgen Lentz\protect\myorcidlink{0009-0000-0531-412X} \myand
  Marco~L\"ubbecke\protect\myorcidlink{0000-0002-2635-0522} \\
  Stephen J. Maher\protect\myorcidlink{0000-0003-3773-6882} \myand
  Paul Matti Meinhold\protect\myorcidlink{0009-0003-5477-9152} \myand
  Gioni Mexi \protect\myorcidlink{0000-0003-0964-9802} \\
  Til Mohr\protect\myorcidlink{0009-0001-9842-210X} \myand
  Erik M\"uhmer\protect\myorcidlink{0000-0003-1114-3800} \myand
  Krunal Kishor Patel\protect\myorcidlink{0000-0001-7414-5040} \\
  Marc E. Pfetsch\protect\myorcidlink{0000-0002-0947-7193} \myand
  Sebastian Pokutta\protect\myorcidlink{0000-0001-7365-3000} \myand
  Chantal Reinartz Groba\protect\myorcidlink{0009-0001-1820-3864} \\
  Felipe Serrano\protect\myorcidlink{0000-0002-7892-3951} \myand
  Yuji Shinano\protect\myorcidlink{0000-0002-2902-882X} \myand
  Mark Turner\protect\myorcidlink{0000-0001-7270-1496} \myand
  Stefan Vigerske\protect\myorcidlink{0009-0001-2262-0601} \\
  Matthias Walter\protect\myorcidlink{0000-0002-6615-5983} \myand
  Dieter Weninger\protect\myorcidlink{0000-0002-1333-8591} \myand
  Liding Xu\protect\myorcidlink{0000-0002-0286-1109}
  \thanks{Extended author information is available at the end of the paper.
    \shortfunding}}


\ifthenelse{\zibreport = 1}{\zibtitlepage}{}

\newgeometry{left=38mm,right=38mm,top=35mm}

\maketitle

\paragraph{\bf Abstract}

The \scipopt provides a collection of software packages for mathematical optimization, centered around the constraint integer programming (CIP) framework \scip.
This report discusses the enhancements and extensions included in \scipoptv.
The updates in \scipv include a new solving mode for exactly solving rational mixed-integer linear programs, a new presolver for detecting implied integral variables, a novel cut-based conflict analysis and separator for flower inequalities, two new heuristics, a novel tool for explaining infeasibility, a new interface for nonlinear solvers as well as improvements in symmetry handling, branching strategies, and \scip's Benders' decomposition framework.
\scipoptv also includes new and improved features in the the presolving library \papilo, the parallel framework \ug, and the decomposition framework \gcg.
Moreover, the \scipoptv contains MIP-DD, the first open-source delta debugger for mixed-integer programming solvers.
These additions and enhancements have resulted in an overall performance improvement of \scip in terms of solving time, number of nodes in the branch-and-bound tree, as well as the reliability of the solver.

\paragraph{\bf Keywords} Constraint integer programming
$\cdot$ linear programming
$\cdot$ mixed-integer linear programming
$\cdot$ mixed-integer nonlinear programming
$\cdot$ optimization solver
$\cdot$ branch-and-cut
$\cdot$ branch-and-price
$\cdot$ column generation
$\cdot$ parallelization
$\cdot$ mixed-integer semidefinite programming

\paragraph{\bf Mathematics Subject Classification} 90C05 $\cdot$ 90C10 $\cdot$ 90C11 $\cdot$ 90C30 $\cdot$ 90C90 $\cdot$ 65Y05

\newpage


\section{Introduction}
\label{sect:introduction}


The \scipopt comprises a set of complementary software packages designed to model and
solve a large variety of mathematical optimization problems:
\begin{itemize}
\item the constraint integer programming solver
\scip~\cite{Achterberg2009}, a solver for
mixed-integer linear and nonlinear programs as well as a flexible framework for
branch-cut-and-price,
\item the simplex-based linear programming solver
  \soplex~\cite{Wunderling1996},
\item the modeling language \zimpl~\cite{Koch2004},
\item the presolving library \presollib for linear and mixed-integer linear programs,
\item the automatic decomposition solver \gcg~\cite{GamrathLuebbecke2010}, and
\item the \ug framework for parallelization of branch-and-bound
  solvers~\cite{Shinano2018}.
\end{itemize}
All six tools are freely available as open-source software packages; \scip~10.0,
\soplex~8.0, \presollib~3.0, and \gcg~4.0 are licensed under the Apache~2.0 license, whereas
\zimpl~3.7.0 and \ug~1.0 make use of the GNU Lesser General Public License.
A notable extension to the \scipopt is the mixed-integer semidefinite
programming solver \scipsdp~\cite{GallyPfetschUlbrich2018}.
Moreover, the delta debugging tool MIP-DD~\cite{HKG24MIPDD} can be used to
assist the development and debugging of optimization software like \scip.%

\paragraph{New Developments and Structure of the Paper}
The goal of this report is to highlight the features becoming available in
the latest release of the \scipopt, where the focus is on the developments
of \scip.
After providing a short background on the problems that can be handled by
\scip, Section~\ref{sect:performance} investigates the change of performance between
\scippv and \scipv for mixed-integer linear and nonlinear programs.
Briefly summarized, \scipv is more performant than \scippv both in terms of
running time and number of solved instances, where the performance
improvement is more pronounced for nonlinear problems.
Section~\ref{sect:scip} forms the main part of this report and highlights the new
features and technical improvements becoming available in \scip~10, namely
\begin{itemize}
\item a new solving mode that allows to solve rational \MILPs exactly,
  i.e., without numerical tolerances and without being affected by
  floating-point roundoff errors (Section~\ref{sec:exact});
\item a new presolver for detecting implied integral variables and a more
  flexible specification of implied integrality information (Section~\ref{sec:presolve});
\item generalizations of symmetry handling methods to reflection symmetries (Section~\ref{sec:symmetry});
\item cut-based conflict analysis (Section~\ref{sec:conflict});
\item a new separator plugin for flower inequalities derived from multilinear
  problems (Section~\ref{sec:separators});
\item two new primal heuristics that exploit decompositions provided by users (Section~\ref{sec:heuristics});
\item improved branching strategies (Section~\ref{sec:branching});
\item enhancing \scip's Benders' decomposition framework by a more flexible
  way to formulate problems and an improved detection of master linking
  variables (Section~\ref{sec:benders});
\item a tool to find explanations for infeasibility via so-called
  irreducible infeasible subsystems (Section~\ref{sec:iis});
\item an interface to the nonlinear programming solver CONOPT (Section~\ref{sec:conopt});
\item two technical improvements: a new constraint handler to avoid infeasibilities due to aggregations from presolving and a JSON export of \scip's statistics (Section~\ref{sec:misc}).
\end{itemize}
Section~\ref{sect:gcg} is devoted to improvements of \gcg, including new data
structures and file formats for storing decompositions, newly
developed solvers for pricing problems, parallelization of pricing,
and support for handling branching and cutting decisions formulated
directly in the extended formulation variables. Section~\ref{sect:scip-sdp}
gives a short update on \scipsdp.
Section~\ref{sect:papilo} highlights the latest developments of \presollib improving
performance and memory management.
Section~\ref{sect:soplex} gives a very brief summary of the capabilities of
\soplex,
and Sections~\ref{sect:ug} and~\ref{sect:zimpl} briefly summarize the latest
changes in \ug and \zimpl, respectively.
In Section~\ref{sect:interfaces}, new features of \scip's various interfaces are
discussed, and Section~\ref{sect:mipdd} concerns the delta debugger MIP-DD.
Section~\ref{sect:applications} discusses the newly developed \solver{PBSolver}, a
\scip-based solver that is tailored for Pseudo-Boolean problems and won
several categories of the Pseudo-Boolean Competition~24.

\paragraph{Background}

\scip is designed as a solver for \emph{constraint integer programs} (\CIPs), a generalization of mixed-integer linear and nonlinear programs (\MILPs and \MINLPs). \CIPs are finite-dimensional optimization problems with arbitrary constraints and a linear objective function that satisfy the following property: if all integer variables are fixed, the remaining subproblem must form a linear or nonlinear program (\LP or \NLP).
To solve \CIPs, \scip constructs relaxations---typically linear relaxations, but also nonlinear relaxations are possible, or relaxations based on semidefinite programming for \scipsdp. If the relaxation solution is not feasible for the current subproblem, an \emph{enforcement} procedure is called that resolves the infeasibility, for example by branching or by separating cutting planes.

The most important subclass of \CIPs that are solvable with \scip are
\emph{mixed-integer programs} (\MIPs) which can be purely linear (\MILPs)
or contain nonlinear terms in the constraints or the objective
function (\MINLPs).
\MILPs are optimization problems of the form
\begin{equation}
  \begin{aligned}
    \min \quad& \linobj^\T x \\
    \text{s.t.} \quad& \linmatrix x \geq \rhs, \\
    &\lb_{i} \leq x_{i} \leq \ub_{i} && \fa i \in \varindex, \\
    &x_{i} \in \Z && \fa i \in \intvarindex,
  \end{aligned}
  \label{eq:generalmilp}
\end{equation}
defined by $c \in \R^n$, $A \in\R^{m\times n}$, $ \rhs\in \R^{m}$, $\lb$, $\ub \in
\Rinf^{n}$, and the index set of integer variables $\mathcal{I} \subseteq \mathcal{N} \defi \{1, \ldots, n\}$.  The usage of $\Rinf \defi \R \cup
\{-\infty,\infty\}$ allows for variables that are free or bounded only in
one direction (we assume that variables are not fixed to~$\pm \infty$).
In contrast, \MINLPs are optimization problems of the form
\begin{equation}
  \begin{aligned}
    \min \quad& \nonlinobj(x) \\
    \text{s.t.} \quad& \nonlincons_{k}(x) \leq 0 && \fa k \in \consindex, \\
    &\lb_{i} \leq x_{i} \leq \ub_{i} && \fa i \in \varindex, \\
    &x_{i} \in \Z && \fa i \in \intvarindex,
  \end{aligned}
  \label{eq:generalminlp}
\end{equation}
where the functions $\nonlinobj\colon \R^n \rightarrow \R$ and $\nonlincons_{k}\colon \R^{n}
\rightarrow \R$, $k \in \consindex \defi \{1,\ldots,m\}$, are possibly nonconvex.
Within \scip, it is assumed that $\nonlinobj$ is linear and that
$\nonlincons_k$, $k \in \consindex$, is specified explicitly in
algebraic form using a known set of base expressions.
Due to its design as a solver for \CIPs, \scip can be extended by plugins for more general or problem-specific classes of optimization problems.
The core of \scip is formed by a central branch-cut-and-price algorithm that utilizes an \LP as the default relaxation which can be solved by a number of different \LP solvers, controlled through a uniform \emph{\LP interface}. To be able to handle any type of constraint, a \emph{constraint handler} interface is provided.
This interface allows to integrate new constraint types, and provides support for many different well-known types of constraints out of the box.
Further solving methods like primal heuristics, branching rules, and cutting plane separators can also be integrated as plugins with a pre-defined interface.
\scip comes with many such plugins needed to achieve a good \MILP and \MINLP performance.
In addition to plugins supplied as part of the SCIP distribution, new plugins can be created by users.
The design approach and solving process is described in detail by
Achterberg~\cite{Achterberg2007a}.

Although \scip is a standalone solver,
it interacts closely with the other components of the \scipopt.
\zimpl is integrated into \scip as a reader plugin, making it possible to read \zimpl problem instances directly by \scip.
\presollib is integrated into \scip as an additional presolver plugin.
The \LPs that need to be solved as relaxations in the branch-and-bound process are by default solved with \soplex,
and it is possible to replace \soplex by other LP solvers.
Interfaces to most actively developed external \LP solvers exist, and new interfaces can be added by users.
\gcg extends \scip to automatically detect problem structure
and generically apply decomposition algorithms based on the Dantzig-Wolfe or the
Benders' decomposition schemes.
Finally, the default instantiations of the \ug framework use \scip as a base
solver in order to perform branch-and-bound in parallel computing
environments with shared or distributed memory architectures.


\section{Overall Performance Improvements for MILP and MINLP with SCIP}
\label{sect:performance}

A major use of the \scipopt is as an out-of-the-box solver
for mixed-integer linear and nonlinear programs.
Therefore, the performance of \scip on \MILP and \MINLP instances is of particular
interest during the development process. In this section, we present
computational experiments to assess the performance of \scip~10.0 in comparison with the previous major release, \scip~9.0, and the latest bugfix release, \scip~9.2.4.
The reason for comparing \scip~10.0 with these two previous versions is to highlight the effect of both increasing robustness of \scip by fixing bugs (comparison with \scip~9.2.4) and new features (comparison with \scip~9.0).


\subsection{Experimental Setup}

We compare \scip~10.0, including \soplex~8.0.0 and \papilo~3.0.0, with \scip~9.0 (including \soplex~7.0.0 and \papilo~2.2.0) as well as \scip~9.2.4 (including \soplex~7.1.6 and \papilo~2.4.4).
All versions of \scip were built with GCC 10.2.1.
Further software packages linked to \scip include the NLP solver 
\solver{Ipopt}~3.14.18 built with \solver{HSL~MA27} as linear system solver,
\solver{Intel MKL} as linear algebra package, and the algorithmic differentiation code
\solver{CppAD}~20180000.0.
For symmetry detection, \scip~10.0 uses the external software \solver{nauty}~2.8.8 combined with the symmetry preprocessor \sassy~2.0;
\scip~9.2.4 uses \solver{nauty}~2.8.8 with \sassy~1.1, and \scip~9.0 uses \solver{bliss}~0.77 with \sassy~1.1.
The time limit was set to~\SI{5400}{\second} for \MILP and to~\SI{3600}{\second} for the
\MINLP runs.
Note that for \MINLP, an instance is considered solved when either a relative
primal-dual gap of $10^{-4}$, or an absolute gap of $10^{-6}$ is reached.

The \MILP instances are selected from the \miplib~2003, 2010, and~2017~\cite{MIPLIB} sets, as well as the
\coral instance set, and include all instances that could be solved by previous
releases, for a total of~349 instances.
The \MINLP instances are selected similarly from the \minlplibtwo\footnote{\url{https://www.minlplib.org}} library,
for a total of~169 instances.
All performance runs are carried out on identical machines with Intel Xeon Gold~5122~@~\SI{3.60}{\giga\hertz}
and~\SI{96}{\giga\byte} RAM. A single run is carried out on each machine in single-threaded mode.
Each \MILP instance is solved with \scip using three different seeds for random number generators, while each \MINLP instance is solved using five different seeds.
This results in a testset of~1047 \MILPs and~845 \MINLPs.

The indicators of interest for our comparison 
on a given subset of instances are
the number of solved instances, the shifted geometric mean
of the number of branch-and-bound nodes, and the shifted geometric mean
of the solving time.
The \emph{shifted geometric mean} of values $t_1, \dots, t_n$ is
\[
\left(\prod_{i=1}^n(t_i + s)\right)^{1/n} - s.
\]
The shift~$s$ is set to 100~nodes and 1~second, respectively.

In Tables~\ref{tab:milp_performance} and~\ref{tab:minlp_performance},
we present the results of the computational experiments for \MILP and \MINLP comparing \scip~10.0 and \scip~9.0;
Tables~\ref{tab:milp_performance_2} and~\ref{tab:minlp_performance_2} show the results of the comparison of \scip~10.0 and \scip~9.2.4.
We present these statistics for several subsets of instances.
The subset \cleaninst contains all instances of the testset 
excluding those with numerically inconsistent results.
The subset ``affected'' contains all instances where solvers differ in the
number of dual simplex iterations.
The brackets $[t,T]$ collect the subsets of instances which were solved by
at least one solver and for which the maximum solving time (among both solver
versions) is at least $t$~seconds and at most~$T$ seconds, where $T$ is
usually equal to the time limit.
With increasing~$t$, this provides a hierarchy of subsets of increasing
difficulty.
The subset ``\alloptimal'' contains instances that can
be solved by both versions within the time limit.
Columns ``instances'' and ``solved'' refer to the number of instances in a subset of instances and the number of solved instances, respectively.
Moreover, columns ``time'' and ``nodes'' provide the shifted geometric mean of the running time and number of nodes, respectively, needed by the different versions.

\subsection{\MILP Performance}

Table \ref{tab:milp_performance} summarizes the results comparing the performance of
\scipv against \scippv on the \MILP testset.
Using \scipv instead of \scippv slightly increases the number of solved instances from~886 to~888 and decreases the mean running time by approximately~\SI{2}{\percent}.
Despite the reduced running time, the number of nodes in \scip's branch-and-bound tree increases by~\SI{3}{\percent}.
Since the mean number of nodes is relatively small though, the increased size of the branch-and-bound tree is only moderate in absolute numbers.

On the various subsets of instances, one can observe that \scip's performance improvement is better for the harder instances that require at least ten seconds to be solved by both versions.
Among all test sets, including both-solved, \scipv is at least~\SI{3}{\percent} faster than \scippv; for the~\bracket{100}{timelim} subset even an~\SI{8}{\percent} performance improvement can be observed and also the number of nodes reduces by~\SI{3}{\percent}.
That is, \scip's performance gain on the harder instances mainly contributes to its overall performance improvement in comparison with \scippv.

Comparing \scipv with \scip~9.2.4 (Table~\ref{tab:milp_performance_2}), similar observations as before can be made, however, with more pronounced performance improvements.
The overall performance of \scipv is \SI{4}{\percent} better than the performance of \scip~9.2.4, where again the biggest performance improvements are achieved for the harder instances.
For example, for the \bracket{100}{timelim} and \bracket{1000}{timelim} subsets, \scipv has a \SI{9}{\percent} and \SI{11}{\percent}, respectively, lower running time than \scip~9.2.4.
This also results in a reduced number of nodes in the branch-and-bound trees on most subsets.

To summarize, \scipv has become more reliable than \scippv based on the bugfixes that have been introduced in \scip~9.2.4.
While some of these bugfixes caused a slowdown of \scip~9.2.4 in comparison with \scippv, the additional features of \scipv yield an overall performance improvement both in comparison with \scippv and \scip~9.2.4.

\begin{table}[t]
\caption{Performance of \scipv and \scippv for \MILP instances.}\label{tab:milp_performance}
\scriptsize

\begin{tabular*}{\textwidth}{@{}l@{\;\;\extracolsep{\fill}}rrrrrrrrr@{}}
\toprule
&           & \multicolumn{3}{c}{\scip~10.0.0+\soplex~8.0.0} & \multicolumn{3}{c}{\scip~9.0.0+\soplex~7.0.0} & \multicolumn{2}{c}{relative} \\
\cmidrule{3-5} \cmidrule{6-8} \cmidrule{9-10}
Subset                & instances &                                   solved &       time &        nodes &                                   solved &       time &        nodes &       time &        nodes \\
\midrule
\cleaninst            &      1046 &                                      888 &      232.6 &         2659 &                                      886 &      236.3 &         2614 &        0.98 &          1.02 \\
\affected             &       827 &                                      796 &      195.6 &         2823 &                                      794 &      198.8 &         2733 &        0.98 &          1.03 \\
\cmidrule{1-10}
\bracket{1}{timelim}    &       890 &                                      859 &      174.7 &         2293 &                                      857 &      177.9 &         2224 &        0.98 &          1.03 \\
\bracket{10}{timelim}   &       815 &                                      784 &      242.1 &         2912 &                                      782 &      248.6 &         2825 &        0.97 &          1.03 \\
\bracket{100}{timelim}  &       573 &                                      542 &      569.1 &         6337 &                                      540 &      617.3 &         6528 &        0.92 &          0.97 \\
\bracket{1000}{timelim} &       228 &                                      197 &     1919.5 &        19663 &                                      195 &     1981.4 &        19049 &        0.97 &          1.03 \\
\alloptimal           &       855 &                                      855 &      119.3 &         1695 &                                      855 &      122.4 &         1671 &        0.97 &          1.01 \\
\bottomrule
\end{tabular*}
\end{table}


\begin{table}[t]
\caption{Performance of SCIP 10.0 and SCIP 9.2.4 for MILP instances.}\label{tab:milp_performance_2}
\scriptsize

\begin{tabular*}{\textwidth}{@{}l@{\;\;\extracolsep{\fill}}rrrrrrrrr@{}}
\toprule
&           & \multicolumn{3}{c}{\scip~10.0.0+\soplex~8.0.0} & \multicolumn{3}{c}{\scip~9.2.4+\soplex~7.1.6} & \multicolumn{2}{c}{relative} \\
\cmidrule{3-5} \cmidrule{6-8} \cmidrule{9-10}
Subset                & instances &                                   solved &       time &        nodes &                                   solved &       time &        nodes &       time &        nodes \\
\midrule
\cleaninst            &      1046 &                                      888 &      232.6 &         2659 &                                      883 &      241.6 &         2669 &        0.96 &          1.00 \\
\affected             &       783 &                                      761 &      202.7 &         3015 &                                      756 &      212.3 &         3043 &        0.96 &          0.99 \\
\cmidrule{1-10}
\bracket{1}{timelim}    &       881 &                                      859 &      168.7 &         2225 &                                      854 &      176.4 &         2243 &        0.96 &          0.99 \\
\bracket{10}{timelim}   &       805 &                                      783 &      234.9 &         2844 &                                      778 &      247.3 &         2862 &        0.95 &          0.99 \\
\bracket{100}{timelim}  &       562 &                                      540 &      561.0 &         6072 &                                      535 &      616.8 &         6413 &        0.91 &          0.95 \\
\bracket{1000}{timelim} &       221 &                                      199 &     1771.0 &        17321 &                                      194 &     2000.8 &        19337 &        0.89 &          0.90 \\
\alloptimal           &       861 &                                      861 &      123.0 &         1770 &                                      861 &      126.4 &         1759 &        0.97 &          1.01 \\
\bottomrule
\end{tabular*}
\end{table}

\subsection{\MINLP Performance}

Table~\ref{tab:minlp_performance} provides an overview of the results for
\scipv and \scippv on the \MINLP testset.
As for the \MILP testset, \scipv slightly increases the number of solved instances from~800 to~802.
For all tested subsets of instances, the number of nodes in the branch-and-bound tree decreases by at least~\SI{2}{\percent}.
On the entire testset, \scipv is \SI{6}{\percent} faster than~\scippv; on the subset of affected instances, the performance improvement is even~\SI{8}{\percent}.

Distinguishing the running time improvements based on the difficulty of instances, the highest performance gain is achieved for instances that require at least~\SI{1000}{\second}, where it materializes in a speed-up of~\SI{20}{\percent}.
This high performance gain can be explained by \scipv solving two more instances than \scippv, i.e., every instance that was solvable by \scippv can also be solved by \scipv within the time limit.
Restricting to the subset both-optimal, Table~\ref{tab:minlp_performance} shows that \scipv is consistently faster on these instances, achieving a performance improvement of~\SI{5}{\percent}.
That is, \scipv both improves the running time and allows to solve more instances.

Comparing \scipv with \scip~9.2.4 on the \MINLP testset (Table~\ref{tab:minlp_performance_2}), we observe a similar trend as for the \MILP testset:
both the performance improvements and reduction of number nodes in the branch-and-bound trees are more pronounced.
On the entire testset, \scipv achieves a performance improvement of \SI{9}{\percent}, and on the \bracket{1000}{timelim} subset the performance even improves by \SI{22}{\percent};
the number of nodes decreases by at least \SI{6}{\percent}.
The additional features of \scipv and bugfixes introduced by \scip~9.2.4 thus also yield for the \MINLP testset a more reliable and faster code base.

\begin{table}[ht]
\caption{Performance of \scip~10.0 and \scip~9.0 for \MINLP instances}
\label{tab:minlp_performance}
\scriptsize

\begin{tabular*}{\textwidth}{@{}l@{\;\;\extracolsep{\fill}}rrrrrrrrr@{}}
\toprule
&           & \multicolumn{3}{c}{\scip~10.0.0+\cplex~12.10.0.0} & \multicolumn{3}{c}{\scip~9.0.0+\cplex~12.10.0.0} & \multicolumn{2}{c}{relative} \\
\cmidrule{3-5} \cmidrule{6-8} \cmidrule{9-10}
Subset                & instances &                                   solved &       time &        nodes &                                   solved &       time &        nodes &       time &        nodes \\
\midrule
\cleaninst            &       809 &                                      802 &       17.9 &         2375 &                                      800 &       19.0 &         2433 &        0.94 &          0.98 \\
\affected             &       696 &                                      696 &       20.1 &         2515 &                                      694 &       21.7 &         2588 &        0.92 &          0.97 \\
\cmidrule{1-10}
\bracket{1}{timelim}    &       720 &                                      720 &       22.9 &         3164 &                                      718 &       24.5 &         3254 &        0.93 &          0.97 \\
\bracket{10}{timelim}   &       444 &                                      444 &       66.7 &         7341 &                                      442 &       74.5 &         7723 &        0.89 &          0.95 \\
\bracket{100}{timelim}  &       168 &                                      168 &      366.6 &        44746 &                                      166 &      413.5 &        48104 &        0.89 &          0.93 \\
\bracket{1000}{timelim} &        52 &                                       52 &     1094.0 &       130587 &                                       50 &     1362.3 &       134476 &        0.80 &          0.97 \\
\alloptimal           &       800 &                                      800 &       17.0 &         2233 &                                      800 &       17.9 &         2274 &        0.95 &          0.98 \\
\bottomrule
\end{tabular*}
\end{table}


\begin{table}
\caption{Performance of SCIP 10.0 and SCIP 9.2.4 for MINLP instances.}
\label{tab:minlp_performance_2}
\scriptsize

\begin{tabular*}{\textwidth}{@{}l@{\;\;\extracolsep{\fill}}rrrrrrrrr@{}}
\toprule
&           & \multicolumn{3}{c}{\scip~10.0.0+\cplex~12.10.0.0} & \multicolumn{3}{c}{\scip~9.2.4+\cplex~12.10.0.0} & \multicolumn{2}{c}{relative} \\
\cmidrule{3-5} \cmidrule{6-8} \cmidrule{9-10}
Subset                & instances &                                   solved &       time &        nodes &                                   solved &       time &        nodes &       time &        nodes \\
\midrule
\cleaninst            &       816 &                                      804 &       18.5 &         2399 &                                      804 &       20.4 &         2547 &        0.91 &          0.94 \\
\affected             &       569 &                                      569 &       21.6 &         2255 &                                      569 &       24.5 &         2460 &        0.88 &          0.92 \\
\cmidrule{1-10}
\bracket{1}{timelim}    &       721 &                                      721 &       22.9 &         3172 &                                      721 &       25.5 &         3384 &        0.90 &          0.94 \\
\bracket{10}{timelim}   &       445 &                                      445 &       66.7 &         7308 &                                      445 &       77.9 &         8080 &        0.86 &          0.90 \\
\bracket{100}{timelim}  &       175 &                                      175 &      330.4 &        38810 &                                      175 &      407.6 &        46111 &        0.81 &          0.84 \\
\bracket{1000}{timelim} &        52 &                                       52 &     1156.7 &       145999 &                                       52 &     1480.6 &       161852 &        0.78 &          0.90 \\
\alloptimal           &       804 &                                      804 &       17.0 &         2246 &                                      804 &       18.8 &         2389 &        0.91 &          0.94 \\
\bottomrule
\end{tabular*}
\end{table}


\section{SCIP}
\label{sect:scip}

\label{sect:exactscip}
\subsection{A Numerically Exact Solving Mode for MILPs}
\label{sec:exact}


As of today, virtually all available MIP solvers rely on  floating-point arithmetic coupled with numerical error tolerances to achieve a tradeoff of numerical accuracy and performance.
While typically a lot of development effort is invested into ensuring high numerical robustness, by design floating-point solvers cannot guarantee that their results are unaffected by accumulated roundoff errors, nor can they produce independently verifiable certificates of correctness.

In \scipv users now for the first time have the option to solve \MILPs with rational input data in a numerically exact solving mode subject to no numerical tolerances.
The exact solving mode relies on the software packages GMP, MPFR, and Boost.
If these dependencies are satisfied while building \scip, the exact solving mode can be activated at runtime simply by setting the parameter \texttt{exact/enabled} to \texttt{true} before the problem instance is loaded.
In addition, \scipv can log a proof of correctness of the LP-based branch-and-bound process by setting the parameter \texttt{certificate/filename}.
\scip then produces a certificate file in \texttt{VIPR}~\cite{VIPR} format that encodes \scip's reasoning for the claimed primal and dual bounds.
The correctness of this certificate can then be verified by independent proof checkers.

As a guiding principle, the exact solving mode follows a hybrid-precision approach and replaces expensive symbolic computations at many places by numerically safe floating-point computations based on directed rounding.
In \scipv, the exact solving mode is restricted to mixed-integer linear programs and supports the following features:
\begin{itemize}
   \item rational presolving using \papilo~\cite{EiflerGleixner2022,GleixnerGottwaldHoen23},
   \item safe dual bounding and reliability pseudocost branching~\cite{CookKochSteffyetal2013,Jarck2020},
   \item safe separation of Gomory mixed-integer cuts~\cite{EiflerGleixner2024},
   \item safe dual proof analysis and constraint propagation~\cite{borst2024certifiedconstraintpropagationdual},
   \item exact post-processing of solutions from all floating-point primal heuristics~\cite{EiflerGleixner2022},
   \item certification of the LP-based branch-and-bound process~\cite{VIPR} except presolving, which can be certified independently for binary programs with integer data~\cite{HoehnOertelGleixnerNordstrom2024}.
\end{itemize}
For a detailed description of the theory and performance of the various exact solving algorithms, we refer to the above references.

The following sections give an overview of the parts of \scip added or modified for implementing the exact solving mode, as well as the general ideas behind the new features.
In order to ensure that addition of an unsafe user plugin does not compromise the correctness of the exact solving mode, all critical plugin types have been equipped with a flag that stores whether a plugin of that type can be safely used in exact solving mode.
As a general design principle, the \scip core executes callbacks only of those plugins that are explicitly marked as safe to use in exact solving mode, which is typically done in the plugin's inclusion method.

\subsubsection{Exact Instance Readers and a Wrapper for Rational Numbers}

The first step in the solving process that can result in numerical errors is the reading of the problem instance. \scipv provides extensions to the readers for MPS, LP, CIP, OPB/WBO, and ZIMPL files to read the problem instance in exact arithmetic. This means that those types of problems can contain any combination of floating-point, rational, and integer coefficients, all of which will be read exactly.
Note that the exact solving mode needs to be enabled prior to the creation of the problem.
The exact readers will then parse all values of bounds and coefficients as rational numbers and add exact variables and constraints to \scip, see Section~\ref{sec:exact:varcons}.
For storing and computing with rational numbers, \scipv comes with a new data structure \texttt{SCIP\_RATIONAL}, which acts as a wrapper to combine functionality of the Boost multiprecision library~\cite{boostmplib}, GMP~\cite{gmplib}, and MPFR~\cite{mpfrpaper,mpfrlib}.

\subsubsection{Exact Variable Data and a New Handler for Exact Linear Constraints}
\label{sec:exact:varcons}

In \scipv, the variable structure contains a new substructure, allocated only in exact solving mode, that keeps rational versions of the objective coefficient, the global and local domains, and data for (multi-)aggregations, the exact LP, and the certificate.
As before, each variable is initially created with floating-point data.
If the floating-point bounds and objective coefficient are correct, the exact data can be initialized with the floating-point data; otherwise, the floating-point data is overwritten.
Note that the floating-point bounds are always maintained as a safe outward rounding of the exact bounds so they can be used algorithmically during constraint propagation and dual bounding.

Following \scip's constraint-based view, also a new constraint handler has been added for exact linear constraints. In the exact solving mode, this constraint handler replaces the standard linear constraint handler, and ensures that there is no loss of accuracy in any of the callbacks. The new handler for exact linear constraints manages the exact representation of linear constraints, their addition to LP data structures, exact feasibility checking and enforcement of solutions, constraint propagation routines, as well as some numerically safe techniques to improve efficiency such as a running error analysis for faster feasibility checks~\cite{EiflerGleixner2022}.

\subsubsection{Exact Presolving using \papilo}

Because the presolving library \papilo~\cite{GleixnerGottwaldHoen23} is fully templatized with respect to the number type of its input data and all arithmetical operations, it is ideally suited for performing presolving in \scip's exact solving mode.
For \scipv, the presolver \texttt{presol\_milp}, which interfaces to \papilo, has been extended to pass a rational representation of the problem instance to \papilo, set \papilo's numerical tolerances to zero, and call \papilo to perform all presolving steps in exact rational arithmetic.
Although presolving in exact rational arithmetic is slower than in floating-point arithmetic, the additional overhead can be reduced by executing \papilo in parallel, which can be enabled via parameter \texttt{presolving/milp/nthreads}, see \cite{EiflerGleixner2022} for a detailed evaluation.

\subsubsection{LP Relaxation, Exact LP Solving, and Numerically Safe Dual Bounds}
\label{sec:exact:safelp}

While \scip's standard LP structure is unmodified and maintains a floating-point approximation of the exact data, \scipv comes with a new exact LP structure that manages a rational version of the LP relaxation.
Similar to the LP interface for floating-point arithmetic, an exact LP solver interface is provided that can currently be used with the \soplex and \qsoptex \cite{Espinoza2006} LP solvers. Both support solving LPs exactly, with \soplex being the default exact LP solver in \scipv.

Because solving an LP exactly for every LP relaxation would be computationally expensive, exact \scip employs a selection of numerically safe dual bounding techniques to avoid calling the exact LP solver for every LP relaxation. These methods are called \emph{bound shift} and \emph{project and shift}, see \cite{CookKochSteffyetal2013,Jarck2020} for detailed descriptions of these methods.
The exact LP solver is called as a fallback if one of these cheaper methods fails, which may happen, e.g., if bound shift encounters unbounded variables.
By default, exact \scip also calls the exact LP solver in one of the following situations in order to improve overall performance:
\begin{itemize}
   \item at the end of the root node,
   \item at depth levels $2^k$ for $k = 1,2,3,\ldots$,
   \item when determining unboundedness,
   \item when the floating-point LP solution is close to being integer feasible,
   \item when the floating-point LP bound is close to the cutoff bound, or
   \item when all integer variables are fixed.
\end{itemize}

\subsubsection{Reliability Pseudocost Branching and Numerically Safe Cutting Planes}
\label{sec:exact:branchcut}

In exact solving mode, if the LP relaxation solution violates integrality by more than the floating-point feasibility tolerance, then the \texttt{integral} constraint handler calls the branching rule with highest priority that is marked exact.
In \scipv, the only branching rule made exact is reliability pseudocost branching, which is also the default branching rule in floating-point mode.
We refer to~\cite{CookKochSteffyetal2013} for a detailed description and computational study.
Note that the current version calls the floating-point LP solver for strong branching and the resulting bounds are not made safe. This means the results are currently not used to fix variables or add bound changes.
If all fractionalities are below the feasibility tolerance, then \scip simply branches on the first fractional variable.

Amongst the separator plugins, currently only the generation of Gomory mixed-integer cuts~\cite{Gomory1960} is available in exact solving mode. As shown in \cite{CookDashFukasawaGoycoolea2009}, these can be created through numerically safe aggregation of rows, yielding an approximation of an optimal LP tableau row. All MIR rounding operations \cite{MarchandWolsey1998} are also performed in double precision, using safe directed rounding. A custom cut postprocessing routine is used to avoid performance issues in exact LP solving~\cite{EiflerGleixner2024}.

\subsubsection{Safe Dual Proof Analysis and Constraint Propagation}
\label{sec:exact:prop}

As described in~\cite{borst2024certifiedconstraintpropagationdual}, the dual proof analysis from~\cite{witzig_computational_2021} has been made safe and learns exact linear constraints from nodes pruned due to an infeasible or bound exceeding LP relaxation.
Similar to the numerically safe Gomory mixed-integer cuts, safe dual proof analysis uses floating-point arithmetic with directed rounding to compute a weighted aggregation of the constraints given by the dual solution of the floating-point LP.
Like in floating-point mode, the learned clauses are marked not to be checked nor separated into the LP relaxation, and only propagated.
To this end, the exact linear constraint handler is equipped with safe constraint propagation using floating-point arithmetic with directed rounding.
This propagation is also applied to all other, model-defining constraints of the exact linear constraint handler.
For details and a computational study, see~\cite{borst2024certifiedconstraintpropagationdual}.

\subsubsection{A Repair Mechanism for Solutions from Floating-Point Heuristics}
\label{sec:splitExact1}

Although primal heuristics can not provide wrong results as long as the solutions they produce are correctly checked for feasibility, some of \scip's default heuristic plugins also derive problem reductions, e.g., constraints that express the infeasibility of a set of variable fixings.
In \scipv, all default floating-point heuristics are marked as safe to use in exact solving mode, because all unsafe reductions have been deactivated in exact solving mode.
However, the chance that a heuristic that is based on floating-point arithmetic and that uses error tolerances can produce an exactly feasible solution is rather low.
Motivated by this, the exact solving mode provides a new constraint handler \texttt{exactsol}, which postprocesses candidate solutions from primal heuristics in order to make them exactly feasible.
Although a constraint handler, the plugin does not constitute a constraint and is not designed to check solutions for feasibility.
Instead, it takes an approximately feasible floating-point solution, rounds and fixes the values of all integer variables, and solves the remaining LP with a call to the exact LP solver.
In order to control and save computational effort, candidate solutions can be buffered and only highly promising candidate solutions are processed immediately.
Periodically, buffered solutions are processed in order of best approximate objective function value.
For a more detailed discussion, see~\cite{EiflerGleixner2022}.

\subsubsection{Certification and Verification}
\label{sec:splitExact2}

Mathematically proven exactness of an algorithm is clearly not sufficient to ensure that the algorithm is implemented without errors. For a program of vast complexity such as \scip, it is currently impossible to formally verify the correctness of the implementation. Therefore, \scipv provides the option to write a certificate file in the \texttt{VIPR} format~\cite{VIPR}, which can be verified independently of the solving process.
While this creates an overhead both for proof logging and proof checking, the overhead does not exceed the time needed for solving. While this does not guarantee correctness of the implementation, it does provide a way to guarantee that the reasoning of the solving process and the optimal solution for a given problem instance are correct.

A certificate file contains the following information:

\begin{itemize}
   \item the problem statement in exact arithmetic,
   \item an optimal solution and claimed primal and dual bounds on its objective value,
   \item the derivation section, which in total certifies optimality of the solution.
\end{itemize}
The derivation section is a sequential encoding of the branch-and-bound tree that was traversed during the solving process. It is restricted to the following arguments:
\begin{enumerate}
   \item Conical combinations of inequalities: Given a system of valid linear inequalities~$Ax \le b$, new constraints $c^{\T}x \le d$ can be added to the certificate by supplying coefficients~$\lambda_i \ge 0$ such that $\sum_i \lambda_i A_i = c$ and $\sum_i \lambda_i b_i = d$, where~$A_i$ is the~$i$th row of~$A$. This argument is used to prove dual bounds obtained from LP relaxations, but also for the aggregations used in creating cuts and for dual proof conflicts, cf. Sections~\ref{sec:exact:branchcut} and~\ref{sec:exact:prop}.
   \item Chvátal-Gomory rounding: Given a valid inequality $a^{\T}x \le b$ with all coefficients~$a_i$ and variables $x_i$ integer, the right-hand side can be rounded to the next integer value~$b' = \lfloor b \rfloor$.
   \item Disjunction logic: New constraints can be added to the certificate as assumptions, which are then used to prove the conditional validity of other constraints. Eventually, a constraint that is valid given two assumptions that form a disjunction can also be made valid without the assumptions.
     This ``unsplitting'' argument is used to derive dual bounds of a parent node from its children and to certify validity of the mixed-integer rounding formula, see~\cite{EiflerGleixner2024}.
\end{enumerate}
For more technical details on this certificate format, we refer to~\cite{VIPR}.

Because correctness of the presolving process can currently not be certified in this format, \scipv produces the above certificate for the LP-based branch-and-cut process of the transformed problem after presolving.
In addition, a certificate is written to allow verifying the primal bound for the original problem by checking feasibility of the final solution in the original problem space.

A C++ proof checker to verify the correctness of these certificate files is available as part of the \scipopt GitHub repository~\cite{viprgithub}.
A more rigorous proof checker formally verified in HOL4 is available through
the CakeML framework \cite{GithubCakeML}.

\subsubsection{Computational Performance of Exact \scip}

In order to analyze the performance of the current state of the exact solving mode, it is insightful to compare it to two different variants of floating-point \scip: the default setting (``fp-default'') and the reduced variant (``fp-reduced''), where all features not available in exact solving mode are disabled. The test was performed on the \miplib~2017 benchmark testset \cite{MIPLIB} with a two hour time limit and three random seeds per instance. All reported times are shifted geometric means with a shift of one; node counts are reported with a shift of $100$. We report results both for the subset of instances that could be solved to optimality by at least one setting (one-solved), and for the subset of instances that could be solved to optimality by all settings (\alloptimal).

The exact setting was able to solve $161$ instances to optimality, while the reduced setting was able to solve $235$ instances (+$46\%$), and the default setting was able to solve~$342$ instances (+$112\%$).
In terms of solving time, we observe a slowdown factor of~$2.6$ for the exact solving mode compared to the reduced mode for one-solved and of~$3.6$ on \alloptimal. When comparing against the default version of \scip, that factor is $10.8$ for one-solved and $6.8$ for \alloptimal.

It is noteworthy that the number of nodes explored in the exact solving mode is significantly higher than in the reduced mode (+$155\%$), which has several causes such as the fact that strong branching reductions are disabled in exact solving mode, that exactly feasible primal solutions are harder to find, or the requirement to reach a primal dual gap of exactly zero.
Overall these results suggest a slowdown factor between $3$ and~$4$ for solving MIPs exactly with the same features and a current slowdown of around $10$ times compared to fp-default. Note that this can be vastly greater or smaller depending on the problem instance.

\begin{table}
   \caption[Price of exactness]{Performance comparison of exact \scip against floating-point \scip with default settings and in a reduced variant, where all features not available in exact solving mode are disabled, taken from~\cite{Eifler2025}. Times are in seconds.}
\label{tbl:final}
\centering
\begin{tabular}{@{}llrr
   r
   l
   r
   l@{}}
   \toprule
    Testset                      & setting & count & solved & {time}  & {(relative)} & {nodes} & {(relative)} \\
   \midrule
   \multirow{3}{*}{one-solved} & exact   & 350   & 161    & 2731.2 & {--}    & 21109.3 & {--}    \\
                               & fp-reduced & 350   & 235    & 1059.3 & 0.39    & 14184.5 & 0.67 \\
                               & fp-default & 350   & 342    & 255.0  & 0.10    & 2182.2  & 0.10  \\
   \cmidrule{1-8}
   \multirow{3}{*}{\alloptimal} & exact   & 153   & 153    & 798.2  & {--}    & 8462.8  & {--}    \\
                               & fp-reduced & 153   & 153    & 223.0  & 0.28    & 3322.1  & 0.39 \\
                               & fp-default & 153   & 153    & 118.2  & 0.15    & 1272.7  & 0.15 \\
   \bottomrule
\end{tabular}
\end{table}


\subsection{Presolvers}
\label{sec:presolve}


\subsubsection{Implied Integrality Detection using Network Submatrices}

\emph{Implied integrality}, also referred to as implicit integrality, refers to implication of integrality of some variables by integrality of others in combination with the constraints.
Already since its inception, \scip detects implied integer variables using a primal detection method that detects implied integrality for variables with $\pm 1$ coefficients in equations with integral coefficients and right-hand side.
In Version 1.1, a dual detection method was added to \scip, which detects implied integrality for variables that appear in inequalities only with $\pm 1$ coefficients.
Although both methods have existed in SCIP and other MILP solvers for a long time already, they were only first described in detail in \cite{AchterbergBixbyGuetal.2019}.
The authors of \cite{AchterbergBixbyGuetal.2019} conduct experiments with Gurobi and show that enabling implied integer detection reduces the solving time by 13\,\% and the node count by 18\,\% on models that take longer than 10 seconds to solve in their internal testing set.

Implied integer variables can be beneficial for performance in several ways.
Implied integer variables that are continuous variables in the original model may be used in branching.
Furthermore, their integrality can be exploited to derive stronger cut coefficients, most notably in Mixed-Integer Rounding cuts~\cite{MarchandWolsey1998}, and can enable new classes of cutting planes that require rows with integer variables only.
Their integrality can also be used in propagation techniques and conflict analysis to derive stronger bounds.
For implied integer variables that are integer in the original model, the redundancy of integrality constraints may be helpful, and can be exploited in primal heuristics such as diving methods or local neighborhood search to limit the number of branching candidates that are to be fixed or explored.
Moreover, it is not necessary to branch on implied integer variables within branch-and-bound.

Both the primal and the dual detection method can only detect implied integrality of one variable at a time.
Two of the contributors recently investigated implied integrality for MILP problems~\cite{HulstW25}.
A sufficient condition for implied integrality was identified to be integrality of the polyhedron that remains after fixing a subset of integer variables to integer values, independent of the fixed values.
To leverage this insight, they formulate a method, called \emph{TU detection}, to detect implied integrality by detecting \emph{totally unimodular} submatrices in the constraint matrix.
The primary advantage of this algorithm compared to the primal detection method and the dual detection method is that it can detect implied integrality of more than one variable at a time.

We now present the TU detection algorithm as it was implemented in \scipv.
For the MILP problem as given in \eqref{eq:generalmilp}, it attempts to find a partition of the variable set $\mathcal{N}$ into three pairwise disjoint sets $\mathcal{S}$, $\mathcal{T}$, and $\mathcal{U}$.
The set $\mathcal{S}$ must consist of integer variables, while $\mathcal{T}$ and $\mathcal{U}$ may consist of both continuous and integer variables.
The TU detection algorithm requires the problem to be of the form
\begin{multline*}
	P = \{ (x, y, z) \in \Z^{\mathcal{S}} \times \R^{\mathcal{T}} \times \R^{\mathcal{U}} : A x + B y \leq d,~ E x + F z \leq h,~ \ell \leq y \leq u,~ \\
  y_i \in \Z ~\text{ for all }i \in \mathcal{T} \cap \mathcal{I},~ z_i \in \Z ~\text{ for all } i\in \mathcal{U} \cap \mathcal{I} \},
\end{multline*}
where $A$ and $d$ are integral, the bounds $\ell$ and $u$ are integral or infinite, and $B$ is totally unimodular.
The variable bounds for $x$ and $z$ can be assumed to be implicit in the other constraints.
Then, after fixing the $x$-variables to any integer vector $\bar{x} \in \Z^{\mathcal{S}}$, the remaining problem $Q(\bar{x}) = \{ (x,y,z) \in P : x = \bar{x}\}$ can be written as
\begin{multline*}
	Q(\bar{x}) = \{ (x,y,z) \in \Z^{\mathcal{S}} \times \R^{\mathcal{T}} \times \R^{\mathcal{U}} : B y \leq d - A \bar{x},~ F z \leq h - E \bar{x},~ \ell \leq y \leq u~, ~x = \bar{x}, \\
 y_i \in \Z ~\text{ for all } i \in \mathcal{T} \cap \mathcal{I},~ z_i \in \Z ~\text{ for all } i\in \mathcal{U} \cap \mathcal{I} \}.
\end{multline*}
Note that $Q(\bar{x})$ is the Cartesian product of two disconnected subproblems over~$y$ and~$z$ (and the fixed vector $\bar{x}$).
The linear programming relaxation of the subproblem over~$y$ is given by $Q^y(\bar{x}) = \{ y \in \R^{\mathcal{T}} : B y \leq d - A \bar{x},~ \ell \leq y \leq u\}$.
Since $B$ is totally unimodular, the bounds $\ell$, $u$ are integral and since the right-hand side $d - A \bar{x}$ is integral, a well-known theorem by Hofmann and Kruskal~\cite{Hoffman1956} shows that $Q^y(\bar{x})$ is an integral polyhedron. 
Provided that the objective considered for $Q^y(\bar{x})$ is linear, there exists an optimal solution~$\bar{y}$ of $Q^y(\bar{x})$ that is integral for all $\bar{x}\in \Z^S$.
Since the $y$- and $z$-variables in $Q(\bar{x})$ are disjoint, such a $\bar{y}$ is also integral, feasible and optimal for $Q(\bar{x})$.
Furthermore, the feasible solutions of $P$ all lie in $Q(\bar{x})$ for some $\bar{x}\in\Z^S$, and each such $Q(\bar{x})$ admits an optimal solution $(\bar{x},\bar{y},\bar{z})$ where $\bar{y}$ is integral. 
This shows that it is safe to assume integrality of the $y$-variables when optimizing a linear objective over $P$. In this sense, the $x$-variables imply the integrality of the $y$-variables for $P$.
In~\cite{HulstW25}, a more rigorous definition of implied integrality and formal proofs of the TU detection method are presented.

Although totally unimodular matrices can be detected in polynomial time~\cite{Truemper1990}, existing implementations of the detection algorithm have high running times~\cite{Walter2013}.
The implied-integer detection algorithm instead only detects \emph{network matrices} and transposed network matrices, which form a large subclass of totally unimodular matrices (see~\cite[Chapter 19]{Schrijver86}).
In order to detect network matrices, the algorithm uses fast augmentation methods that can determine whether a network matrix can be extended by a new column~\cite{BixbyWagner1988} or row~\cite{HulstW24}.
These methods are then used iteratively to grow a (transposed) network matrix one column at a time.

Note that in the description of $P$, $y$ and $z$ must be in disconnected blocks, and $x$ must be integral. Thus, the detection method first considers each connected block of the submatrix formed by the continuous columns, and determines if it is network, transposed network or neither, and if it satisfies the given integrality requirements.

After determining implied integrality of the continuous variables, the submatrix $B$ is greedily grown with columns of variables with integrality constraints by using the network matrix augmentation methods. 

\begin{table}[tb]
	\caption{%
		Implied integer detection on 1035 instances of the MIPLIB 2017 collection set~\cite{MIPLIB}.
		For a set $X\subseteq \mathcal{N}$, $X_{\mathrm{int}}$ ($X_{\mathrm{con}}$) indicates the set of variables in $X$ that (do not) have integrality constraints in the original model. $\mathcal{T}$ indicates the set of implied integer variables.
		The numbers and ratios reported are based on the presolved model, and averaged over all instances tested. 
		The number of variables reported is the shifted geometric mean of the number of variables with shift 10.
		The detection time reported is the shifted geometric mean with a shift of \SI{1}{\milli\second}.
		All other means are arithmetic means.
	}
	\label{tab:implint_detection_stats}
        \centering
		\begin{tabular}{@{}lrr@{}}\toprule
			Method & SCIP~9.2.1 & SCIP~10.0.0\\
			\midrule
			affected instances & 203 & \textbf{712} \\
			mean of ratio $\frac{|\mathcal{T}|}{|\mathcal{N}|}$ & \SI{3.3}{\percent} & \textbf{\SI{18.8}{\percent}}\\
			mean of ratio $\frac{|\mathcal{T}_{\mathrm{con}}|}{|\mathcal{N}|}$& \SI{3.2}{\percent} & \SI{7.5}{\percent} \\
			mean of ratio $\frac{|\mathcal{T}_{\mathrm{int}}|}{|\mathcal{N}|}$& \SI{0.1}{\percent} & \SI{11.3}{\percent} \\
			mean of ratio $\frac{|\mathcal{N}_{\mathrm{con}}\setminus \mathcal{T}|}{|\mathcal{N}|}$& \SI{26.5}{\percent} & \SI{22.2}{\percent} \\
			mean of ratio $\frac{|\mathcal{N}_{\mathrm{int}}\setminus \mathcal{T}|}{|\mathcal{N}|}$& \SI{70.2}{\percent} & \textbf{\SI{59.0}{\percent}} \\ 
			number of variables & 6937 & 6944 \\
			mean detection time & -- & \textbf{13 ms} \\\bottomrule
		\end{tabular}
\end{table}

Some statistics of the new detection method can be found in \cref{tab:implint_detection_stats}. It presents results of the TU detection method on MIPLIB 2017, where 5 instances were excluded due to memory limits and 25 further instances were excluded because they were solved during presolve by both SCIP versions.
On average, the new TU detection method is very effective, and detects implied integrality for \SI{18.8}{\percent} of the variables, compared to~\SI{3.3}{\percent} detected using the primal and dual methods. Moreover, the new TU detection method is quite fast and typically runs in a few milliseconds.
The changed variable type distribution poses a few challenges.
Since (implied) integrality of a variable is a key property of a MIP problem, nearly every algorithm within \scip uses it, and the performance of \scip is tuned relative to the variable distribution.
Thus, adjusting \scip to effectively use the detected implied integrality is time- and resource-intensive and requires an ongoing effort.
With the release of \scipv, many important parameters and decisions within the cutting plane and branch-and-bound frameworks have been re-tuned to account for the changed variable type distribution.
However, the performance impact of these changes is still only relatively neutral, which is in stark contrast to the results obtained by Gurobi~\cite{AchterbergBixbyGuetal.2019}.  Consequently, implied integer detection is not enabled by default yet in \scipv. In future releases, we hope to further improve the exploitation of implied integrality within SCIP.

\scipv has several API changes that are closely connected to implied integers.
First of all, the implied integer variable type flag is deprecated.
Instead, an additional flag that, separately from the variable type, indicates whether a variable is implied integer or not.
The new flag can take three values; one for no implied integrality, one for implied integrality such as derived by the TU detection method and dual detection method, and one for the implied integrality detected by the primal detection method.
The latter is distinguished because it guarantees that a variable takes an integral value in every solution after fixing other integers.
This strong property can be used by \scip in some algorithms, most notably when aggregating variables.
In contrast, the dual detection method and the TU detection method only guarantee the existence of an integral solution value, but do not enforce it.
This may be disadvantageous, as there are cases where the exact integral values of the $y$ variables are desired, such as primal solutions. For example, it may happen that a primal heuristic uses implied integrality and finds a solution with fractional $y$ that otherwise satisfies the integrality constraints. Then, one needs to find a vertex solution of the LP relaxation of $Q^y(\bar{x})$ to recover the integral solution of the $y$ variables, which can be costly to compute.
One additional advantage of the new implied integer flag is that it enables running the dual detection method for variables with integrality constraints, which was previously not possible in \scip.


\subsection{Symmetry Handling}
\label{sec:symmetry}


A symmetry of an \MILP or \MINLP transforms feasible
solutions into feasible solutions while preserving their objective values.
When solving problems containing symmetries, (spatial)
branch-and-bound algorithms arguably create symmetric subproblems,
which contain equivalent information.
Branch-and-bound can thus be accelerated by removing symmetric subproblems,
which is the goal of symmetry handling.

Since Version~5.0, many methods for handling so-called permutation
symmetries have been added to \scip.
\scip~\scipversion extends many of these methods to
reflection symmetries, see below for a definition of
both types of symmetries, and introduces new methods for special symmetries.
Moreover, \scip's capabilities to detect symmetries have been extended to
problems containing Pseudo-Boolean or disjunctive constraints.

\subsubsection{Handling Reflection Symmetries}

\scip can detect and handle two types of symmetries:
permutation symmetries and reflection symmetries.
A map~$\gamma\colon \R^n \to \R^n$ is a \emph{permutation symmetry} of an
\MILP or \MINLP if it defines a symmetry of the respective problem, as
defined above, and there exists a permutation~$\pi$ of~$[n]$ such that~$\gamma(x) =
(x_{\pi^{-1}(1)},\dots,x_{\pi^{-1}(n)})$ for all $x \in \R^n$.
That is, permutation symmetries reorder the entries of a variable vector.
A symmetry~${\rho\colon\R^n \to \R^n}$ of an \MILP or \MINLP is called a
\emph{reflection symmetry} if there exists a permutation~$\pi$ of~$[n]$ and
a vector~$s \in \{-1,1\}^n$ such that~$\rho(x) = (s_1
x_{\pi^{-1}(n)},\dots,s_n x_{\pi^{-1}(n)})$ for all $x \in \R^n$.
That is, reflection symmetries reorder the entries of a variable vector
like a permutation symmetry, but additionally, entries can be negated
(reflected at the origin).
Actually, \scip can detect more general reflection symmetries
that reflect variables at a point different than the origin, e.g.,
binary variables can be reflected at~$\frac{1}{2}$ to exchange 0-entries
and 1-entries.
For the ease of exposition, the presentation is restricted to reflections at
the origin though; see~\cite{Hojny2024+} for a discussion of the
general case.
The remainder of this section provides an overview of the symmetry
handling methods that have been made compatible with reflection symmetries
in \scip~\scipversion, and new methods that have been added to \scip.

\paragraph{SST Cuts}
Let~$\Gamma$ be a group of permutation symmetries.
Liberti and Ostrowski~\cite{LibertiOstrowski2014} and
Salvagnin~\cite{Salvagnin2018} introduced \emph{SST cuts} as a family of
symmetry handling inequalities, which take the form~$x_i \geq \gamma(x)_i$
for all~$i \in [n]$ and~$\gamma \in \Gamma_i$, where~$(\Gamma_i)_{i \in [n]}$ is a
carefully selected family of subgroups of~$\Gamma$.
In \scip~\scipversion, a version of SST cuts has been implemented that
replaces the permutation symmetry~$\gamma$ by a reflection symmetry~$\rho$
(not necessarily with respect to the origin).
That is, the generalized SST cut is defined as~$x_i \geq \rho(x)_i$.

\paragraph{Orbitopes}
In many applications, such as graph coloring~\cite{KaibelPfetsch2008}, the
variables of an \MILP or \MINLP can be arranged in a matrix~$X \in \R^{s
  \times t}$ such that the symmetries can reorder the columns of the matrix
arbitrarily.
To handle such column symmetries, which are called orbitopal symmetries, a
standard approach is to enforce that a solution has lexicographically
sorted columns.
To achieve this sorting, previous versions of \scip implemented
both symmetry handling inequalities~\cite{Hojny2020,HojnyPfetsch2019} as
well as the propagation algorithm orbitopal
reduction~\cite{BendottiEtAl2021,DoornmalenHojny2024a}.
Moreover, if the matrix~$X$ consists of binary variables and every feasible
solution has at most or exactly one~1-entry per row, specialized cutting
planes~\cite{KaibelPfetsch2008} and propagation
algorithms~\cite{KaibelEtAl2011} can be used.

In \scip~\scipversion, the methods for orbitopes have been extended by also
taking reflection symmetries into account.
When \scip detects orbitopal symmetries for a variable matrix~$X$, it is now
also checked whether there is a reflection symmetry that simultaneously
reflects all variables within a single column of~$X$.
If this check evaluates positively, by symmetry of all columns of~$X$,
every column admits such a reflection.
Hence, by possibly applying these reflection symmetries, one can always
guarantee that there is a feasible solution in which all entries in the first row
of~$X$ are nonnegative.
The nonnegativity of the first row is then enforced by adding the
inequalities~$X_{1,j} \geq 0$ for all~$j \in [t]$.

\paragraph{Double-Lex Matrices}
In applications like the disk packing problem~\cite{Szabo2005},
richer symmetries than orbitopal symmetries arise.
The goal of one of this problem's variants is to find the maximum value~$r$
such that~$n$ nonoverlapping disks of radius~$r$ can be packed
into~$[-1,1]^2$.
This problem can be modeled by introducing a variable for the radius and a
variable matrix~$X \in [-1,1]^{n \times 2}$, in which entry~$X_{i,j}$
models the~$j$-th coordinate of the center of disk~$i$.
Given a solution matrix~$\bar{X}$, equivalent solutions can be found by
permuting the rows of~$\bar{X}$ arbitrarily (all disks are identical).
Moreover, since the disks are packed into a square, the order of the
columns of~$\bar{X}$ can be exchanged, and each individual column can be
reflected.
That is, there are orbitopal symmetries for both exchanges of columns and
rows, and columns can be reflected individually.

In \scip~\scipversion, a heuristic has been added to detect a
generalization of such symmetries.
This heuristic tries to identify variable matrices with block shape
\[
  B
  =
  \begin{pmatrix}
    B_{1,1} & \dots & B_{1,t}\\
    \vdots & \ddots & \vdots\\
    B_{s,1} & \dots & B_{s,t}
  \end{pmatrix},
\]
where each block~$B_{i,j}$ is a submatrix of compatible dimensions and
there are orbitopal symmetries on the row and column blocks of this matrix.
That is, for each row and column block
\[
  B_{i,\cdot} = (B_{i,1},\dots,B_{i,t}),\; i \in [s],
  \text{ and }
  B_{\cdot,j} = (B^\T_{1,j},\dots,B^\T_{s,j})^\T,\; j \in [t],
\]
respectively, there is an orbitopal symmetry on the rows of~$B_{i,\cdot}$
and columns of~$B_{\cdot,j}$.
The row and column symmetries of the disk packing problem arise then as a
special case, in which~$s = t = 1$.

To handle orbitopal symmetries of row and column blocks, \scip applies
orbitopal reduction on the submatrices~$B_{i,\cdot}$ and~$B_{\cdot,j}$,
respectively.
That is, both the rows and columns can be sorted lexicographically, leading
to a so-called double-lex matrix.
For the special case~$s = t = 1$, \scip additionally checks whether there
exists a reflection symmetry that simultaneously reflects all variables
within a single column.
If this is the case, due to the presence of row symmetries, one can assume
that the first half of the variables in the first column only take
nonnegative values; see Khajavirad~\cite{Khajavirad2024}.
Moreover, some entries of columns~$2,\dots,t$ can be assumed to be
nonnegative, too; see~\cite{Hojny2024+}.
This is enforced by adding the inequalities~$X_{i,j} \geq 0$ for the
nonnegative entries~$(i,j)$.

\paragraph{Simple Symmetry Handling Inequality}
Suppose there exists a reflection symmetry that simultaneously
reflects all variables of a problem, but does not change the order of the
variables, i.e., the permutation~$\pi$ of this reflection is the identity.
This symmetry can be handled by the simple inequality~$\sum_{i = 1}^n x_i
\geq 0$; see~\cite{Hojny2024+}.
Because this inequality is rather weak, it is only applied when \scip does
not detect that any other symmetry handling method seems to be useful.

\subsubsection{Extensions of Symmetry Detection}
Symmetries of an \MILP or \MINLP are usually detected by computing
automorphisms of a suitably defined graph; see
Salvagnin~\cite{Salvagnin2005} and Liberti~\cite{Liberti2012a}.
Since \scip~9.0, symmetry detection graphs are created via
an optional callback for constraint handlers; see~\cite{Hojny2024+} for a
detailed discussion.
Symmetries of a problem can then be detected if all active constraints
handlers implement this callback.

In \scip~\scipversion, these symmetry detection callbacks have been
implemented for Pseudo-Boolean and disjunctive constraints, thus extending
the range of problems for which \scip can detect symmetries.

\subsubsection{Further Enhancements}

To compute automorphisms of symmetry detection graphs, previous versions of
\scip provide interfaces to the external graph automorphism software \bliss~\cite{JunttilaKaski2007,JunttilaKaski2011}
and \nauty~\cite{Nauty}.
The release of~\scip~\scipversion adds another interface to the software \dejavu~\cite{Anders2024,AndersS21,AndersSS2023}.
In contrast to previous versions, \bliss is no longer shipped with \scip;
instead, \nauty is provided with the \scip release.
Since we observed that the software \nauty can spend a substantial amount of
time on computing automorphisms, a work limit has been introduced to
terminate \nauty prematurely.
This limit can be adapted via the parameter
\param{propagating/symmetry/nautymaxlevel}, which controls the maximum
depth level of \nauty's search tree.

Furthermore, some parts of the code have been refactored.
Previous versions of \scip contained the constraint handler
plugin \plugin{cons\_orbitope} to handle both plain orbitopes and those
incorporating set packing and partitioning constraints on the rows.
To save memory, \scip~\scipversion contains two new plugins
\plugin{cons\_orbitope\_full} and \plugin{cons\_orbitope\_pp} to handle
the different types of orbitopes, respectively.
To ensure backward compatibility, the old plugin still exists and serves
as an interface to the two new ones, i.e., the old API can still be used to
add orbitope constraints.


\subsection{Cut-based Conflict Analysis}
\label{sec:conflict}


In \scip, conflict analysis is applied whenever a local infeasible constraint is
detected at a node of the branch-and-bound tree. A common source of such
detections is propagation, which may occur after multiple rounds on the current
node. The goal of conflict analysis is to ``learn'' a constraint that would
have revealed the infeasibility earlier. Such a constraint would have
propagated at an earlier node, preventing \scip from visiting the infeasible
node in the first place. Up to \scip 9.0, infeasibilities were analyzed with a
SAT-inspired, graph-based conflict analysis \cite{achterberg2007conflict}. This
approach represents the sequence of branching decisions and propagated
deductions that led to the contradiction as a directed acyclic graph. From this
graph, a valid constraint can be inferred by identifying a subset of bound
changes that separates the decisions from the node where the contradiction was
reached. The learned constraints are disjunctive constraints, are implied by the original
problem, and are therefore globally valid.

As of \scip 10.0, cut-based conflict analysis was introduced \cite{mexi2024cut}.
The method is inspired by conflict-analysis techniques from Pseudo-Boolean
solvers \cite{chai2003fast,le2011sat4j,elffers2018divide,mexi23improving} and operates
directly on more expressive linear inequalities rather than conflict graphs. It proceeds via
linear combinations, integer roundings, and cut generation. During
branch-and-bound, each branching decision or propagation tightens the local domain~$L_T$. 
If propagation detects infeasibility (i.e., when
propagating a constraint $C_\texttt{conf}$ on the local domain $L_T$), conflict analysis is
triggered. In this context, $C_\texttt{conf}$ is the conflict constraint, while the
constraint $C_\texttt{reason}$ that defined $L_T$ (that is, $C_\texttt{reason}$ propagated under the
previous local domain~$L_{T-1}$) is the reason constraint. Together, $L_{T-1}$,
$C_\texttt{reason}$, and $C_\texttt{conf}$ form an infeasible system. Each iteration aims to derive a
single constraint $C$ that is globally valid and infeasible in $L_{T-1}$.
Once such a constraint $C$ is identified, it becomes the new conflict constraint $C_\texttt{conf}$. 
In each iteration, the reason constraint $C_\texttt{reason}$ is set to the constraint that
propagated the most recent bound change contributing to 
the infeasibility of $C_\texttt{conf}$ and the algorithm continues until a globally valid constraint is
obtained that could have prevented the infeasibility before the latest branching.

\begin{figure}
    \centering
    \includegraphics[width=1.0\textwidth]{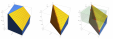}
    \caption{Left: Region defined by $C_\texttt{reason}$, $C_\texttt{conf}$, and the global domain.
    Middle: Mixed-integer hull, where \texttt{x1} is the only integer variable.
    Right: Mixed-integer hull intersected with the local domain. In particular, two integer-feasible points (leftmost and rightmost red points)
    and a convex combination of them that intersects the local domain (middle red point).}
    \label{fig:mbp-reason-reduction}
\end{figure}

In an iteration of conflict analysis, it may not be possible to derive a
globally valid conflict constraint $C$ using only the inequalities
$C_\texttt{reason}$ and $C_\texttt{conf}$. 
Figure~\ref{fig:mbp-reason-reduction} illustrates this situation. The left plot
shows the region defined by the reason constraint $C_\texttt{reason}$, the
conflict constraint $C_\texttt{conf}$, and the global domain. The middle plot
shows the mixed-integer hull of this region, where \texttt{x1} is the only integer variable.
The right plot shows the intersection of the mixed-integer hull with the local domain.
In particular, the leftmost and rightmost red points are integer-feasible, while the middle red point is a
convex combination of them that intersects the local domain.
Because the local domain intersects the mixed-integer hull, no globally valid
inequality can be derived that separates the local domain from the hull.
In such cases, the reason constraint $C_\texttt{reason}$ needs to be
modified to ensure that the new conflict constraint $C_\texttt{conf}$ remains
infeasible in the previous local domain. This is done by applying so-called
\emph{reduction rules} to the reason constraint that can be activated by
setting the parameter \texttt{conflict/reduction} to \texttt{m} (which stands for \texttt{MIR} reduction) 
and, for the case of mixed-binary constraints, set
\texttt{conflict/mbreduction} to \texttt{true}.
To activate cut-based conflict analysis, the parameter
\texttt{conflict/usegenres} has to be set to \texttt{true}.




\subsection{Separators}
\label{sec:separators}

%
%

A new separator \texttt{sepa\_flower} dealing with products of nonnegative variables was added to \scip.
It works with logical AND constraints $z = x_1 \cdot x_2 \cdot \dotsc \cdot x_\ell$ for binary variables~${x_i \in \{0,1\}}$ and products appearing in nonlinear constraints.
In the former, $z$ is called \emph{resultant} and the $x_i$ are called \emph{operands}.
The purpose of the new separator is to generate cutting planes that strengthen the LP relaxation in the presence of multiple such constraints.
The collection of all constraints can be represented via the \emph{multilinear hypergraph} $G = (V,E)$, which has a node $v \in V$ per variable $x_v$ that appears as an operand in at least one AND constraint and a hyperedge $e \in E$ per constraint, where the resultant variable satisfies $z_e = \prod_{v \in e} x_v$.

This hypergraph is constructed after presolving the problem by inspecting all AND constraints and scanning expression trees of nonlinear constraints.
Moreover, \emph{overlap sets}, which are sets of the form $e \cap f$ for $e$, $f \in E$, are gathered using hash tables, from which \emph{compressed sparse row} and \emph{column} representations of the incidences among hyperedges, overlap sets and nodes are computed.

While there exist many classes of valid inequalities~\cite{CramaR17,DelPiaK17,DelPiaK18,DelPiaD21,DelPiaW22}, we restricted the current implementation to \emph{$k$-flower inequalities}~\cite{DelPiaK18} for $k=1$, $2$, since the separation problem can be solved very efficiently once the overlap sets are available.
Such an inequality is parameterized by one \emph{center} edge $e \in E$ and $k$ \emph{neighbor edges} ${f_1,f_2,\dotsc,f_k \in E}$, all of which must intersect $e$.
It reads
\[
  z_e + \sum_{i=1}^k (1-z_{f_i}) + \sum_{v \in R} (1-x_v) \geq 1,
\]
where $R \coloneqq e \setminus \bigcup_{i=1}^k f_i$ denotes the nodes of $e$ that are contained in no neighbor edge.
By simple enumeration, the separation problem can easily be solved in time $\mathcal{O}(|E|^{k+1})$.
Under the reasonable assumption that the size $|e|$ of every edge is bounded by a constant, this can be reduced to $\mathcal{O}(|E|^2)$ as done in the implementation~\cite{DelPiaKS20}.
By exploiting overlap sets using the implemented data structures, this is reduced to $\mathcal{O}(|E|)$.

Hyperedges for the new separator can, in principle, also be derived from product expressions involving continuous variables.
However, this requires nonnegative lower bounds that are relaxed to $0$, and in our tests we observed a performance loss, which is why this is disabled by default.
Users can enable this feature via the parameter \texttt{separating/flower/scanproduct}.

\begin{table}[htpb]
  \caption{%
    Performance comparison for instances from MILP, MINLP and Pseudo-Boolean testsets.
    MILP and Pseudo-Boolean experiments were run using SoPlex~8.0.0 as an LP solver and those for MINLP using CPLEX~12.10.0.0.
    Handling products of continuous variables remained disabled for all experiments.
  }
  \label{tbl:sepa_flower}
\scriptsize

\begin{tabular*}{\textwidth}{@{}ll@{\;\;\extracolsep{\fill}}rrrrrrrrr@{}}
\toprule
               &                       &           & \multicolumn{3}{c}{Without \texttt{sepa\_flower}} & \multicolumn{3}{c}{With \texttt{sepa\_flower}} & \multicolumn{2}{c}{relative} \\
\cmidrule{4-6} \cmidrule{7-9} \cmidrule{10-11}
Testset        & subset                & instances &                                   solved &       time &        nodes &                                   solved &       time &        nodes &       time &        nodes \\
\midrule
MILP           & \cleaninst            &       698 &                                      604 &      291.7 &         2815 &                                      605 &      290.7 &         2810 &        1.00 &          1.00 \\
MILP           & \affected             &        14 &                                       13 &      369.4 &          261 &                                       14 &      313.2 &          223 &        0.85 &          0.85 \\
\midrule
MINLP          & \cleaninst            &       822 &                                      810 &       18.9 &         2407 &                                      810 &       18.7 &         2377 &        0.99 &          0.99 \\
MINLP          & \affected             &        21 &                                       21 &        8.6 &         1712 &                                       21 &        5.1 &         1023 &        0.60 &          0.60 \\
\midrule
Pseudo-Boolean & \cleaninst            &       450 &                                      423 &        7.8 &          333 &                                      423 &        7.7 &          331 &        1.00 &          0.99 \\
Pseudo-Boolean & \affected             &        11 &                                       11 &       12.0 &          359 &                                       11 &        9.7 &          258 &        0.81 &          0.72 \\
\bottomrule
\end{tabular*}
\end{table}

Computational results are shown in \cref{tbl:sepa_flower} for MILPs, MINLPs and Pseudo-Boolean instances, which give a clear advantage on affected instances although there are not too many.
However, thanks to the very efficient computation of the hypergraph $G$, the separator is performance-neutral in case multilinear structures are absent.


\subsection{Heuristics}
\label{sec:heuristics}


\subsubsection{Kernel Search}
\label{subsubsec:ks}

In~\scipv the Decomposition Kernel Search (\texttt{DKS}) heuristic is introduced.
This includes the basic framework of Kernel Search (KS) explained in the following.
The extensions of KS to respect additional decomposition information is described in the next subsection.
Originally, the heuristic was introduced in~\cite{Angelelli2010,Angelelli2012} to leverage a substructure on binary variables, which occurs in certain applications.
It was later refined to integrate general integer variables into its functionality~\cite{Guastaroba2017}, which serves as basis for the present implementation of KS, see~\cite{Halbig2025} for details.

In KS, a subset of the integer variables is defined as the \emph{kernel} $\mathcal{K}$.
It is supposed to contain variables that are likely to be active in a solution.
Here and in the remainder of the section, a variable is called \emph{active} in a solution when its value is nonzero and not at the lower bound.
In the opposite case, a variable is called \emph{inactive}.
Note that this definition of (in)activity is a generalization to \MILP problems of the definition in the classical literature, where typically only lower bounds of zero are considered.
One exemplary way to identify such a property is to consider active variables in the \LP solution.
The remaining integer variables $\intvarindex\setminus\mathcal{K}$ are divided into \emph{buckets}.
In the core procedure of KS, the kernel is united with a bucket~$\mathcal{B}$ and Problem~\labelcref{eq:generalmilp} is restricted to this union.
The restriction is obtained by fixing every integer variable not in the union to be truly inactive, in particular, to the lower bound, if it is finite, or to zero otherwise.
Note that variables without a finite lower bound but with a negative upper bound are considered active in the current implementation and thus remain in the kernel.

This identification of kernel and buckets as well as their handling is also performed for continuous variables, which is in contrast to the original sources but is based on the extensions below.
To force a binary or general integer variable from the bucket to be not at its lower bound and thus potentially be active, the parameter \param{heuristics/dks/advanced/addUseConss} controls the addition of the constraint
\begin{equation*}
	\sum_{i \in \mathcal{K}\cup \mathcal{B}} x_i \geq \sum_{i \in \mathcal{K}\cup \mathcal{B}}\lb_i + 1.
\end{equation*}
Here, only variables with finite lower bounds $\lb_i$ are considered.
After solving the restricted problem, the kernel is updated by adding active variables from the current bucket and deleting inactive kernel variables with respect to the obtained solution.
The update is performed whenever an improving solution is found.
This procedure is iterated until solutions for all kernel-bucket combinations are computed or resource limits are reached.

In~\cite{Halbig2025}, several extensions to the original work are introduced.
There may already exist a feasible solution to the original problem when the heuristic is called.
This can be exploited by using the parameter \param{heuristics/dks/usebestsol} to define the initial kernel $\mathcal{K}$ by considering the active variables in the currently best solution.
The remaining inactive variables are typically assigned to distinct buckets of similar size, after being sorted by their reduced cost in the \LP relaxation.
However, the reduced costs may be distributed in different orders of magnitude.
This motivates an assignment to buckets of different size, yet based on uniform ranges of the logarithm of the reduced costs.
The parameter \param{heuristics/dks/redcostlogsort} is provided to switch between these two approaches.

When solving the sequence of kernel-bucket-subproblems, the trade-off between solving times and objective improvements must be controlled carefully.
To do this, the gap limit for each of these subproblems is chosen adaptively based on whether the node limit of the previously solved problem is reached.

Two points are noteworthy with regard to the algorithmic scheme described:
First of all, there is an initial solve restricted to the kernel only.
Second, although it was introduced and defined for \MILPs, the above is no restriction for \MINLPs in theory.
However, extending the implementation to this general problem class requires special care.

\subsubsection{Decomposition Kernel Search}
\label{subsubsec:dks}

\citet{Halbig2025} extend KS explained in the previous section by the possibility to utilize decomposition information, called Decomposition Kernel Search (\texttt{DKS}).
Details on the representation of decomposition information can be found in the release report for \scip~7.0~\cite[Section~4.2]{SCIP7}.
Such decomposition information indicates relations among variables by specifying a block structure of the constraint matrix $A$ in~\labelcref{eq:generalmilp}.

Since empirical tests yielded active variables in \MILP-solutions for variables across all blocks and variable types, this is taken into account when creating the kernel and the buckets.
For this reason, kernel and buckets are determined for each block and variable type and are then united.
Since this strategy is performed across several levels, in particular for the blocks and variable types, the resulting structure of the final kernel and buckets is referred to as \emph{multi-level}.

There are different types of decomposition information depending on how the information was generated, e.g., through graph partitioning or flow modeling.
Computational experiments on instances associated with decompositions do not show a favor towards a special type of decomposition information but for specific problem classes.
In particular, \texttt{DKS} has a positive effect on solving binary problems, which KS was originally presented for.
On top of that, the solving times for problems with binary and integer variables without continuous variables are reduced.
The effect is attributed to the multi-level structure, as it seems to appropriately capture the occurrence of active variables in a solution.
Since it has a positive impact for certain problem types, but requires finetuning to be applicable in general, the heuristic is disabled by default.
For a detailed description please refer to~\cite{Halbig2025}.


\subsection{Branching Rules}
\label{sec:branching}


\subsubsection{Ancestral Pseudocosts}

Classical pseudocosts \cite{benichou1971pscost} estimate the immediate LP relaxation improvement from branching but ignore downstream effects. Our proposed approach, ancestral pseudocosts, aims to approximate this longer-term influence by aggregating improvements from subsequent levels, weighted by a discount factor. The following example illustrates how these updates are recorded during branching.

Assume we are at a node $N_0$ with LP relaxation objective value $L_0 \in \mathbb{R}$. Consider a candidate variable $x$ with LP value $x^{lp}_{0} \in \mathbb{R}$. When we branch on $x \leq \floor{x^{lp}_0}$, we create a new node $N_1$ with LP relaxation objective value $L_1 \in \mathbb{R}$. At node $N_1$, we update the traditional pseudocost $PS_0(x)$ for the variable $x$ by adding a record $(L_1 - L_0)/\{x^{lp}_0\}$, where $\{x^{lp}_0\} \defin x^{lp}_{0} - \floor{x^{lp}_0}$. Now consider the future-case where $N_1$ is the considered node and $y$ the candidate variable with LP solution value $y^{lp}_{1} \in \mathbb{R}$. When we branch on $y \leq \floor{y^{lp}_1}$, a new node $N_2$ is created with LP relaxation objective value $L_2 \in \mathbb{R}$. At $N_2$, we now update two records instead of one:
\begin{enumerate}
\item The traditional pseudocost $PS_0(y)$ of variable $y$ by adding datapoint $(L_2 - L_1)/\{y^{lp}_1\}$,
  where $\{y^{lp}_1\} \defin y^{lp}_{1} - \floor{y^{lp}_1}$.
\item First-level pseudocost $PS_1(x)$ of variable $x$ by adding datapoint $(L_2 - L_1)/\{x^{lp}_0\}$.
\end{enumerate}

In this example, part of the improvement in the LP relaxation bound at node \(N_2\) is attributed not only to the variable \(y\), which was branched on at \(N_1\), but also to the earlier decision to branch on variable \(x\) at \(N_0\). This attribution is made through the first-level pseudocost \(PS_1(x)\), which receives a discounted share of the bound improvement from~\(L_1\) to \(L_2\). In general, for a variable \(x\), we can maintain multiple levels of pseudocosts \(PS_k(x)\), where \(k\) denotes the number of levels between the branching decision on \(x\) and the resulting node where the LP improvement is observed. 

These contributions are aggregated using a geometric discount factor \(\gamma \in (0,1]\), leading to the following formulation for the \textit{ancestral pseudocost} of variable \(x\):

\begin{equation}
  APS(x) = PS_0(x) + \gamma PS_1(x) + \gamma^2 PS_2(x) + \dots
\end{equation}
Instead of selecting a branching variable based solely on traditional pseudocosts~\(PS_0(x)\), we compare variables using their ancestral pseudocosts \(APS(x)\). This encourages the solver to favor decisions that not only offer immediate LP bound improvements, but also lead to stronger relaxations deeper in the search tree.

Currently, ancestral pseudocosts are implemented using one-level. Such an implementation can be interpreted as an approximation of lookahead branching \cite{glankwamdee2006lookahead}. While traditional pseudocosts estimate the effect of strong branching at the immediate child node, incorporating one-level ancestral pseudocosts mimics lookahead branching by capturing downstream LP relaxation improvements, especially when the discount factor~\(\gamma = 1\). For ancestral pseudocosts to work well, both levels of pseudocosts---\(PS_0(x)\) and \(PS_1(x)\)---must be reliable. In \scip, we adopt a conservative strategy: if either level is deemed unreliable for any candidate variable, the branching rule falls back to the default (hybrid branching \cite{achterberg2005reliability} using only $PS_0(x)$) at that node.

To control the discount factor~$\gamma$ applied to ancestral pseudocosts
for the ``pscost'' and ``relpscost'' branching rules, the parameters
\param{branching/pscost/discountfactor} and
\param{branching/relpscost/discountfactor} have been added to \scip, respectively (default = \(0.2\)). An additional flag, \param{branching/collectancpscost}, enables or disables the recording of ancestral pseudocosts (default = ``false''). In the default settings of \scip~10.0, ancestral pseudocosts are disabled until conclusive evidence of improvement is found for general MILPs.

\subsubsection{Probabilistic Lookahead for Strong Branching}

The hybrid branching rule that \scip uses by default relies on strong branching calls to select branching variables, especially early in the tree before the pseudocosts of many variables are well initialized and considered reliable.
Strong branching is, however, expensive and careful working limits were added to avoid solving a prohibitive amount of LPs.
One such working limit is the so-called lookahead number, i.e., the number of consecutive variables that are evaluated without improvement to the dual bound estimate, after which we stop evaluating more candidates.
This maximum number is typically static, and a new dynamic criterion was introduced by \citet{mexi2023probabilistic}, based on an abstract probabilistic model of the branch-and-bound tree.
Importantly, experiments showed that the new dynamic criterion achieves a reduction in both runtime and tree size.
It reallocates strong branching calls at nodes where it is needed while stopping the candidate evaluation early when the probabilistic model estimates a low probability of finding a better candidate than the incumbent one.
The new criterion can be activated with the \texttt{branching/relpscost/dynamiclookahead} parameter.

\subsubsection{Mix Branching on Integer and Nonlinear Variables}

When solving problems with nonlinear nonconvex constraints, the constraint handler \texttt{cons\_nonlinear} would previously consider branching on a variable in a nonlinear nonconvex term only if there were no integer variables with fractional value in the solution to enforce.
Such a strict preference for integer variables did not seem appropriate, but it was also unclear when branching on a nonlinear variable was preferable.
To support ongoing investigations in this topic, with \scip~9.1.0 the enforcement priority of \texttt{cons\_nonlinear} has been increased to be above that of \texttt{cons\_integral}. Additionally, the option \texttt{constraints/nonlinear/branching/mixfractional} has been added to enable considering nonlinear variables ahead of integer variables with fractional value for branching.
The new option specifies the minimal average pseudocost count across integer variables that is required to consider branching on a nonlinear variable before branching on a fractional integer variable.
By default the option is set to infinity, but when set to a small value, e.g., 0, 1, or 2, fractional integer and nonlinear variables are jointly
considered for branching by \texttt{cons\_nonlinear} as soon as
sufficiently many pseudocosts have been computed for integer variables.
For this purpose, the scoring of branching variables in \texttt{cons\_nonlinear} \citep[Section 4.2.12]{SCIP8} has been extended to consider the fractionality of an integer variable as a score similar to the violation score that is computed for a nonlinear variable.
If the scoring on this joint set picks a fractional integer variable however, then the regular branching rules for integer variables are employed.


\subsection{Benders' Decomposition}
\label{sec:benders}


The Benders' decomposition framework~\cite{maher2021benders} has received some small feature updates for \scip~\scipversion.
First is the introduction of objective types for the subproblems.
This aims to allow more flexibility to the user in formulating problems when using the Benders' decomposition framework.
Second is the improved identification and handling of master linking variables.
The latter feature is designed to enable the integer and no-good cuts on a broader range of problem types.
Finally, the decomposition process when supplying a DEC file has been updated so that the original solution is now returned from the Benders' decomposition solve.

\subsubsection{Objective type}

Classically, the application of Benders' decomposition results in the inclusion of an auxiliary variable to provide an underestimation of the subproblem objective value.
If the subproblem can be separated into disjoint problems, then this auxiliary variable can be substituted with the sum of auxiliary variables (one for each subproblem).
Prior to \scip~\scipversion, the Benders' decomposition framework only supported this classical handling of the subproblem objective.

While the summation of auxiliary variables is theoretically possible for all applications of Benders' decomposition, there are problems where an alternative objective may be beneficial.
An example is a multiple machine scheduling problem with a makespan objective.
In this case, an application of Benders' decomposition typically involves splitting the subproblems by machine.
Then in the master problem, a makespan variable is defined, which is then constrained to take the maximum of the makespan from each individual subproblem.
In such a setting, it is more convenient to define the subproblem objective as the minimum of the maximum makespan across all subproblem, compared to some summation objective.

In this release, the Benders' decomposition framework has been extended to support two different objective types: the classical summation and the minimum of the maximum subproblem auxiliary variables.
The objective type can be set using the method \texttt{SCIPsetBendersObjectiveType()}.
The objective type must be set during the problem creation stage, since this setting has an impact on the problem being solved.
Note that the different objective types only have an impact if more than one subproblem is used in the Benders' decomposition.

\subsubsection{Handling of Master Linking Variables}

The linking master variables are a critical part of the Benders' decomposition algorithm.
In the \scip implementation, the linking variables are identified by the existence of a pair of variables with the same name in the master and subproblem.
It is assumed that all master variables could potential be linking variables.
This assumption led to the incorrect disabling of the integer and no-good cuts.

Prior to \scip~\scipversion, the integer and no-good cuts would only be applied if the master problem was pure binary.
However, it is possible to apply these cuts when the master problem is an \MIP, provided that the linking variables are all binary.
This can treated in a more fine-grained manner, where the integer and no-good cuts can be applied to specific subproblems, because these only have binary linking variables.

To facilitate the better handling of integer and no-good cuts, the linking master variables are identified and stored during the initialization stage of the Benders' decomposition framework.
Statistics on the number of integer and binary linking variables are recorded for each subproblem.
These statistics are then used to determine whether the integer and no-good cuts can be enabled for the subproblems.

\subsubsection{Returning Original Problem Solution}

One method for applying Benders' decomposition using \scip is to supply a DEC file specifying the decomposition for an instance.
Additionally setting the parameters \param{decomposition/benderslabels} and \param{decomposition/applybenders} will decompose the instance and apply Benders' decomposition to solve the problem.
In previous versions of \scip, only the variable values for the master problem solution were directly available.
The subproblem solutions could be retrieved by calling \texttt{display subsolution <subproblem-index>} from the interactive shell or by calling \method{SCIPbendersSubproblem} (to retrieve the subproblem), \method{SCIPsetupSubproblem} (to apply the best known master problem solution), and then \method{SCIPsolveBendersSubproblem} (to solve the subproblem from scratch).
While this approach for retrieving the original problem solution is still necessary when an instance is supplied in its decomposed form, such as using the SMPS instance format or using a custom Benders' decomposition implementation, it was not practical when the original problem instance is known.

The application of Benders' decomposition when supplying a DEC file has been redeveloped in \scipv.
A relaxator has been added to the default plugins, which handles the decomposition and solving of an instance using Benders' decomposition.
The motivation for using a relaxator is so that the original instance remains untouched, while the decomposition can be applied within the relaxator.
The relaxator acts as a sandbox where it is possible to make the appropriate, destructive modifications to the problems in order execute the Benders' decomposition algorithm.

When the solving process is started in \scip, the original \scip instance will execute presolving.
During the processing of the root node from the original problem, the Benders' decomposition relaxator will be called first.
The Benders' decomposition algorithm attempts to solve the problem to optimality.
At the completion of the Benders' decomposition algorithm, the best found primal and dual bounds are returned to the original \scip instance.
If possible, the solution from the decomposed problem is mapped back to the original instance variables.

By default, the original \scip instance will terminate with an optimal solution or infeasible status.
If a limit has been set, a time, gap, or bound limit status could be returned, or a ``user interrupt'' indicated.
If the Benders' decomposition algorithm fails to solve the instance---due to reaching working limits---or the solution could not be returned, a ``user interrupt'' is triggered.
To continue solving the original \scip instance after the conclusion of the Benders' decomposition algorithm, parameter \texttt{relaxing/benders/continueorig} can be set to TRUE.


\subsection{IIS Finder}
\label{sec:iis}

It is a common issue for integer programming practitioners to (unexpectedly) encounter infeasible problems. Often it is desirable to better understand exactly why the problem is infeasible. Was it an error in the input data, was the underlying formulation incorrect, or was the model simply infeasible by construction? A common tool for helping diagnose the reason for infeasibility is an Irreducible Infeasible Subsystem (IIS) finder. An IIS is a subset of constraints and variable bounds from the original problem that when considered together remain infeasible and cannot be further reduced without the subsystem becoming feasible.

Practitioners can use IIS finders, since they narrow their focus onto a smaller, more manageable problem. Note, however, that there are potentially many different IISs for a given infeasible problem, and that IIS finders generally provide no guarantee of an IIS of minimum size. For a complete overview of IISs, see \cite{Chinneck2008}, and for an alternate reference detailing the implementation of an IIS finder, see \cite{Puranik2017}.

In \scip~\scipversion, users now have the functionality to generate an IIS for any infeasible CIP. This is achieved by calling \texttt{SCIPgenerateIIS} via the API or \texttt{iis} via the CLI, with the resultant IIS accessible afterwards by calling \texttt{SCIPgetIIS} via the API or both \texttt{write/iis} and \texttt{display/iis} via the CLI.  Moreover, users now have the ability to easily develop their own algorithms for generating an IIS, where we call these algorithms \textit{IIS finders}. These have been added as a new \scip plugin, aptly named \texttt{SCIP\_IISFINDER}, and therefore a user need only implement the required callbacks to write their own IIS algorithms.
Currently, \scip has implemented the additive and deletion based methods from \cite{Guieu1999}. The resultant IIS finder containing the additive and deletion based methods is called \textit{greedy}. Many parameters are also made available to the user, with two examples being: (1) Allowing only constraint deletion and leaving the original set of bounds as is. (2) Allowing (potentially faster) generation of an infeasible subsystem with no guarantee of irreducibility.


\subsection{CONOPT Interface}
\label{sec:conopt}


When solving MINLPs, SCIP makes use of local NLP solvers for several purposes. Primal heuristics solve NLPs, obtained by fixing integer variables, instead of LPs in order to find feasible solutions, convex NLP relaxations are used in cut generation and bound propagation, the Benders algorithm makes use of NLP subproblems, and some constraint handlers solve specialized NLPs.

This release introduces an interface to CONOPT (\url{conopt.com}), a nonlinear programming solver implementing a feasible path algorithm~\cite{parker1981approximation} that is based on active set methods~\cite{goldfarb1983numerically,gill1984weighted}. The solver is particularly suited for large, sparse models and settings where a series of mostly similar problems need to be solved, the solution of one problem providing a good starting point for the next.

Table~\ref{tbl:conopt} compares the performance of SCIP when using Ipopt and CONOPT as the nonlinear programming solver. The experiments were conducted on a subset of the MINLPLib instances which can be solved by SCIP~8.0.0 under one hour. The time limit was one hour.
The results show that using CONOPT instead of Ipopt increases the mean time by \SI{7}{\percent} and the mean number of nodes by \SI{4}{\percent}. However, on the subsets \bracket{100}{timelim} and \bracket{1000}{timelim} it reduces the mean time by \SI{4}{\percent} and \SI{40}{\percent}, respectively, and the mean number of nodes by \SI{1}{\percent} and \SI{26}{\percent}, respectively, and the overall number of solved instances is larger by 3 with CONOPT. We observed that the largest differences in performance occurred when one NLP solver was capable of solving an NLP problem and the other reached a time or iteration limit or arrived at a locally infeasible point.

\begin{table}
\caption{Performance comparison vs Ipopt}
\label{tbl:conopt}
\scriptsize

\begin{tabular*}{\textwidth}{@{}l@{\;\;\extracolsep{\fill}}rrrrrrrrr@{}}
\toprule
&           & \multicolumn{3}{c}{SCIP + CONOPT} & \multicolumn{3}{c}{SCIP + Ipopt} & \multicolumn{2}{c}{relative} \\
\cmidrule{3-5} \cmidrule{6-8} \cmidrule{9-10}
Subset                & instances &                                   solved &       time &        nodes &                                   solved &       time &        nodes &       time &        nodes \\
\midrule
\cleaninst            &       812 &                                      803 &       18.8 &         2526 &                                      800 &       18.2 &         2361 &        1.04 &          1.07 \\
\affected             &       523 &                                      520 &       18.2 &         3978 &                                      517 &       17.3 &         3587 &        1.06 &          1.11 \\
\cmidrule{1-10}
\bracket{10}{timelim}   &       426 &                                      423 &       81.5 &         8541 &                                      420 &       77.8 &         8132 &        1.05 &          1.05 \\
\bracket{100}{timelim}  &       187 &                                      184 &      338.3 &        43711 &                                      181 &      354.1 &        43976 &        0.96 &          0.99 \\
\bracket{1000}{timelim} &        49 &                                       46 &      703.3 &       140280 &                                       43 &     1167.9 &       190723 &        0.60 &          0.74 \\
\difftimeouts         &         9 &                                        6 &       31.4 &         2009 &                                        3 &     2718.6 &        42060 &        0.01 &          0.05 \\
\alloptimal           &       797 &                                      797 &       18.0 &         2419 &                                      797 &       16.4 &         2180 &        1.09 &          1.11 \\
\cmidrule{1-10}
continuous            &       184 &                                      180 &       11.2 &         4033 &                                      177 &       13.1 &         4460 &        0.86 &          0.90 \\
integer               &       623 &                                      618 &       21.8 &         2225 &                                      618 &       20.1 &         1995 &        1.08 &          1.12 \\
\bottomrule
\end{tabular*}
\end{table}


\subsection{Technical Improvements}
\label{sec:misc}


\subsubsection{Checking bounds on aggregated variables}

In some situations, \scip terminates with a solution that is feasible in the presolved problem, but has small infeasibilities (above the feasibility tolerance) in the original problem.
One possible source for such a situation are aggregated or multiaggregated variables.
For example, if a linear constraint $z=x+y$ with $x,y,z\geq 0$ is removed and $z$ is replaced by $x+y$ in the objective function and any other constraint, then a solution with $x=y=-10^{-6}$ can be feasible in the presolved problem (assuming default \texttt{numerics/feastol=1e-6}), but $z=-2\cdot 10^{-6}$ violates the lower bound of~$z$ by more than~$10^{-6}$.

Since \scip 9.1.0, a new constraint handler \texttt{cons\_fixedvar} is available that checks whether the bounds of aggregated variables are satisfied in a given solution.
If this is not the case, but the solutions needs to be enforced, then the aggregation is added as a cut to the LP relaxation.
In the given example, $z=x+y$ would be added.
If a cut could not be generated for some reason, then the feasibility tolerance for the LP relaxation is reduced for the current node.

\subsubsection{\scip Statistics Serialization}

\scipv introduces a new callback method for the table plugin used to write out statistics tables, called \texttt{TABLE\_COLLECT} that collects the data in the newly introduced \texttt{DATATREE} object that represents generic serialization of data.
This allows to generate the statistics in a generic way and to write them in different file formats.
This new callback method is implemented for all the statistics tables in \scip, and they can currently be exported to a JSON file using the \texttt{SCIPprintStatisticsJSON} function through the C API or
through the command line if the written statistics file has a JSON
extension.

\subsubsection{Writing AMPL \texttt{.nl} files}

The reader for AMPL \texttt{.nl} files has been extended by writing capabilities.
It utilizes the NL writer of the \texttt{ampl/mp} project and currently supports general and specialized linear constraints and nonlinear constraints.



\section{The \gcg Decomposition Solver}
\label{sect:gcg}


\noindent
\gcg extends \scip by providing the functionality of a
decomposition-based solver. That is, \gcg reformulates a given \MILP
according to an automatically detected Benders or Dantzig-Wolfe
decomposition, and then solves the reformulation by branch-and-cut or
branch-cut-and-price, respectively. For the user, no expertise in
these techniques is needed, but for an understanding of the following,
we assume that the reader is familiar with decomposition methods.
We refer to~\cite{desrosiers2024branch} for a recent reference.

The central design paradigm of \gcg for its main use case,
branch-cut-and-price based on a Dantzig-Wolfe reformulation, is the
synchronized work on two \scip instances, reflecting the original
(user-given) and the master (reformulated) models. This design allows
for benefiting from \scip's functionality on both models, which are
intimately linked by theory, and thus both needed for a modern
branch-cut-and-price implementation.  Most decisions and information
available in one tree, like branching, cutting, dual bounds, primal
solutions, propagations, etc.\ are mirrored to the other tree. 
The (relaxed) Dantzig-Wolfe master model lives in a relaxator of the
original \scip instance and is solved by column generation in every
node.

\subsection{License}

\gcg has been open-source since its initial release.
In order to harmonize the licensing with \scip, all contributing code
authors have approved to switch from the GNU
LGPL to the Apache~2.0 license, starting with \gcgv.

\subsection{\gcg Data Structure}
\gcgv introduces a \gcg data structure similar to \scip's main data structure.
The main purpose is to simplify the C API as many functions now accept a \gcg object instead of a SCIP object.
Users do no longer need to check whether they have to pass the original or a master problem to the function.
The \gcg object stores important and frequently accessed pointers.
Among these are pointers to the original problem, the master problems (Dantzig-Wolfe and Benders), the relaxator, and the pricer.
This also eliminates the need of storing the same pointers in many user data objects or the need of \texttt{SCIPfind*} calls.
A \gcg object can be created using \texttt{GCGcreate} and freed using \texttt{GCGfree}.
We plan to extend this data structure in future releases.

\subsection{Extended Master Constraints}

The Dantzig-Wolfe reformulation of an original model ${\min \{c^\T x :
Ax \geq b, Dx \geq d, x \in \Z^{n}\}}$  works with $X\defi\conv\{x \in
\Z^{n} : Dx \geq d\}$. Denote the finite sets of integer extreme (and
potentially some integer interior) points and integer extreme rays of
$X$ by $\{x_p\}_{p\in P}$ and~$\{x_r\}_{r\in R}$, respectively. The
Dantzig-Wolfe (discretized) integer master problem
\begin{equation}
  \label{eq:sec_gcg_master}
  \begin{aligned}
    \min  \quad&  \sum_{p\in P} (c^\T x_p)\lambda_p + \sum_{r\in R} (c^\T
    x_r)\lambda_r \\ 
    \text{s.t.} \quad&   \sum_{p\in P} (A x_p) \lambda_p + \sum_{r\in R} (A x_r)\lambda_r
    \geq b &\qquad[\pi]\\
    & \sum_{p\in P} \lambda_p = 1 &\qquad[\pi_0]\\
    & \lambda \in \Z_+^{|P|+|R|}
  \end{aligned}
\end{equation}
is (one option for) an equivalent reformulation of the original model,
using (among other conditions) the linking
\begin{equation}
  \label{eq:sec_gcg_link_x_and_lambda}
  x= \sum_{p\in P} x_p \lambda_p + \sum_{r\in R} x_r \lambda_p
\end{equation}
between the original $x$ and master $\lambda$ variables. 
The $\pi, \pi_0$ in brackets denote the dual variables corresponding
to the constraints in the LP relaxation. This
relaxation is potentially stronger than that of the original problem.
To solve the problem by column generation, it is necessary to repeatedly solve the 
pricing problem ${\min \{(c^\T -\pi^\T A) x -\pi_0: x \in X\}}$, where the 
objective function
expresses the reduced costs of a master variable.

In~\eqref{eq:sec_gcg_link_x_and_lambda} it can be seen that integrality of
$\lambda$-variables implies integrality of $x$-variables, but the
converse does not hold.
Branching (and cutting) can be carried out on
both kinds of variables. A branching decision on a fractional original
$x$-variable translates via~\eqref{eq:sec_gcg_link_x_and_lambda} to a
constraint in the (relaxation of the) master
problem~\eqref{eq:sec_gcg_master}. Even easier, the branching decision
can be enforced by a variable bound in the pricing problem. Therefore,
tailored pricing algorithms are usually still applicable after
branching. This branching (and cutting) on original variables,
sometimes called \emph{robust}, is part of \gcg.

In certain situations, in particular when the master problem is
\emph{aggregated}, users want to 
branch on fractional sums of master variables. Branching constraints, 
that are added to the master problem, are of the following form:
\begin{equation}
  \label {eq:sec_gcg_master_branching}
\sum_{p \in P} f(p) \lambda_p + \sum_{r \in R} f(r) \lambda_r \geq h, \qquad \left[\gamma\right]
\end{equation}
where, again, $\gamma$ denotes the corresponding dual variable in
the master relaxation. The coefficients $f(j)$, $j\in P\cup R$, are
usually binary and represent a well-defined subset of master
variables. The Ryan-Foster branching~\cite{gcg:RyanFoster:81} on two
rows $r$ and $s$ is a classical example where $f(j)=1 \iff
x_{rj}=x_{sj}=1$, $j\in P$, that is, both components $r$ and~$s$ of vector $x_j$ are nonzero.
Ryan-Foster is part of \gcg, but many more
options are conceivable, see
Section~\ref{sec:gcg_component_bound_branching} for an example.

Note that cutting planes in the master variables, like clique cuts
or Chv\`atal-Gomory cuts, are also of the 
form~\eqref{eq:sec_gcg_master_branching}.
%

The coefficient $f(j)$ needs to be computed in the pricing problem as
well.  Unlike in robust branching, $f(j)$ \emph{cannot} be stated as linear
expressions in the original $x$-variables already present in the
pricing problem \emph{alone}: Instead, we introduce a \emph{coefficient 
variable}~$y$, one for every branching constraint or cutting plane. This 
variable
is constrained to a domain $Y$ that ensures the correct correspondence of
the coefficient $y=f(j)$ and the semantics of the branching decision
or cutting plane. The set $Y$ is defined by linear 
constraints in $x$, $y$, and
potentially auxiliary (integer) variables $z$.
The $z$-variables are called \emph{inferred} variables in the code.
The necessary modifications in the pricing problem 
\begin{displaymath}
\begin{aligned}
\min\quad & \left( c^\T - \pi^\T A \right) x - \gamma y - \pi_0 & \\
\text{s.t.}\quad 
 &x  \in X \\
 & y   \in Y(x,z)
\end{aligned}
\end{displaymath}
may entail a nontrivial additional computational burden and/or
destroy the applicability of a tailored pricing algorithm. Such
branching rules or cutting planes are therefore also called \emph{non-robust}.
They can be stronger than robust branching and cutting.
%

Concluding, in nonrobust branching or cutting the master
constraint~\eqref{eq:sec_gcg_master_branching} 
\emph{does not} result from a
Dantzig-Wolfe
reformulation of a constraint in the original model, in contrast to
the robust case.
In \gcgv we therefore introduce the concept of \textit{extended master
  constraints} which live \emph{only} in the master problem and therefore do 
  not correspond to an 
original constraint.
We provide a new interface for creating and managing extended master
constraints to use within 
branching rules and separators.
This interface is deeply integrated into \gcg.  The modifications to
the pricing problem are automatically and dynamically applied and
undone in each node as necessary.  Dual value stabilization of
extended master constraints is supported.  When users develop their
own master branching rules, it suffices to create extended master
constraints using the new interface in \gcgv.
Moreover, basic functionality for separators is available.
However, this is an experimental feature that is still under development.

\subsubsection{Component Bound Branching Rule}
\label{sec:gcg_component_bound_branching}

When branching on master variables, one needs to define a function
$f(j)$, $j\in P\cup R$, that guarantees that the left-hand side
of~\eqref{eq:sec_gcg_master_branching} is fractional whenever the
solution to the restricted master problem is fractional.
Typically, $f(j)$ just identifies a subset of master variables.
The \emph{component bound branching} rule~\cite{gcg:Vanderbeck:00} is
a generic example for such a function.
In essence, the rule branches on sums of master variables $\lambda_j$
for which the components of the corresponding $x_j$, $j\in P\cup R$,
satisfy given lower and upper bounds, see~\cite{desrosiers2024branch} for details.

This rule is implemented in \gcgv as a code example for how to use the
newly introduced extended master constraints.
The documentation contains a section on ``how to add master-only
constraints'' with details about the interface methods that need to be
implemented.

\subsection{JDEC Decomposition File Format}
\gcg supports various decomposition file formats.
Users can provide \gcg with their own predefined decompositions or let \gcg write detected ones to files.
The \texttt{DEC} file format is easy to use as \texttt{.dec} files are simple text files.
However, this makes adding additional and advanced features difficult.
\gcgv supports a new JSON-based decomposition file format, called \texttt{JDEC} using the file extension \texttt{.jdec}.
JSON is a straightforward and widely adopted file format that is human-readable.
Data is stored using arrays and key-value pairs.
Many modern programming languages support reading and writing JSON files either natively or by using an open-source library.
Hence, a great advantage is that the basic structure of \texttt{.jdec} files can be read in without any specialized parser.
Currently, \texttt{.jdec} files support the following main features, which we plan to extend in the 
future:
\begin{itemize}
    \item partial and complete decompositions,
    \item declaring constraints as master or block constraints,
    \item specifying nested structures, i.e., decompositions can be provided 
    for blocks,
    \item storing of symmetry information.
\end{itemize}

The root object of a \texttt{.jdec} file stores metadata of the decomposition and supports the following data fields:
\begin{table}[H]
    \centering
    \begin{tabularx}{\textwidth}{@{}lX@{}}
        \toprule
        name & description \\
        \midrule
        \texttt{version}           & version of the file format \\
        \texttt{problem\_name}      & name of the problem the decomposition belongs to \\
        \texttt{description}       & description of the decomposition \\
        \texttt{decomposition}     & the decomposition (see below) \\
        \texttt{decomposition\_id} & internal ID of the decomposition, only written by GCG \\
        \bottomrule
      \end{tabularx}
\end{table}
The root object stores the actual (root) decomposition using a decomposition object.
Decomposition objects always have the same data fields and can be used to specify decompositions of blocks as well.
A decomposition object of a \texttt{.jdec} file supports the following data fields:
\begin{table}[H]
    \centering
    \begin{tabularx}{\textwidth}{@{}lX@{}}
      \toprule
      name & description \\
      \midrule
      \texttt{n\_blocks}             & number of blocks of the decomposition, ignored by GCG \\
      \texttt{presolved}            & indicates whether the decomposition refers to a presolved problem \\
      \texttt{master\_constraints}   & list of master constraints \\
      \texttt{blocks}               & list of block objects (see below) \\
      \texttt{symmetry\_var\_mapping} & a mapping that maps all variables to their representatives \\
      \bottomrule
    \end{tabularx}
\end{table}
A block object of a \texttt{.jdec} file is used to store information about a block.
As already mentioned, users may provide a decomposition for blocks, leading to nested structures.
A block object supports the following data fields:
\begin{table}[H]
    \centering
    \begin{tabularx}{\textwidth}{@{}lX@{}}
      \toprule
      name & description \\
      \midrule
      \texttt{index}             & index of the block \\
      \texttt{constraints}            & list of constraints assigned to the block \\
      \texttt{decomposition}   & decomposition (object) of the block \\
      \texttt{symmetry\_representative\_block}               & index of the representative block this block should be mapped to \\
      \bottomrule
    \end{tabularx}
\end{table}

When reading a \texttt{.jdec} file GCG ignores unknown data fields.
Thus, users can add additional information, that they need to work with, to 
their custom decompositions. 
Similarly, we plan to extend the corresponding reader plugin such that \gcg includes additional information about the detection 
process in the written file.
\cref{lst:gcg:jdecfile} shows an example \texttt{.jdec} file that stores a 
nested decomposition for the \texttt{.lp} file depicted in 
\cref{lst:gcg:lpfile}.
For more details we refer to \gcg's documentation\footnote{\url{gcg.or.rwth-aachen.de}}.

\begin{lstlisting}[
    caption={truncated \texttt{.lp} file example},
    label={lst:gcg:lpfile},
    basicstyle=\scriptsize,
]
Minimize
  obj:  + y#0 + y#1 + y#2 + y#3
Subject to
  assign_0:
   + x#0#0 + x#1#0 + x#2#0 + x#3#0 >= 1
  assign_1:
   + x#0#1 + x#1#1 + x#2#1 + x#3#1 >= 1
  link_0:
   - y#0 + y#1 >= 0
  link_1:
   - y#2 + y#3 >= 0
  cap_0:
   -10 y#0 + 2 x#0#0 + 3 x#0#1 <= 0
  cap_1:
   -5 y#1 + 2 x#1#0 + 3 x#1#1 <= 0
  cap_2:
   -10 y#2 + 2 x#2#0 + 3 x#2#1 <= 0
  cap_3:
   -5 y#3 + 2 x#3#0 + 3 x#3#1 <= 0
\end{lstlisting}

\begin{lstlisting}[
    caption={\texttt{.jdec} file example},
    label={lst:gcg:jdecfile},
    basicstyle=\scriptsize,
]
{
  "version": 1,
  "name": "example_nested_dec",
  "problem_name": "example.lp",
  "description": "nested decomposition with aggregated blocks",
  "decomposition": {
    "presolved": false,
    "n_blocks": 2,
    "master_constraints": [
      "assign_0"
    ],
    "blocks": [
      {
        "index": 0,
        "constraints": [
          "cap_0",
          "cap_1",
          "link_0"
        ],
        "decomposition": {
          "presolved": false,
          "n_blocks": 2,
          "master_constraints": [
            "link_0"
          ],
          "blocks": [
            {
              "index": 0,
              "constraints": [
                "cap_0"
              ]
            },
            {
              "index": 1,
              "constraints": [
                "cap_1"
              ]
            }
          ]
        },
        "symmetry_representative_block": 0
      },
      {
        "index": 1,
        "constraints": [
          "cap_2",
          "cap_3",
          "link_1"
        ],
        "decomposition": {
          "presolved": false,
          "n_blocks": 2,
          "master_constraints": [
            "link_1"
          ],
          "blocks": [
            {
              "index": 0,
              "constraints": [
                "cap_2"
              ]
            },
            {
              "index": 1,
              "constraints": [
                "cap_3"
              ]
            }
          ]
        },
        "symmetry_representative_block": 0,
      }
    ],
    "symmetry_var_mapping": {
      "y_0": "y_0",
      "x_0_0": "x_0_0",
      "x_0_1": "x_0_1",
      "y_1": "y_0",
      "x_1_0": "x_0_0",
      "x_1_1": "x_0_1",
      "y_2": "y_0",
      "x_2_0": "x_0_0",
      "x_2_1": "x_0_1",
      "y_3": "y_0",
      "x_3_0": "x_0_0",
      "x_3_1": "x_0_1"
    }
  }
}
\end{lstlisting}

\subsection{Pricing Problem Solvers}
This release adds two new pricing problem solvers: the \gcg pricing problem solver and the HiGHS pricing problem solver.
In addition, we revised the Cliquer pricing solver.

\subsubsection{\gcg Pricing Problem Solver}
The \gcg pricing problem solver solves general \MILP pricing problems using \gcg.
The common use case is to specify a nested structure using the new \texttt{JDEC} file format.
Then, if enabled, the \gcg pricing solver uses this structure to reformulate and solve the pricing problems.
Currently, only Dantzig-Wolfe reformulations are supported.
The solver can be enabled by setting the parameter 
\texttt{pricingsolver/gcg/maxdepth}, which refers to the maximal 
recursion depth, to a value greater than zero.
Until this depth is reached, the solver either uses the structure information provided by the current (parent) decomposition or tries to detect a structure.
The detection takes only place if no nested structure information is provided 
for all blocks of the parent decomposition, and the maximum depth is not 
reached yet.

\subsubsection{HiGHS Pricing Problem Solver}
\label{sec:gcg:highspricingsolver}

This pricing solver uses the HiGHS \MILP solver \cite{HuangfuHall15} to solve the pricing problems.
The use of the HiGHS pricing solver is an optional feature and has to be enabled at compile time.
It is necessary to have a separate installation of HiGHS\footnote{\url{https://github.com/ERGO-Code/HiGHS}}, which is open-source.
If enabled, it will be used to solve general \MILP pricing problems instead of 
the default \MILP pricing problem solver using \scip.

\subsubsection{Cliquer Pricing Problem Solver}
The Cliquer pricing problem solver is a specialized heuristic pricing solver.
It is applicable to weighted independent set 
pricing problems and closely related variants.
The solver uses the Cliquer library\footnote{\url{https://users.aalto.fi/~pat/cliquer.html}} that provides functionality to find maximum weighted cliques.
We revised large parts of the implementation and improved variable aggregation, the handling of fixed variables, and the check of pricing problem compatibility.
Moreover, we added more parameters that allow to control whether the solver should run or reject the solving request.
These parameters are based on properties of the constructed graph, such as 
number of nodes or density.

\subsection{Pricing Parallelization}
In the past, \gcg's parallelization was disabled by default at compile time.
With \gcgv, we have revised the implementation of \gcg's pricing parallelization.
If supported by the build environment, this feature is now enabled by default.
The parameter \texttt{pricing/masterpricer/nthreads} controls
the maximal number of threads \gcg can use for parallelization.
Setting the parameter to $0$ means that \gcg is allowed to use all available logical cores.
At runtime, \gcg informs the user about 
how many threads are used.
Note that using more than one thread may lead to nondeterministic behavior, 
i.e., the solving path for the same input may differ.
Thus, by default, \gcg uses only one thread.
We intend to add a deterministic mode for pricing parallelization in the future.

\subsection{IPColGen Primal Heuristic}
\label{sec:gcg:ipcolgen}
IPColGen is a matheuristic proposed by \citet{maher2023ipcolgen} and is implemented in \gcgv.
It is based on an LNS heuristic framework and uses an adapted pricing scheme.
The employed column generation strategy aims at finding columns that lead to 
high-quality primal feasible solutions for the original problem.
The heuristic is applicable if the master problem is a set covering, set packing, or set partitioning problem.
For the implementation of IPColGen, \gcg's API and data structures were extended.
Since IPColGen adapts the pricing scheme, new callback methods are introduced with \gcgv.
These allow users to influence \gcg's pricing procedure and are called before 
or after \gcg's pricing loop.
The callback methods are used in the IPColGen heuristic to inject custom 
weights used by the subsequent pricing procedure.
For more details of the IPColGen heuristic we refer to \citep{maher2023ipcolgen}.
The heuristic is currently 
disabled by default, but can be enabled by changing the parameter 
\texttt{heuristics/ipcolgen/freq} in the master problem.

\subsection{Decomposition Score Plugin}
\gcg computes scores to compare decompositions.
Based on the resulting ranking a decomposition is picked for the reformulation of the original problem.
The implementation of the scores has been refactored\footnote{this happened already in \gcg version 3.6, but was not documented}.
Scores are now implemented as plugins using the data structure \texttt{GCG\_SCORE} and can be included by calling \texttt{GCGincludeScore}.
This simplifies the addition and handling of scores.
When working with the CLI, users can display implemented scores with the \texttt{display scores} command.
The currently selected score can be changed by setting the parameter \texttt{detection/scores/selected}.




\section{SCIP-SDP}
\label{sect:scip-sdp}

\scipsdp is a solver for handling mixed-integer semidefinite programs
(MISDPs), without loss of generality, written in the following form
\begin{equation}\label{MISDP}
  \begin{aligned}
    \inf \quad & b^\T y \\
    \text{s.t.} \;\;\, & \sum_{k=1}^m A^k\, y_k - A^0 \succeq 0, \\
    & \ell_i \leq y_i \leq u_i && \text{for all } i \in [m], \\
    & y_i \in \Z && \text{for all } i \in \mathcal{I},
  \end{aligned}
\end{equation}
with symmetric matrices $A^k \in \R^{n \times n}$ for
$i \in \{0, \dots, m\}$, $b \in \R^m$, $\ell_i \in \R \cup \{- \infty\}$,
$u_i \in \R \cup \{\infty\}$ for all $i \in [m] \defi \{1, \dots, m\}$. The
set of indices of integer variables is given by $\mathcal{I} \subseteq [m]$, and
$M \succeq 0$ denotes that a matrix $M$ is positive semidefinite. The
development of \scipsdp has been described in the earlier SCIP reports.

The current version of \scipsdp is 4.4.0, which incorporates updates for
SCIP 10. Moreover, it features the possibility to write the original and
transformed problem in the \emph{Conic Benchmark Format} (CBF), which has
been developed in conjunction with the Conic Benchmark Library, see
\url{https://cblib.zib.de/}. In this way, also presolved problems can be
written.

In addition, we mention the Python code of Johanna Skåntorp at
\url{https://github.com/J-Skantorp/CBF-pythonic}, which allows to formulate
MISDPs and write them in CBF, which then can be read by \scipsdp.


\section{\papilo}
\label{sect:papilo}


The solver-independent presolving library PaPILO~\cite{GleixnerGottwaldHoen23} for MIP and LP has been shipped with the \scipopt since \scip~7.0.
It hooks into the solving process of \scip through the presolver plugin \texttt{milp}.
As the only presolving method in \scip, it also supports multi-precision and rational arithmetic and is, therefore, an important plugin when \scip is run in exact solving mode, see \Cref{sect:exactscip}.
However, as a solver-independent library, \papilo is also used outside of the \scipopt.
Examples are its use as a presolver for the first-order linear programming method \solver{PDLP}~\cite{PDLP}, and notably for presolving in the PB solver \solver{mixed-bag}, which won the opt-lin category in the Pseudo-Boolean Competition 24~\cite{PB24}.
In the same competition, \papilo also participated as part of the \scip entry, reaching first or second places in linear and nonlinear categories.

The latest version \papilov comes with significant improvements in performance and numerical stability for floating-point arithmetic.
In the following, we describe the improvements of performance and memory consumption of the \emph{dominated columns} presolver and introduce the newly included presolver \emph{clique merging}.

The presolvers in \papilo are parallelized using a transaction-based design~\cite{GleixnerGottwaldHoen23}, in which changes are first recorded as so-called transactions and applied later in a sequential synchronization phase.
This introduces challenges, requiring developers to ensure that transactions
\begin{enumerate}
    \item are stored in a reasonable order (effectiveness) and
    \item are generated without excessive redundancy (efficiency),
\end{enumerate}
especially to avoid conflicts among reductions of the same presolving module.
Redundant transactions may be generated, for example, by the dominated columns presolver~\cite{AchterbergBixbyGuetal.2019} due to the transitivity of dominations.
This means, dominations $x_1 \xrightarrow{\text{dom}} x_2$ and $x_2 \xrightarrow{\text{dom}} x_3$ imply the domination $x_1 \xrightarrow{\text{dom}}x_3$ but not all these dominations are required to fix the dominating columns $x_2$ and $x_3$.
Previous versions of \papilo store all of the detected dominations without regarding potential redundancies.
This can cause memory issues and substantial conflicts especially on feasibility problems, on which usually more redundant dominations exist than on problems with a nonzero objective.
In \papilov, we apply a \emph{topological compression} to the set of column domination arcs in order to avoid storing redundant dominations as explained next.
This keeps the required memory linear in the problem size and substantially reduces the computation time, while maintaining the same number of applied reductions in practice.

One requirement of the transaction-based design is that the bounds of a dominated column must be locked so that applying the domination by fixing the dominating column remains a valid reduction.
Due to the lock, a dominating column can only be fixed once.
Therefore, we aim at a directed forest on the column set of applicable dominations because in this structure no conflicts can arise.
To keep the detection of dominations parallelized and deterministic, the detected dominations are locally collected in thread-separate buckets.
The detection is continued until the number of collected dominations reaches the number of columns because then there exists a redundant domination as explained above.
Subsequently, this set of dominations is compressed by adding the dominations one by one to a directed forest structure while filtering out dominations that would lead to a directed cycle.
The detection and compression is repeated until all possible dominated columns are considered.
Afterwards, the resulting set of dominations is converted into a transaction sequence without conflicts by iterating through the domination trees in topological order by a breadth-first search.

The following examples highlight the resulting performance improvements:
\begin{itemize}
    \item On the instance \solver{neos-94313}~\cite{MIPLIB} the number of transactions decreases from \num{223200}, of which only \SI{6.2}{\percent} are applied, to now \num{13950} transactions, which are all applied.
    \item On the instance \solver{normalized-aries-da\_network\_1000\_5\_\_369\_766\_\_256}~\cite{PB9} the number of transactions decreases from \num{163200000}, of which only \SI{0.8}{\percent} are applied, to now \num{1275000} transactions, which are all applied.
    \item On the instance series \solver{normalized-PB09/OPT-LIN/aries-da\_nrp/normalized-aries-da\_network\_1000\_5\_\_369\_766}~\cite{PB9} for indices 512, 1024, and 2048 with up to \num{10240998} variables, presolver DominatedColumns in \papilo exceeded a memory limit of \SI{48}{\giga\byte} while \scip with DominatedColumns in \papilo disabled was not able to solve them within reasonable time. With these changes even such large-scale instances are solved by default to optimality within the given memory limit.
\end{itemize}
In all cases, the number of transactions is significantly reduced, and the resulting transactions can be applied without conflicts, thereby drastically reducing memory demand.

The new clique merging presolver is implemented similarly to the approach described by Achterberg in~\cite{Achterberg2022}.
Cliques are constraints which consist of binary variables where at most one variable can be assigned the value one.
In the presolver, we first construct a so-called \emph{clique
  table}~\cite{Achterberg2007a}, which is a graph structure, 
where each binary variable is represented by a vertex, and an edge connects each pair of variables that belong to the same clique.
Then a greedy maximum clique heuristic expands cliques that are already represented by other cliques.
Extending cliques can make some cliques redundant, which are then labeled as such and are eventually removed.

In order to limit the memory consumption, the algorithm is restricted to cliques with a certain size.
Since clique merging in \papilov is executed in parallel, it is considerably faster than the pre-existing clique merging presolver in \scip.
The current implementation is limited by the fact that extended cliques are copied back to the data structure of the reduced problem.
However, the current data structure in \papilo does not support large extensions of constraints, leading to the necessity of discarding some reductions. 
This issue could be addressed in the future by maintaining a separate clique table.



\section{SoPlex}
\label{sect:soplex}


\soplex~\cite{Wunderling1996} is a simplex-based open-source LP solver and serves as the default underlying solver for LP relaxations within \scip.
In particular, since \soplex supports solving LPs exactly over the rational numbers, it plays an important role in the exact solving mode of \scip, cf.\ Section~\ref{sec:exact}.
Although no algorithmic improvements have been added since the release of \scipopt~9.0, a new major version \soplexv is now released to signify improvements of the build system and needed adjustments of the public API.
For details, see the \soplex changelog.


\section{The UG Framework}
\label{sect:ug}


This release includes a UG application for solving Pseudo-Boolean problems in parallel with FiberSCIP~\cite{Mexi2025}.
A boolean parameter \texttt{PBCompetionOutputFormat} was added to adjust the output to that of the 2024 Pseudo-Boolean Competition\footnote{\url{https://www.cril.univ-artois.fr/PB24/}}
(default value is ``FALSE''). 

\section{ZIMPL}
\label{sect:zimpl}


With this release we updated \zimpl to Version~3.7.0. Apart from some small bug fixes,
\zimpl can now compute permutations. Given a set \texttt{A}, \texttt{permutate(A)} generates a set consisting of tuples with all
permutations of the elements of \texttt{A}. This simplifies modeling for problems where one needs
all the permutations, e.g., Birkhoff decompositions.
Also, variables declared as implicit integer are now recognized by \scip and handled appropriately. 


\section{Interfaces}
\label{sect:interfaces}

\subsection{PySCIPOpt}
\label{subsect:pyscipopt}
The main new feature of PySCIPOpt~\cite{pyscipopt}, SCIP's Python interface, is the support of matrix variables.
In some optimization subfields and other adjacent fields, it is standard practice to both view and model the optimization problem from a matrix perspective. 
The previous requirement of using individual variables often created a disconnect for such practitioners between the model and its implementation.
With this in mind, the concept of matrix variables has been added to PySCIPOpt. The matrix variable structure is built upon the popular \texttt{numpy.ndarray}, and thus \href{https://numpy.org/doc/stable/index.html}{NumPy} has become a required dependency.
Matrix variables act in a similar manner to standard variables, and can be used in an expected fashion to create intermediate expressions involving other variables, as well as for creating both linear and general nonlinear constraints. 
The syntax resembles that of single variables, using the methods \texttt{addMatrixVar} instead of \texttt{addVar} and \texttt{addMatrixCons} instead of \texttt{addCons}.
This allows for an easy transition for users familiar with the existing interface.
For tutorials on how to use this functionality, visit the \href{https://pyscipopt.readthedocs.io/en/latest/tutorials/matrix.html}{Matrix API} tutorial in the documentation.
 
In addition to matrix variables, work was put into helping new users get started with the software.
This was achieved by three main means: improving the documentation, creating exercises, and including recipes.
The documentation~\cite{pyscipoptDocu2025} was reworked to have a more aesthetic style, and to include many new tutorials and examples.
In the same vein, scipdex~\cite{scipdex} was created, with the goal of offering a more interactive
way for users to familiarize themselves with PySCIPOpt. 
It contains exercises of varying difficulty, facilitating the exploration
of different aspects of SCIP.

Further, in an attempt to provide additional functionality to users, the recipes sub-package was created. 
It consists of common features that newcomers often have trouble implementing. 
Two such examples are (1) requiring an epigraph reformulation to add a nonlinear objective and (2) plotting the evolution of the primal and dual bounds, for example.
These were not added directly to PySCIPOpt's source code as to reduce divergence from SCIP.

\subsection{SCIP.jl}
\label{subsect:scipjl}

\scip is available in {\sc Julia} through the SCIP.jl package, which exposes both a low-level interface mirroring exactly the C API and a high-level interface based on {\sc MathOptInterface}~\cite{legat2022mathoptinterface}.
The high-level interface now exposes the creation of \scip event handlers, allowing users to customize the execution of the solving process with more flexibility.
The high-level interface also includes the possibility to compute a set of Minimum Unsatisfiable Constraints (MinUC), i.e., a set of constraints which cannot be satisfied in the current problem.
The SCIP.jl implementation of MinUC computes the set of unsatisfiable
constraints and then allows users to query the status of each constraint to
determine whether they are in the conflict, see, e.g., \cite{Chinneck2008,Pfe08}.



\subsection{JSCIPOpt}

The Java interface to \scip, \solver{JSCIPOpt}~\cite{javagithub}, has been extended to make the following functionality of the {\scip} C API newly available from Java:
\begin{itemize}
\item the \texttt{SCIPgetStage} getter method (contributed by the GitHub user fuookami (Sakurakouji Sakuya)),
\item the getters \texttt{SCIPgetBoolParam}, \texttt{SCIPgetIntParam}, \texttt{SCIPgetLongintParam},\\ \texttt{SCIPgetRealParam}, \texttt{SCIPgetCharParam}, \texttt{SCIPgetStringParam},\\ \texttt{SCIPgetPrimalbound}, \texttt{SCIPgetSolvingTime}, \texttt{SCIPgetNOrigVars},\\ and \texttt{SCIPgetOrigVars},
\item algebraic expressions (\texttt{SCIPcreateExpr*} and \texttt{SCIPreleaseExpr} methods),
\item all built-in constraint types (\texttt{SCIPcreateConsBasic*} methods), e.g., nonlinear constraints (using the algebraic expressions), SOS constraints, logical constraints, etc.
\end{itemize}

In addition, experimental \solver{JSCIPOpt} support for writing custom expression handlers in Java, using the \texttt{ObjExprHdlr} interface added to the \texttt{objscip} C++ API in SCIP 10.0, is available in a branch.

Please note that support for SCIP versions prior to 8.0, which was already broken in practice, was officially dropped from \solver{JSCIPOpt}, and the minimum CMake version increased to 3.3, which was the minimum required by SCIP 8.0. This allowed fixing some deprecation warnings.

\subsection{russcip}
\label{subsect:rust}
At the time of this release, the Rust interface to \scip, \texttt{russcip}, is at Version $0.8.2$ bringing many quality-of-life improvements since the previous release.
It introduces the new \texttt{bundled} feature, which allows users to access a precompiled version of \scip and its dependencies on all major platforms.
Additionally, the \texttt{russcip} library has been updated to provide a more ergonomic API, providing a more expressive way to build models.
The new API has default values for variable and constraint data, and only the needed ones are passed.
This also ensures correct types are used at compile time, e.g., an integer variable can only be created with integer bounds.
Next is an example of the new API in comparison to the old one.

\lstdefinelanguage{Rust}{
    morekeywords={abstract,alignof,as,become,box,break,const,continue,crate,do,else,enum,extern,false,final,fn,for,if,impl,in,let,loop,macro,match,mod,move,mut,offsetof,override,priv,proc,pub,pure,ref,return,self,Self,static,struct,super,trait,true,type,typeof,unsafe,unsized,use,virtual,where,while,yield},
    sensitive=true,
    morecomment=[l]{//},
    morecomment=[s]{/*}{*/},
    morestring=[b]{"},
}

\lstset{
    language=Rust,
    basicstyle=\ttfamily\footnotesize,
    keywordstyle=\color{blue},
    commentstyle=\color{gray},
    stringstyle=\color{red},
    showstringspaces=false,
    breaklines=true,
    tabsize=4
}

\begin{lstlisting}
use russcip::prelude::*;

let mut model = Model::default().minimize();

let x = model.add_var(0.0, 1.0, 1.0, "x", VarType::Binary); // previous API
let y = model.add(var().bin().obj(1.0)); // add binary variable with objective coefficient 1.0

model.add(cons().coef(&x, 1.0).coef(&y, 1.0).eq(1.0)); // x + y = 1
model.add_cons(vec![&x, &y], &[1.0, 1.0], 1.0, 1.0, "cons"); // previous API
\end{lstlisting}

The safe rust API has also grown significantly, providing access to the separator and constraint handler plugins, and methods to build and add cuts.
Access to a \texttt{Model} object in the solving stage is now only available through the plugin callbacks (trait methods), which enables separation of methods only available in the solving stage to only be called in plugin implementations.
Finally, the release brings a simpler memory management model.
All \scip wrapper objects now contain a reference-counted pointer to the \scip instance they are generated from ensuring that they are valid for the lifetime of the instance.


\subsection{Further Interfaces}

Next to the previously mentioned interfaces, there also exist interfaces of
\scip for Matlab~\cite{matlabgithub}, AMPL, and the C++ interface \solver{SCIPpp}~\cite{scipppgithub}.
The LP solver \soplex features the Python interface \solver{PySoPlex}~\cite{pysoplexgithub}.
The presolver \papilo now also has a file-based interface in Julia for presolving and postsolving instances~\cite{papilojl}.

\section{\mipdd}
\label{sect:mipdd}


\mipdd is a C++-package that helps to simplify debugging optimization software by reducing the instance triggering a bug. \mipdd represents the first open-source and solver-independent \emph{delta debugger} for mixed-integer programming solvers~\citep{HKG24MIPDD}.

Debugging optimization solvers, such as \scip, can be quite challenging with many modules interacting in a complex way.
A general technique to track down bugs in any kind of software is the implementation of assertions, which are boolean expressions expected to be true signaling an error if realized to be false.
Another approach specific to optimization software is the utilization of a so-called debug solution mechanism, which for a given feasible solution signals an error if the solver withdraws it as optimal solution candidate.
Depending on how familiar a developer is with the relevant code, identifying the reason of an error on large instances can be difficult and time-consuming.
Furthermore, the above debugging methods can only be applied if the failing instance is available.
Due to legal reasons or concerns about sensitive information, users may be unable to share an instance with the development team.
Therefore, a benefit of an automatic simplification method is that it can destroy relevant sensitive information if it is irrelevant to reproduce an issue, allowing external developers to investigate the issue.
Another advantage of simplified instances is that they can be simple enough to be included in regular test sets to ensure that an error is not reintroduced later accidentally.
As an example, the instances generated for \scip can be found in its repository at \texttt{check/instances/Issue}, which was added to the regular ctest pipeline.

Delta debugging~\cite{ZNL99Deltadebugging} follows a hypothesis-trial-result approach to isolate the cause of a solver failure by simplifying the input data. This approach has been successfully applied to facilitate debugging SAT and SMT solvers~\cite{BrunmayerBiere-FuzzingAndDeltadebuggingSMT-Solver,BrunmayerLonsingBiere10AutomatedTestingAndDebuggingSAT,KaufmannBiere-TAP22,NiemetzBiere-SMT13,PaxianBiere-POS23}.

The primary goal of \mipdd is to simplify debugging.
The main approach is to iteratively simplify the problem instance while ensuring that the error remains reproducible.
This can lead to substantially shorter debug logs and also simpler numerics, which results in a streamlined debugging process.
As demonstrated in the case studies in Hoen et al.~\citep{HKG24MIPDD}, instances triggering fundamental bugs in \scip can often be reduced by \mipdd within a reasonable time to an instance comprising just a few variables and constraints.

The latest version \mipdd~2.0 features an automatic adaption of modification batches, which aims at keeping the total solving effort spent in a modifier run constant over the \mipdd run based on the solving effort of the current simplified instance.
Additionally, it features an automatic restriction of solving limits to terminate solving earlier if an issue should have occurred already based on the previous solving run.
Especially these features can significantly accelerate the delta debugging process.
Furthermore, it is shipped with an interface to \soplex as well as \scip, both for their real and exact solving modes, while there is a separate compilation parameter to select the arithmetic type of the internal data representation in order to improve reproducibility outside \mipdd.
Moreover, there is an option to use the implemented delta debugging algorithm to heuristically compute an IIS, similar to the functionality available in \scip (see Section~\ref{sec:iis}).

In the \scip releases 8.0.1 to 8.1.1 \mipdd significantly supported fixing 24 out of all 51 MIP-related bugs documented in the changelog~\citep{HKG24MIPDD}.
Up to now, \mipdd has evolved into a highly useful debugging tool that is actively used in current \scip development.

\mipdd is publicly available at \url{https://github.com/scipopt/MIP-DD}.
For a more detailed description of the basic features of \mipdd, please refer \citep{HKG24MIPDD}.


\section{Applications}
\label{sect:applications}

\scipv contains a new application that provides tailored functionality for
Pseudo-Boolean optimization problems: the PBSolver.
In the following, we provide the main features of this application.

To participate in the Pseudo-Boolean competition, a solver is required to strictly comply with the given regulations as prescribed by the organizers \cite{PB24}.
Especially, the solver has to produce a DIMACS CNF styled log output to correctly communicate with the competition environment.
For this reason, an application called PBSolver is created to make this interface public and maintained for simpler future submissions.
This application mainly provides the following basic features:
\begin{itemize}
\item message handler to prepend each standard log line by the comment specifier ``c'' in order to retrieve relevant solving statistics for posterior investigations
\item event handler to catch events of type \solver{SCIP\_EVENTTYPE\_BESTSOLFOUND} in order to directly signal a computation of a better primal solution by the solution specifier ``o'' followed by its precise objective value
\end{itemize}
Fixed-size arithmetic can only handle numbers up to a certain magnitude without loss of integral precision.
Therefore, the OPB reader was extended by the parameter \param{reading/opbreader/maxintsize}.
If an intsize is given in the instance header which exceeds the parameter value, the status specifier ``s UNSUPPORTED'' is displayed and the run is immediately terminated.
With this application both OPB and WBO instances can be solved out of the box in accordance to competition rules.
In the Pseudo-Boolean competition 2024, solvers incorporating \scip won the categories \textsc{opt-lin}, \textsc{partial-lin}, \textsc{soft-lin}, \textsc{opt-nlc}, and \textsc{dec-nlc}.
Out of a total of 1,207 instances, \scip successfully solved 759, while its parallelized variant FiberSCIP even solved 776.
Details on the algorithms useful for Pseudo-Boolean solving which are now implemented in \scip can be found in \cite{Mexi2025}.
Yet missing but planned to be added in the near future are a problem type detection to apply specific solving parameter settings, a dedicated feasibility definition that allows to solve also numerically challenging instances reliably, and an exact solving mode to also participate in verified categories next time.


\subsection*{Acknowledgments}

The authors want to thank all previous developers and contributors to the \scip Optimization Suite and
all users that reported bugs and often also helped reproducing and fixing the bugs.
Ambros Gleixner and Leon Eifler wish to thank Kati Jarck and Daniel Steffy for their prototype of the exact solving mode described in~\cite{CookKochSteffyetal2013,Jarck2020}, in particular their implementation of safe dual bounding described in Section~\ref{sec:exact:safelp}.
We like to thank Zuse Institute Berlin and RWTH Aachen for the computing infrastructure that was used for countless test runs and hosting the projects of the suite.

\subsection*{Contributions of the Authors}

The material presented in the article is highly related to code and software.
In the following we try to make the corresponding contributions of the authors and possible contact points more transparent.

\begin{enumerate*}
\item[Mathieu Besançon] implemented parts of the Probabilistic lookahead strong branching (see Section~\ref{sec:branching}) and maintains the Julia interface.
\item[Ksenia Bestuzheva] implemented the interface to the CONOPT nonlinear programming solver (see Section~\ref{sec:conopt}).
\item[Sander Borst] and Leon Eifler contributed safe dual proof analysis and constraint propagation (see Section~\ref{sec:exact:prop}).
\item[Antonia Chmiela] provided the first version of the repair mechanism described in Section~\ref{sec:splitExact1}.
\item[Jo\~{a}o Dion\'isio] has been developing and maintaining PySCIPOpt.
\item[Johannes Ehls] revised the implementation of \gcg's Cliquer pricing problem solver.
\item[Leon Eifler] was the main author of \scip's new exact solving mode, with contributions and revisions by Ambros Gleixner and Dominik Kamp, and a concerted code review effort of the whole team (see Section~\ref{sec:exact}).
\item[Mohammed Ghannam] implemented the serialization of SCIP statistics and is one of the developers of PySCIPOpt.
\item[Adrian G{\"o}{\ss}] implemented, tested, and finetuned the original heuristic (Section~\ref{subsubsec:ks}) and its extension (Section~\ref{subsubsec:dks}).
\item[Alexander Hoen,] Jacob von Holly-Ponientzietz and Dominik Kamp maintained \papilo and implemented the features described in Section~\ref{sect:papilo}.
Alexander Hoen and Dominik Kamp also developed \mipdd (see Section~\ref{sect:mipdd}).
\item[Christopher Hojny] implemented the generalization of symmetry handling methods to reflection symmetries as well as symmetry detection callbacks for Pseudo-Boolean and disjunctive constraints.
\item[Rolf van der Hulst] revised detection and handling of implied integrality in \scip (see Section~\ref{sec:presolve}).
\item[Dominik Kamp] provided numerous bugfixes all across \scip, \soplex, and \papilo.
\item[Thorsten Koch] added the permutation feature to ZIMPL (see Section~\ref{sect:zimpl}).
\item[Kevin Kofler] maintains the Java interface \solver{JSCIPOpt}.
\item[Jurgen Lentz] developed a draft of the \gcg structure and implemented the decomposition score plugin of \gcg. He is also a maintainer of PyCGCOpt.
\item[Marco L\"ubbecke] coordinated and supervised the development of GCG.
\item[Stephen J.\ Maher] implemented the updates to the Benders' decomposition framework (see Section~\ref{sec:benders}), the HiGHS pricing solver (Section~\ref{sec:gcg:highspricingsolver}), and the IPColGen heuristic in \gcg (Section~\ref{sec:gcg:ipcolgen}).
\item[Paul Matti Meinhold] implemented the dynamic batch sizing for the Greedy IIS finder of \scip.
\item[Gioni Mexi] implemented the cut-based conflict analysis (see Section~\ref{sec:conflict}), the probabilistic lookahead strong branching (see Section~\ref{sec:branching}), and together with Dominik Kamp and Alexander Hoen the Pseudo-Boolean application (see Section~\ref{sect:applications}).
\item[Til Mohr] implemented \gcg's extended master constraints interface and added the component bound branching rule to \gcg.
\item[Erik M\"uhmer] implemented the \gcg pricing problem solver, \gcg's new \texttt{JDEC} decomposition file format, and the \gcg structure. Furthermore, he revised the pricing parallelization and the extended master constraints interface.
\item[Krunal Kishor Patel] implemented the Ancestral pseudocosts branching rule.
\item[Marc E.\ Pfetsch] maintains SCIP-SDP, wrote the interface of \scip to \dejavu, and contributed at several places, including the symmetry code and diverse bug-fixes.
\item[Sebastian Pokutta] coordinated and supervised multiple developers at ZIB.
\item[Chantal Reinartz Groba] adapted \gcg's extended master constraints interface and added support for separators that work only on the master problem.
\item[Felipe Serrano] contributed to the development of the cut-based conflict analysis (see Section~\ref{sec:conflict}).
\item[Yuji Shinano] worked on the application for Pseudo-Boolean solving in UG (see Section~\ref{sect:ug}).
\item[Mark Turner] implemented the IIS finder plugin and was one of the maintainers of PySCIPOpt.
\item[Stefan Vigerske] implemented mixing branching on integer and nonlinear variables, the bound check for aggregated variables, and writing of \texttt{.nl} files.
\item[Matthias Walter] implemented the new multilinear separator (see Section~\ref{sec:separators}).
\item[Dieter Weninger] and Adrian G{\"o}{\ss} developed the ideas of the Decomposition Kernel Search heuristic (Section~\ref{subsubsec:dks}).
\item[Liding Xu] extended the LPI interface of PySCIPOpt.
\end{enumerate*}


\renewcommand{\refname}{\normalsize References}
\setlength{\bibsep}{0.25ex plus 0.3ex}
\DeclareRobustCommand{\VAN}[3]{#3}
\bibliographystyle{abbrvnat}

\begin{small}
\bibliography{scipopt}

@InProceedings{pyscipopt,
    author="Maher, Stephen and Miltenberger, Matthias and Pedroso, Jo{\~a}o Pedro and Rehfeldt, Daniel and Schwarz, Robert and Serrano, Felipe",
    editor="Greuel, Gert-Martin and Koch, Thorsten and Paule, Peter and Sommese, Andrew",
    title="{PySCIPOpt}: Mathematical Programming in {Python} with the {SCIP} {Optimization} {Suite}",
    booktitle="Mathematical Software -- ICMS 2016",
    year="2016",
    publisher="Springer International Publishing",
    address="Cham",
    pages="301--307"
}

@article{GleixnerGottwaldHoen23,
    author     = {Gleixner, Ambros and Gottwald, Leona and Hoen, Alexander},
    title      = {{PaPILO}: A Parallel Presolving Library for Integer and Linear Optimization with Multiprecision Support},
    year       = {2023},
    issue_date = {November-December 2023},
    publisher  = {INFORMS},
    address    = {Linthicum, MD, USA},
    volume     = {35},
    number     = {6},
    doi        = {10.1287/ijoc.2022.0171},
    journal    = {INFORMS Journal on Computing},
    month      = nov,
    pages      = {1329--1341},
    numpages   = {13}
}

@article{HKG24MIPDD,
    author     = {Hoen, Alexander and Kamp, Dominik and Gleixner, Ambros},
    title      = {{MIP-DD}: Delta Debugging for Mixed-Integer Programming Solvers},
    year       = {2025},
    publisher  = {INFORMS},
    address    = {Linthicum, MD, USA},
    doi        = {10.1287/ijoc.2024.0844},
    journal    = {INFORMS Journal on Computing}
}

@article{MIPLIB,
    author                   = {Gleixner, Ambros and Hendel, Gregor and Gamrath, Gerald and Achterberg, Tobias and Bastubbe, Michael and Berthold, Timo and Christophel, Philipp M. and Jarck, Kati and Koch, Thorsten and Linderoth, Jeff and L\"ubbecke, Marco and Mittelmann, Hans D. and Ozyurt, Derya and Ralphs, Ted K. and Salvagnin, Domenico and Shinano, Yuji},
    title                    = {{MIPLIB 2017: Data-Driven Compilation of the 6th Mixed-Integer Programming Library}},
    journal                  = {Mathematical Programming Computation},
    year                     = {2021},
    doi                      = {10.1007/s12532-020-00194-3}
}

@misc{PB24,
    author = {Olivier Roussel},
    title = {Pseudo-{B}oolean competition 2024},
    year = {2024},
    url = {http://www.cril.univ-artois.fr/PB24/},
    urldate = {2024-10-24},
}

@misc{PB9,
    author = {Olivier Roussel},
    title = {Pseudo-{B}oolean competition 2009},
    year = {2010},
    url = {http://www.cril.univ-artois.fr/PB09/},
    urldate = {2024-10-24},
}

@inproceedings{PDLP,
    author = {Applegate, David and Diaz, Mateo and Hinder, Oliver and Lu, Haihao and Lubin, Miles and O\textquotesingle Donoghue, Brendan and Schudy, Warren},
    booktitle = {Advances in Neural Information Processing Systems},
    editor = {M. Ranzato and A. Beygelzimer and Y. Dauphin and P.S. Liang and J. Wortman Vaughan},
    pages = {20243--20257},
    publisher = {Curran Associates, Inc.},
    title = {Practical Large-Scale Linear Programming using Primal-Dual Hybrid Gradient},
    url = {https://proceedings.neurips.cc/paper_files/paper/2021/file/a8fbbd3b11424ce032ba813493d95ad7-Paper.pdf},
    volume = {34},
    year = {2021}
}

@article{AchterbergBixbyGuetal.2019,
    author =	 {Achterberg, Tobias and Bixby, Robert E. and Gu, Zonghao and Rothberg, Edward and
                  Weninger, Dieter},
    title =	 {Presolve Reductions in Mixed Integer Programming},
    journal =	 {INFORMS Journal on Computing},
    volume =	 {32},
    number =	 {2},
    pages =	 {473--506},
    year =	 {2020},
    doi =		 {10.1287/ijoc.2018.0857}
}

@misc{Achterberg2022,
    author    = {Tobias Achterberg},
    title     = {Combinatorial Algorithms Used Inside a {MIP} Solver},
    journal   = {SEA 2022, 20th Symposium on Experimental Algorithms},
    url       = {https://sea2022.ifi.uni-heidelberg.de/talk_achterberg.pdf},
    year      = {2022},
    pages     = {49--54}
}

@inproceedings{BrunmayerBiere-FuzzingAndDeltadebuggingSMT-Solver,
    author = {Brummayer, Robert and Biere, Armin},
    title = {Fuzzing and delta-debugging {SMT} solvers},
    year = {2009},
    publisher = {Association for Computing Machinery},
    address = {New York, NY, USA},
    doi = {10.1145/1670412.1670413},
    booktitle = {Proceedings of the 7th International Workshop on Satisfiability Modulo Theories},
    pages = {1--5},
    numpages = {5},
    location = {Montreal, Canada},
    series = {SMT '09},
}

@InProceedings{BrunmayerLonsingBiere10AutomatedTestingAndDebuggingSAT,
    author="Brummayer, Robert and Lonsing, Florian and Biere, Armin",
    editor="Strichman, Ofer and Szeider, Stefan",
    title="Automated Testing and Debugging of {SAT} and {QBF} Solvers",
    booktitle="Theory and Applications of Satisfiability Testing -- SAT 2010",
    year="2010",
    publisher="Springer",
    address="Berlin, Heidelberg",
    pages="44--57"
}

@inproceedings{KaufmannBiere-TAP22,
    author    = {Daniela Kaufmann and Armin Biere},
    editor    = {Laura Kov{\'{a}}cs and	Karl Meinke},
    title     = {Fuzzing and Delta Debugging And-Inverter Graph Verification Tools},
    booktitle = {Tests and Proofs - 16th International Conference, {TAP} 2022, Held
        as Part of {STAF} 2022, Nantes, France, July 5, 2022, Proceedings},
    series    = {Lecture Notes in Computer Science},
    volume    = {13361},
    pages     = {69--88},
    publisher = {Springer},
    year      = {2022},
    doi       = {10.1007/978-3-031-09827-7\_5},
}

@inproceedings{PaxianBiere-POS23,
    author    = {Tobias Paxian and Armin Biere},
    editor    = {J\"{a}rvisalo, Matti and Le Berre, Daniel},
    title     = {Uncovering and Classifying Bugs in {MaxSAT} Solvers through Fuzzing and Delta Debugging},
    booktitle = {Proceedings of the 14th Internantional Workshop on Pragmatics of SAT
        Co-located with the 26th International Conference on Theory
        and Applicationas of Satisfiability Testing (SAT 2003),
        Alghero, Italy, July, 4, 2023},
    series    = {{CEUR} Workshop Proceedings},
    pages     = {59--71},
    volume    = {3545},
    publisher = {CEUR-WS.org},
    year      = {2023},
    url       = {http://ceur-ws.org/Vol-3545/paper5.pdf},
}

@inproceedings{NiemetzBiere-SMT13,
    author    = {Aina Niemetz and Armin Biere},
    editor    = {Roberto Bruttomesso and Alberto Griggio},
    title     = {{ddSMT: A Delta Debugger for the SMT-LIB v2 Format}},
    booktitle = {Proceedings of the 11th International Workshop on
    Satisfiability Modulo Theories, {SMT} 2013),
        affiliated with the 16th International Conference on
        Theory and Applications of Satisfiability Testing, {SAT} 2013,
        Helsinki, Finland, July 8-9, 2013},
    pages     = {36--45},
    year      = {2013},
}

@InProceedings{ZNL99Deltadebugging,
    author="Zeller, Andreas",
    editor="Nierstrasz, Oscar and Lemoine, Michel",
    title="Yesterday, my Program Worked. {Today}, it Does Not. {Why}?",
    booktitle="Software Engineering --- ESEC/FSE '99",
    year="1999",
    publisher="Springer Berlin Heidelberg",
    address="Berlin, Heidelberg",
    pages="253--267",
    isbn="978-3-540-48166-9"
}

@article{Hojny2024+,
  author =	 {Christopher Hojny},
  title =	 {Detecting and handling reflection symmetries in mixed-integer (nonlinear) programming and beyond},
  journal =	 {Mathematical Programming Compution},
  volume =	 {},
  pages =	 {},
  year =	 {2025},
  doi =		 {10.1007/s12532-025-00289-9}
}

@Article{LibertiOstrowski2014,
  author                   = {Leo Liberti and James Ostrowski},
  title                    = {Stabilizer-based symmetry breaking constraints for mathematical programs},
  journal                  = {Journal of Global Optimization},
  year                     = {2014},
  volume                   = {60},
  pages                    = {183--194}
}

@InProceedings{Salvagnin2018,
  author =	 "Salvagnin, Domenico",
  editor =	 "van Hoeve, Willem-Jan",
  title =	 "Symmetry Breaking Inequalities from the {Schreier-Sims} Table",
  booktitle =	 "Integration of Constraint Programming, Artificial Intelligence, and Operations
                  Research",
  year =	 "2018",
  publisher =	 "Springer International Publishing",
  address =	 "Cham",
  pages =	 "521--529",
}

@article{KaibelPfetsch2008,
  year                     = {2008},
  journal                  = {Mathematical Programming},
  volume                   = {114},
  number                   = {1},
  doi                      = {10.1007/s10107-006-0081-5},
  title                    = {Packing and partitioning orbitopes},
  publisher                = {Springer},
  author                   = {Volker Kaibel and Marc E. Pfetsch},
  pages                    = {1--36}
}

@Article{HojnyPfetsch2019,
  author =	 {Christopher Hojny and Marc E. Pfetsch},
  title =	 {Polytopes associated with symmetry handling},
  journal =	 {Mathematical Programming},
  year =	 {2019},
  volume =       {175},
  number =       {1},
  pages =        {197--240},
  doi =		 {10.1007/s10107-018-1239-7}
}

@Article{Hojny2020,
  author =	 {Christopher Hojny},
  title =	 {Packing, Partitioning, and Covering Symresacks},
  journal =	 {Discrete Applied Mathematics},
  year =	 {2020},
  volume =       {283},
  pages =        {689--717},
  doi =		 {10.1016/j.dam.2020.03.002}
}

@article{BendottiEtAl2021,
  Title =	 {Orbitopal fixing for the full (sub-)orbitope and application to the Unit
                  Commitment Problem},
  Author =	 {Pascale Bendotti and Pierre Fouilhoux and C\'ecile Rottner},
  Year =	 {2021},
  Journal =	 {Mathematical Programming},
  Volume =	 {186},
  Pages =	 {337--372},
  DOI =		 {10.1007/s10107-019-01457-1}
}

@article{KaibelEtAl2011,
  title                    = {Orbitopal fixing},
  journal                  = {Discrete Optimization},
  volume                   = {8},
  number                   = {4},
  pages                    = {595--610},
  year                     = {2011},
  doi                      = {10.1016/j.disopt.2011.07.001},
  author                   = {Volker Kaibel and Matthias Peinhardt and Marc E. Pfetsch}
}

@article{DoornmalenHojny2024a,
  author =	 {{\VAN{Doornmalen}{Van}{van}} Doornmalen, Jasper and Hojny, Christopher},
  title =	 {A unified framework for symmetry handling},
  journal =	 {Mathematical Programming},
  year =	 {2025},
  volume =	 {212},
  pages =	 {217--271},
  doi =		 {10.1007/s10107-024-02102-2}
}

@Inbook{Szabo2005,
  author =       "Szab{\'o}, P{\'e}ter G{\'a}bor and Mark{\'o}t, Mih{\'a}ly Csaba and Csendes,
                  Tibor",
  title =        "Global Optimization in Geometry --- Circle Packing into the Square",
  bookTitle =    "Essays and Surveys in Global Optimization",
  year =         "2005",
  publisher =    "Springer US",
  address =      "Boston, MA",
  pages =        "233--265",
  doi =          "10.1007/0-387-25570-2_9"
}

@article{Khajavirad2024,
  title =        {The circle packing problem: A theoretical comparison of various convexification
                  techniques},
  journal =      {Operations Research Letters},
  volume =       {57},
  year =         {2024},
  doi =          {10.1016/j.orl.2024.107197},
  author =       {Aida Khajavirad}
}

@article{Salvagnin2005,
  title   = {A dominance procedure for integer programming},
  author  = {Salvagnin, Domenico},
  journal = {Master’s thesis, University of Padova, Padova, Italy},
  year    = {2005}
}

@article{Liberti2012a,
  year                     = {2012},
  journal                  = {Mathematical Programming},
  volume                   = {131},
  number                   = {1--2},
  doi                      = {10.1007/s10107-010-0351-0},
  title                    = {Reformulations in mathematical programming: automatic symmetry detection and exploitation},
  publisher                = {Springer},
  author                   = {Liberti, Leo},
  pages                    = {273--304}
}

@PhdThesis{Achterberg2007a,
  Title =	 {Constraint Integer Programming},
  Author =	 {Tobias Achterberg},
  School =	 {Technische Universit{\"a}t Berlin},
  type =	 {Dissertation},
  Year =	 {2007}
}

@PhdThesis{Eifler2025,
  Title =	 {Algorithms and Certificates for Exact Mixed Integer Programming},
  Author =	 {Leon Eifler},
  School =	 {Technische Universit{\"a}t Berlin},
  type =	 {Dissertation},
  Year =	 {2025},
  url = {https://doi.org/10.14279/depositonce-23941}
}

@Article{Achterberg2009,
  Title                    = {{SCIP:} {S}olving {C}onstraint {I}nteger {P}rograms},
  Author                   = {Tobias Achterberg},
  Journal                  = {Mathematical Programming Computation},
  Year                     = {2009},
  Number                   = {1},
  Pages                    = {1--41},
  Volume                   = {1},
  Doi                      = {10.1007/s12532-008-0001-1},
  Address                  = {Berlin, Heidelberg},
  Publisher                = {Springer}
}

@InProceedings{GamrathLuebbecke2010,
  Title                    = {Experiments with a Generic {D}antzig-{W}olfe Decomposition for Integer Programs},
  Author                   = {Gerald Gamrath and Marco E. L{\"u}bbecke},
  Booktitle                = {Experimental Algorithms},
  Year                     = {2010},
  Editor                   = {Festa, Paola},
  Pages                    = {239--252},
  publisher                = {Springer Berlin Heidelberg},
  Series                   = {Lecture Notes in Computer Science},
  Volume                   = {6049},
  Doi                      = {10.1007/978-3-642-13193-6_21}
}

@InProceedings{Shinano2018,
  Title =	 {The {Ubiquity Generator} Framework: 7 Years of Progress in Parallelizing
                  Branch-and-Bound},
  Author =	 {Shinano, Yuji},
  Booktitle =	 {Operations Research Proceedings 2017},
  Year =	 {2018},
  Editor =	 {Kliewer, Natalia and Ehmke, Jan Fabian and Bornd{\"o}rfer, Ralf},
  Pages =	 {143--149},
  Publisher =	 {Springer},
  doi =		 "10.1007/978-3-319-89920-6_20"
}

@PhdThesis{Koch2004,
  author =	 {Thorsten Koch},
  title =	 {Rapid Mathematical Prototyping},
  school =	 {Technische Universit{\"a}t Berlin},
  type =	 {Dissertation},
  year =	 {2004}
}

@PhdThesis{Wunderling1996,
  Title =	 {Paralleler und objektorientierter {S}implex-{A}l\-go\-rith\-mus},
  Author =	 {Roland Wunderling},
  School =	 {Technische Universit{\"a}t Berlin},
  type =	 {Dissertation},
  Year =	 {1996}
}

@Article{GallyPfetschUlbrich2018,
  author =	 {Tristan Gally and Marc E. Pfetsch and Stefan Ulbrich},
  title =	 {A Framework for Solving Mixed-Integer Semidefinite Programs},
  journal =	 {Optimization Methods and Software},
  year =	 {2018},
  volume =	 {33},
  number =	 {3},
  pages =	 {594--632},
  doi =		 {10.1080/10556788.2017.1322081}
}

@article{legat2022mathoptinterface,
  title={{MathOptInterface}: a data structure for mathematical optimization problems},
  author={Legat, Beno{\^\i}t and Dowson, Oscar and Garcia, Joaquim Dias and Lubin, Miles},
  journal={INFORMS Journal on Computing},
  volume={34},
  number={2},
  pages={672--689},
  year={2022},
  publisher={INFORMS}
}

@Article{Angelelli2010,
	title={Kernel search: a general heuristic for the multi-dimensional knapsack problem},
	author={Angelelli, Enrico and Mansini, Renata and Speranza, M. Grazia},
	journal={Computers \& Operations Research},
	volume={37},
	number={11},
	pages={2017--2026},
	year={2010},
	publisher={Elsevier},
	doi = {10.1016/j.cor.2010.02.002},
}

@Article{Angelelli2012,
	title={Kernel search: A new heuristic framework for portfolio selection},
	author={Angelelli, Enrico and Mansini, Renata and Speranza, M. Grazia},
	journal={Computational Optimization and Applications},
	volume={51},
	number={1},
	pages={345--361},
	year={2012},
	publisher={Springer},
	doi = {10.1007/s10589-010-9326-6}
}

@Article{Guastaroba2017,
	title={Adaptive kernel search: a heuristic for solving mixed integer linear programs},
	author={Guastaroba, Gianfranco and Savelsbergh, M. and Speranza, M. Grazia},
	journal={European Journal of Operational Research},
	volume={263},
	number={3},
	pages={789--804},
	year={2017},
	publisher={Elsevier},
	doi = {10.1016/j.ejor.2017.06.005},
}

@Article{Halbig2025,
	title={Exploiting user-supplied Decompositions inside Heuristics},
	author= {Halbig, Katrin and G{\"o}{\ss}, Adrian and Weninger, Dieter},
	journal={Journal of Heuristics},
	volume={31},
	number={36},
	year={2025},
	publisher={Springer},
	doi={10.1007/s10732-025-09572-3}
}

@article{Nauty,
  title =	 {Practical graph isomorphism, {II}},
  journal =	 {Journal of Symbolic Computation},
  volume =	 {60},
  pages =	 {94--112},
  year =	 {2014},
  doi =		 {10.1016/j.jsc.2013.09.003},
  author =	 {Brendan D. McKay and Adolfo Piperno}
}

@InProceedings{JunttilaKaski2011,
  author =	 {Tommi Junttila and Petteri Kaski},
  title =	 {Conflict Propagation and Component Recursion for Canonical Labeling},
  booktitle =	 {Theory and Practice of Algorithms in (Computer) Systems -- First International
                  {ICST} Conference, {TAPAS} 2011, Rome, Italy, April 18--20, 2011. Proceedings},
  year =	 {2011},
  editor =	 {Alberto Marchetti{-}Spaccamela and Michael Segal},
  volume =	 {6595},
  series =	 {Lecture Notes in Computer Science},
  pages =	 {151--162},
  publisher =	 {Springer},
  doi =		 {10.1007/978-3-642-19754-3\_16},
}

@InProceedings{JunttilaKaski2007,
  author =	 {Tommi Junttila and Petteri Kaski},
  title =	 {Engineering an efficient canonical labeling tool for large and sparse graphs},
  booktitle =	 {Proceedings of the Ninth Workshop on Algorithm Engineering and Experiments and the
                  Fourth Workshop on Analytic Algorithms and Combinatorics},
  pages =	 {135--149},
  year =	 {2007},
  editor =	 {David Applegate and Gerth St{\o}lting Brodal and Daniel Panario and Robert
                  Sedgewick},
  publisher =	 {SIAM},
  doi =		 {10.1137/1.9781611972870.13},
}

@inproceedings{AndersS21,
  author =	 {Markus Anders and Pascal Schweitzer},
  title =	 {Parallel Computation of Combinatorial Symmetries},
  booktitle =	 {29th Annual European Symposium on Algorithms, {ESA} 2021, September 6-8, 2021,
                  Lisbon, Portugal (Virtual Conference)},
  series =	 {LIPIcs},
  volume =	 {204},
  pages =	 {6:1--6:18},
  publisher =	 {Schloss Dagstuhl - Leibniz-Zentrum f{\"{u}}r Informatik},
  year =	 {2021},
  doi =		 {10.4230/LIPIcs.ESA.2021.6},
}

@InProceedings{AndersSS2023,
  author =	 {Anders, Markus and Schweitzer, Pascal and Stie{\ss}, Julian},
  title =	 {{Engineering a Preprocessor for Symmetry Detection}},
  booktitle =	 {21st International Symposium on Experimental Algorithms (SEA 2023)},
  pages =	 {1:1--1:21},
  series =	 {Leibniz International Proceedings in Informatics (LIPIcs)},
  year =	 {2023},
  volume =	 {265},
  editor =	 {Georgiadis, Loukas},
  publisher =	 {Schloss Dagstuhl -- Leibniz-Zentrum f{\"u}r Informatik},
  address =	 {Dagstuhl, Germany},
  doi =		 {10.4230/LIPIcs.SEA.2023.1},
}

@PhdThesis{Anders2024,
  author = 	 {Markus Anders},
  title = 	 {Efficient Algorithms for Symmetry Detection},
  school = 	 {TU Darmstadt},
  year = 	 {2024}
}

@Article{DelPiaK17,
  author    = {Del Pia, Alberto and Khajavirad, Aida},
  title     = {{A Polyhedral Study of Binary Polynomial Programs}},
  doi       = {10.1287/moor.2016.0804},
  number    = {2},
  pages     = {389--410},
  volume    = {42},
  comment   = {cite for multilinear polytope},
  journal   = {Mathematics of Operations Research},
  year      = {2017},
}

@Article{DelPiaD21,
  author   = {Del Pia, Alberto and Di Gregorio, Silvia},
  title    = {Chv{\'a}tal Rank in Binary Polynomial Optimization},
  doi      = {10.1287/ijoo.2019.0049},
  comment  = {cite for multilinear polytope},
  journal  = {INFORMS Journal on Optimization},
  year     = {2021},
  pages    = {315--443}
}

@InProceedings{DelPiaW22,
  author    = {Del Pia, Alberto and Walter, Matthias},
  booktitle = {Integer Programming and Combinatorial Optimization},
  title     = {Simple Odd $\beta$-Cycle Inequalities for Binary Polynomial Optimization},
  doi       = {10.1007/978-3-031-06901-7_14},
  editor    = {Aardal, Karen and Sanit{\`a}, Laura},
  pages     = {181--194},
  publisher = {Springer International Publishing},
  comment   = {cite for multilinear polytope},
  year      = {2022},
}

@Article{CramaR17,
  author   = {Crama, Yves and Rodríguez-Heck, Elisabeth},
  title    = {A class of valid inequalities for multilinear 0--1 optimization problems},
  doi      = {10.1016/j.disopt.2017.02.001},
  pages    = {28--47},
  volume   = {25},
  comment  = {cite for multilinear polytope},
  journal  = {Discrete Optimization},
  year     = {2017},
}

@Article{DelPiaK18,
  author    = {Del Pia, Alberto and Khajavirad, Aida},
  title     = {The Multilinear Polytope for Acyclic Hypergraphs},
  doi       = {10.1137/16M1095998},
  number    = {2},
  pages     = {1049--1076},
  volume    = {28},
  comment   = {cite for multilinear polytope},
  journal   = {SIAM Journal on Optimization},
  year      = {2018},
}

@Article{DelPiaKS20,
  author    = {Del Pia, Alberto and Khajavirad, Aida and Sahinidis, Nikolaos V.},
  title     = {On the impact of running intersection inequalities for globally solving polynomial optimization problems},
  doi       = {10.1007/s12532-019-00169-z},
  number    = {2},
  pages     = {165--191},
  volume    = {12},
  comment   = {cite for multilinear polytopecite for multilinear polytope},
  journal   = {Mathematical Programming Computation},
  publisher = {Springer},
  year      = {2020},
}

@Book{Chinneck2008,
  Title =	 {Feasibility and infeasibility in optimization: algorithms and computational
                  methods},
  Author =	 {Chinneck, John W},
  Publisher =	 {Springer},
  Year =	 {2008},
  Series =	 {International Series in Operations Research and Management Sciences},
  Volume =	 {118}
}

@article{Guieu1999,
	title={Analyzing infeasible mixed-integer and integer linear programs},
	author={Guieu, Olivier and Chinneck, John W},
	journal={INFORMS Journal on Computing},
	volume={11},
	number={1},
	pages={63--77},
	year={1999},
	publisher={INFORMS}
}

@article{Puranik2017,
	title={Deletion presolve for accelerating infeasibility diagnosis in optimization models},
	author={Puranik, Yash and Sahinidis, Nikolaos V},
	journal={INFORMS Journal on Computing},
	volume={29},
	number={4},
	pages={754--766},
	year={2017},
	publisher={INFORMS}
}

@Article{Pfe08,
  Title =	 {Branch-And-Cut for the Maximum Feasible Subsystem Problem},
  Author =	 {Marc E. Pfetsch},
  Journal =	 {SIAM J. Optim.},
  Year =	 {2008},
  Number =	 {1},
  Pages =	 {21--38},
  Volume =	 {19},
  doi =		 {10.1137/050645828}
}

@TechReport{SCIP7,
  author =	 {Gerald Gamrath and Daniel Anderson and Ksenia Bestuzheva and Wei-Kun Chen and Leon
                  Eifler and Maxime Gasse and Patrick Gemander and Ambros Gleixner and Leona
                  Gottwald and Katrin Halbig and Gregor Hendel and Christopher Hojny and Thorsten
                  Koch and Pierre Le Bodic and Stephen J. Maher and Frederic Matter and Matthias
                  Miltenberger and Erik Mühmer and Benjamin Müller and Marc E. Pfetsch and Franziska
                  Schlösser and Felipe Serrano and Yuji Shinano and Christine Tawfik and Stefan
                  Vigerske and Fabian Wegscheider and Dieter Weninger and Jakob Witzig},
  title =	 {{The SCIP Optimization Suite 7.0}},
  institution =	 {Optimization Online},
  year =	 {2020},
  note =	 {\url{http://www.optimization-online.org/DB_HTML/2020/03/7705.html}}
}

@TechReport{SCIP8,
  author      = {Ksenia Bestuzheva and Mathieu Besan{\c{c}}on and Wei-Kun Chen and Antonia Chmiela and Tim Donkiewicz and Jasper van Doornmalen and Leon Eifler and Oliver Gaul and Gerald Gamrath and Ambros Gleixner and Leona Gottwald and Christoph Graczyk and Katrin Halbig and Alexander Hoen and Christopher Hojny and Rolf van der Hulst and Thorsten Koch and Marco L{\"u}bbecke and Stephen J. Maher and Frederic Matter and Erik M{\"u}hmer and Benjamin M{\"u}ller and Marc E. Pfetsch and Daniel Rehfeldt and Steffan Schlein and Franziska Schl{\"o}sser and Felipe Serrano and Yuji Shinano and Boro Sofranac and Mark Turner and Stefan Vigerske and Fabian Wegscheider and Philipp Wellner and Dieter Weninger and Jakob Witzig},
  title       = {{The SCIP Optimization Suite 8.0}},
  institution = {Zuse Institute Berlin},
  type = {ZIB-Report},
  number      = {21--41},
  year        = {2021}
}

@article{HuangfuHall15,
	author = {Huangfu, Q. and Hall, Julian},
	year = {2015},
	pages = {119--142},
	title = {Parallelizing the dual revised simplex method},
	volume = {10},
	journal = {Mathematical Programming Computation},
	doi = {10.1007/s12532-017-0130-5}
}

@TechReport{desrosiers2024branch,
  author =	 {Jacques Desrosiers and Marco L{\"u}bbecke and Guy
                  Desaulniers and Jean Bertrand Gauthier},
  title =	 {Branch-and-Price},
  institution =	 {GERAD},
  year =	 2024,
  type =	 {Les Cahiers du GERAD},
  number =	 {G-2024-36},
  note =	 {Forthcoming with Springer}
}

@Article{gcg:Vanderbeck:00,
  author =	 {Vanderbeck, Fran{\c{c}}ois},
  title =	 {On {D}antzig-{W}olfe Decomposition in Integer
                  Programming and Ways to Perform Branching in a
                  Branch-and-Price Algorithm},
  journal =	 {Operations Research},
  year =	 2000,
  volume =	 48,
  number =	 1,
  pages =	 {111--128}
}

@Article{maher2023ipcolgen,
  author =	 {Maher, S.J. and R\"{o}nnberg, E.},
  title =	 {Integer programming column generation: accelerating
                  branch-and-price using a novel pricing scheme for
                  finding high-quality solutions in set covering,
                  packing, and partitioning problems},
  journal =	 {Mathematical Programming Computation},
  year =	 2023,
  volume =	 15,
  pages =	 {509--548},
  doi =		 {10.1007/s12532-023-00240-w}
}

@article{maher2021benders,
  title={Implementing the branch-and-cut approach for a general purpose {Benders'} decomposition framework},
  author={Maher, Stephen J.},
  journal={European Journal of Operational Research},
  volume={290},
  number={2},
  pages={479--498},
  year={2021},
  publisher={Elsevier}
}

@incollection{Hoffman1956,
author = {Hoffman, A. J. and Kruskal, J. B.},
booktitle = {Linear Inequalities and Related Systems. (AM-38)},
doi = {10.1515/9781400881987-014},
month = {dec},
pages = {223--246},
publisher = {Princeton University Press},
title = {{Integral Boundary Points of Convex Polyhedra}},
volume = {38},
year = {1957}
}

@Unpublished{HulstW24,
author = {van der Hulst, Rolf and Walter, Matthias},
title = {A Row-wise Algorithm for Graph Realization},
note =	 {\url{https://optimization-online.org/?p=27423}},
year = {2024}
}

@Unpublished{HulstW25,
	title={Implied Integrality in Mixed-Integer Optimization}, 
	author={Rolf van der Hulst and Matthias Walter},
	year={2025},
	eprint={2504.07209},
	archivePrefix={arXiv},
	note =	 {\url{https://arxiv.org/abs/2504.07209}, To appear in "Integer Programming and Combinatorial Optimization 2025"},
}

@article{BixbyWagner1988,
author = {Bixby, Robert E and Wagner, Donald K},
journal = {Mathematics of Operation Research},
number = {1},
title = {{An almost linear-time algorithm for graph realization}},
volume = {13},
year = {1988}
}

@article{Truemper1990,
author = {Truemper, Klaus},
doi = {10.1016/0095-8956(90)90030-4},
journal = {Journal of Combinatorial Theory, Series B},
number = {2},
pages = {241--281},
title = {{A decomposition theory for matroids. V. Testing of matrix total unimodularity}},
volume = {49},
year = {1990}
}

@article{Walter2013,
author = {Walter, Matthias and Truemper, Klaus},
doi = {10.1007/s12532-012-0048-x},
journal = {Mathematical Programming Computation},
number = {1},
pages = {57--73},
title = {{Implementation of a unimodularity test}},
volume = {5},
year = {2013}
}

@BOOK{Schrijver86,
  title = {{Theory of Linear and Integer Programming}},
  publisher = {John Wiley \& Sons, Inc.},
  year = {1986},
  author = {Schrijver, Alexander},
  address = {New York, NY, USA},
}

@Article{EiflerGleixner2024,
  author   = {Eifler, Leon and Gleixner, Ambros},
  title    = {Safe and Verified {Gomory} Mixed-Integer Cuts in a Rational Mixed-Integer Program Framework},
  doi      = {10.1137/23M156046X},
  number   = {1},
  pages    = {742--763},
  volume   = {34},
  journal  = {SIAM Journal on Optimization},
  year     = {2024},
}

@misc{borst2024certifiedconstraintpropagationdual,
      title={Certified Constraint Propagation and Dual Proof Analysis in a Numerically Exact {MIP} Solver}, 
      author={Sander Borst and Leon Eifler and Ambros Gleixner},
      year={2024},
      eprint={2403.13567},
      archivePrefix={arXiv},
      primaryClass={math.OC},
      url={https://arxiv.org/abs/2403.13567}
}

@article{EiflerGleixner2022,
	author = {Eifler, Leon and Gleixner, Ambros},
	title = {A Computational Status Update for Exact Rational Mixed Integer Programming},
	journal = {Mathematical Programming},
	publisher = {Springer},
	year = {2023},
	volume = {197},
	issue = {2},
	pages = {793--812},
	doi = {10.1007/s10107-021-01749-5}
}

@phdthesis{Jarck2020,
  title  = {Exact mixed-integer programming},
  author = {Kati Jarck},
  year   = {2020},
  school = {TU Berlin},
  doi    = {10.14279/depositonce-9955}
}

@article{CookKochSteffyetal2013,
  author  = {William Cook and Thorsten Koch and Daniel E. Steffy and Kati Wolter},
  title   = {A hybrid branch-and-bound approach for exact rational mixed-integer programming},
  journal = {Mathematical Programming Computation},
  volume  = {5},
  number  = {3},
  pages   = {305--344},
  doi     = {10.1007/s12532-013-0055-6},
  year    = {2013}
}

@phdthesis{Espinoza2006,
  title  = {On Linear Programming, Integer Programming and Cutting
            Planes},
  author = {Daniel G. Espinoza},
  year   = {2006},
  school = {Georgia Institute of Technology}
}

@article{MarchandWolsey1998,
  author  = {Marchand, Hugues and Wolsey, Laurence},
  year    = {2001},
  pages   = {325--468},
  title   = {Aggregation and Mixed Integer Rounding to Solve {MIPs}},
  volume  = {49},
  journal = {Operations Research},
  doi     = {10.1287/opre.49.3.363.11211}
}

@misc{Mexi2025,
  title =	 {State-of-the-art Methods for Pseudo-Boolean Solving with {SCIP}},
  author =	 {Gioni Mexi and Dominik Kamp and Yuji Shinano and Shanwen Pu and Alexander Hoen and
                  Ksenia Bestuzheva and Christopher Hojny and Matthias Walter and Marc E. Pfetsch
                  and Sebastian Pokutta and Thorsten Koch},
  year =	 {2025},
  eprint =	 {2501.03390},
  archivePrefix ={arXiv},
  primaryClass = {math.OC},
  howpublished = {\url{https://doi.org/10.48550/arXiv.2501.03390}}
}

@inproceedings{VIPR,
  title        = {Verifying {Integer} {Programming} {Results}},
  author       = {Cheung, Kevin and Gleixner, Ambros and Steffy, Daniel E.},
  booktitle    = {International Conference on Integer Programming and Combinatorial Optimization},
  pages        = {148--160},
  year         = {2017},
  organization = {Springer},
  doi          = {10.1007/978-3-319-59250-3_13}
}

@InProceedings{HoehnOertelGleixnerNordstrom2024,
author="Hoen, Alexander
and Oertel, Andy
and Gleixner, Ambros
and Nordstr{\"o}m, Jakob",
editor="Dilkina, Bistra",
title="Certifying {MIP}-Based Presolve Reductions for 0-1 Integer Linear Programs",
booktitle="Integration of Constraint Programming, Artificial Intelligence, and Operations Research",
year="2024",
publisher="Springer Nature Switzerland",
address="Cham",
pages="310--328",
doi="10.1007/978-3-031-60597-0_20"
}

@techreport{Gomory1960,
  author      = {Gomory, Ralph E.},
  year        = {1960},
  institution = {RAND Corporation},
  title       = {An algorithm for the mixed integer problem},
  location    = {Santa Monica, CA}
}

@article{CookDashFukasawaGoycoolea2009,
  author  = {Cook, William and Dash, Sanjeeb and Fukasawa, Ricardo and Goycoolea, Marcos},
  year    = {2009},
  pages   = {641--649},
  title   = {Numerically Safe {Gomory} Mixed-Integer Cuts},
  volume  = {21},
  journal = {INFORMS Journal on Computing},
  doi     = {10.1287/ijoc.1090.0324}
}

@misc{GithubCakeML,
  author       = {Yong Kiam and Magnus Myreen},
  title        = {{Formalisation of VIPR in CakeML}},
  howpublished = {\url{https://github.com/CakeML/cakeml/tree/master/examples/vipr} (accessed May 1, 2024)},
  year         = {2024}
}

@misc{scipdex,
  author       = {João Dionísio and Mohammed Ghannam},
  title        = {{scipdex}},
  howpublished = {\url{https://github.com/mmghannam/scipdex} (accessed May 2, 2025)},
  year         = {2025}
}

@misc{pyscipoptDocu2025,
  author       = {Mark Turner},
  title        = {{PySCIPOpt Documentation}},
  howpublished = {\url{https://pyscipopt.readthedocs.io/en/latest/index.html} (accessed May 2, 2025)},
  year         = {2025}
}

@misc{scipppgithub,
  key = {SCIP++},
  title = {The {C++} Interface for {SCIP}},
  howpublished={\url{https://www.github.com/scipopt/SCIPpp}},
  year = {2023}
}

@misc{pysoplexgithub,
  key = {PySoPlex},
  title = {{The} {Python} Interface for {SoPlex}},
  howpublished={\url{https://www.github.com/scipopt/PySoPlex}},
  year = {2023}
}

@misc{papilojl,
  key = {PaPILO.jl},
  title = {{The} {Julia} Interface for {PaPILO}},
  howpublished={\url{https://www.github.com/scipopt/PaPILO.jl}},
  year = {2025}
}

@misc{matlabgithub,
  key = {MatlabSCIPInterface},
  title = {{The} {Matlab} Interface to {SCIP} and {SCIP-SDP}},
  howpublished={\url{https://github.com/scipopt/MatlabSCIPInterface}},
  year = {2023}
}

@misc{javagithub,
  key = {JSCIPOpt},
  title = {{JSCIPOpt}: {The} {Java} Interface to {SCIP}},
  howpublished={\url{https://github.com/scipopt/JSCIPOpt}},
  year = {2023}
}

@inproceedings{mexi2023probabilistic,
      title={Probabilistic Lookahead Strong Branching via a Stochastic Abstract Branching Model}, 
      author={Gioni Mexi and Somayeh Shamsi and Mathieu Besançon and Pierre Le Bodic},
      year={2024},
      booktitle={International Conference on Integration of Constraint Programming, Artificial Intelligence, and Operations Research},
      eprint={2312.07041},
      archivePrefix={arXiv},
      primaryClass={math.OC}
}

@article{achterberg2007conflict,
  title={Conflict analysis in mixed integer programming},
  author={Achterberg, Tobias},
  journal={Discrete Optimization},
  volume={4},
  number={1},
  pages={4--20},
  year={2007},
  publisher={Elsevier}
}

@article{mexi2024cut,
  title={Cut-based Conflict Analysis in Mixed Integer Programming},
  author={Mexi, Gioni and Serrano, Felipe and Berthold, Timo and Gleixner, Ambros and Nordstr{\"o}m, Jakob},
  journal={arXiv preprint arXiv:2410.15110},
  year={2024}
}

@inproceedings{glankwamdee2006lookahead,
	title={Lookahead branching for mixed integer programming},
	author={Glankwamdee, Wasu and Linderoth, Jeff},
	booktitle={Twelfth INFORMS Computing Society Meeting},
	pages={130--150},
	year={2006}
}

@article{achterberg2005reliability,
	title={Branching rules revisited},
	author={Achterberg, Tobias and Koch, Thorsten and Martin, Alexander},
	journal={Operations Research Letters},
	volume={33},
	number={1},
	pages={42--54},
	year={2005},
	publisher={Elsevier}
}

@article{benichou1971pscost,
	title={Experiments in mixed-integer linear programming},
	author={B{\'e}nichou, Michel and Gauthier, Jean-Michel and Girodet, Paul and Hentges, Gerard and Ribi{\`e}re, Gerard and Vincent, Olivier},
	journal={Mathematical Programming},
	volume={1},
	pages={76--94},
	year={1971},
	publisher={Springer}
}

@inproceedings{mexi23improving,
  title={Improving {Conflict} {Analysis} in {MIP} {Solvers} by {Pseudo-Boolean} {Reasoning}},
  author={Gioni Mexi and Timo Berthold and Ambros Gleixner and Jakob Nordstr{\"o}m},
  booktitle={Proceedings of the 29th International Conference on Principles and Practice of Constraint Programming (CP’23)},
  volume={280},
  pages={27},
  year={2023}
}

@article{le2011sat4j,
  title={The {Sat4j} library, release 2.2: System description},
  author={Le Berre, Daniel and Parrain, Anne},
  journal={Journal on Satisfiability, Boolean Modelling and Computation},
  volume={7},
  number={2-3},
  pages={59--64},
  year={2011},
  publisher={SAGE Publications Sage UK: London, England}
}

@inproceedings{elffers2018divide,
  title={Divide and conquer: Towards faster pseudo-boolean solving.},
  author={Elffers, Jan and Nordstr{\"o}m, Jakob},
  booktitle={IJCAI},
  volume={18},
  pages={1291--1299},
  year={2018}
}

@inproceedings{chai2003fast,
  title={A fast pseudo-boolean constraint solver},
  author={Chai, Donald and Kuehlmann, Andreas},
  booktitle={Proceedings of the 40th annual Design Automation Conference},
  pages={830--835},
  year={2003}
}

@misc{gmplib,
  key = {GMP},
  title = {{GMP: The GNU multiple precision arithmetic library}},
  howpublished={\url{https://gmplib.org} (accessed June 21, 2025)},
  year = {2025}
}

@misc{mpfrlib,
  key = {MPFR},
  title = {{The GNU MPFR library}},
  howpublished={\url{https://www.mpfr.org/} (accessed June 21, 2025)},
  year = {2025}
}

@article{mpfrpaper,
author = {Fousse, Laurent and Hanrot, Guillaume and Lef\`{e}vre, Vincent and P\'{e}lissier, Patrick and Zimmermann, Paul},
title = {{MPFR}: A multiple-precision binary floating-point library with correct rounding},
year = {2007},
issue_date = {June 2007},
publisher = {Association for Computing Machinery},
address = {New York, NY, USA},
volume = {33},
number = {2},
note = {\url{https://doi.org/10.1145/1236463.1236468}},
journal = {ACM Trans. Math. Softw.},
month = jun,
pages = {13–es},
numpages = {15}
}

@misc{boostmplib,
  key = {BOOSTMP},
  title = {{The Boost Multiprecision Library}},
  howpublished={\url{https://github.com/boostorg/multiprecision} (accessed June 21, 2025)},
  year = {2025}
}

@misc{viprgithub,
  key = {VIPR},
  title = {{VIPR: Verifying Integer Programming Results}},
  howpublished={\url{https://github.com/scipopt/vipr} (accessed June 21, 2025)},
  year = {2025}
}

@article{witzig_computational_2021,
  title = {Computational Aspects of Infeasibility Analysis in Mixed Integer Programming},
  author = {Witzig, Jakob and Berthold, Timo and Heinz, Stefan},
  year = {2021},
  month = mar,
  journal = {Mathematical Programming Computation},
  doi = {10.1007/s12532-021-00202-0},
  langid = {english},
}

@InProceedings{gcg:RyanFoster:81,
  author =	 {D.M. Ryan and B.A. Foster},
  title =	 {An Integer Programming Approach to Scheduling},
  booktitle =	 {Computer Scheduling of Public Transport Urban
                  Passenger Vehicle and Crew Scheduling},
  pages =	 {269--280},
  year =	 1981,
  editor =	 {A. Wren},
  address =	 {Amsterdam},
  publisher =	 {North-Holland}
}

@article{parker1981approximation,
  title={Approximation programming of chemical processes—2: computational difficulties},
  author={Parker, Arthur L and Hughes, Richard R},
  journal={Computers \& Chemical Engineering},
  volume={5},
  number={3},
  pages={135--141},
  year={1981},
  publisher={Elsevier}
}

@article{goldfarb1983numerically,
  title={A numerically stable dual method for solving strictly convex quadratic programs},
  author={Goldfarb, Donald and Idnani, Ashok},
  journal={Mathematical programming},
  volume={27},
  number={1},
  pages={1--33},
  year={1983},
  publisher={Springer}
}

@article{gill1984weighted,
  title={A weighted {G}ram-{S}chmidt method for convex quadratic programming},
  author={Gill, Philip E and Gould, Nicholas IM and Murray, Walter and Saunders, Michael A and Wright, Margaret H},
  journal={Mathematical programming},
  volume={30},
  number={2},
  pages={176--195},
  year={1984},
  publisher={Springer}
}

\end{small}

\subsection*{Author Affiliations}

\hypersetup{urlcolor=black}
\newcommand{\myorcid}[1]{ORCID: \href{https://orcid.org/#1}{#1}}
\newcommand{\myemail}[1]{E-mail: \href{#1}{#1}}
\newcommand{\myaffil}[2]{{\noindent #1}\\{#2}\bigskip}

\small

\myaffil{Christopher Hojny}{%
  Technische Universiteit Eindhoven, Department of Mathematics and Computer Science, P.O.\ Box 513, 5600 MB Eindhoven, The Netherlands\\
  \myemail{c.hojny@tue.nl}\\
  \myorcid{0000-0002-5324-8996}}

\myaffil{Mathieu Besançon}{%
  Université Grenoble Alpes, Inria, Laboratoire d'Informatique de Grenoble, 38100~Grenoble, France\\
  \myemail{mathieu.besancon@inria.fr}\\
  \myorcid{0000-0002-6284-3033}}

\myaffil{Ksenia Bestuzheva}{%
  GAMS Software GmbH\\
  \myemail{kbestuzheva@gams.com}\\
  \myorcid{0000-0002-7018-7099}}

\myaffil{Sander Borst}{%
  Max Planck Institute for Informatics, Saarland Informatics Campus, Campus
  E1 4, 66123 Saarbr\"ucken, Germany \\
  \myemail{sborst@mpi-inf.mpg.de}\\
  \myorcid{0000-0003-4001-6675}}

\myaffil{Antonia Chmiela}{%
  Zuse Institute Berlin, Department AIS$^2$T, Takustr.~7, 14195~Berlin, Germany\\
  \myemail{chmiela@zib.de}\\
  \myorcid{0000-0002-4809-2958}}

\myaffil{Jo{\~a}o~Dion{\'i}sio}{%
  University of Porto, Faculty of Sciences, Rua do Campo Alegre, 4169-007 Porto, Portugal\\
  and Zuse Institute Berlin, Takustr.~7, 14195~Berlin, Germany\\
  \myemail{joao.goncalves.dionisio@gmail.com}\\
  \myorcid{0009-0005-5160-0203}}

\myaffil{Johannes Ehls}{%
  RWTH Aachen University, Lehrstuhl f\"ur Operations Research, Kackertstr.~7, 52072 Aachen, Germany\\
  \myemail{johannes.ehls@rwth-aachen.de}\\
  \myorcid{0009-0005-1130-6683}}

\myaffil{Leon Eifler}{%
  Zuse Institute Berlin, Department AIS$^2$T, Takustr.~7, 14195~Berlin, Germany\\
  \myemail{eifler@zib.de}\\
  \myorcid{0000-0003-0245-9344}}

\myaffil{Mohammed Ghannam}{%
    Zuse Institute Berlin, Takustr.~7, 14195~Berlin, Germany\\
    \myemail{ghannam@zib.de}
}

\myaffil{Ambros Gleixner}{%
  Hochschule für Technik und Wirtschaft Berlin, 10313~Berlin, Germany\\
  and Zuse Institute Berlin, Takustr.~7, 14195~Berlin, Germany\\
  \myemail{gleixner@htw-berlin.de}\\
  \myorcid{0000-0003-0391-5903}}

\myaffil{Adrian G{\"o}{\ss}}{%
  University of Technology Nuremberg, Dr.-Luise-Herzberg-Str.~4, 90461 Nuremberg, Germany\\
  \myemail{adrian.goess@utn.de}\\
  \myorcid{0009-0002-7144-8657}}

\myaffil{Alexander Hoen}{%
  Hochschule für Technik und Wirtschaft Berlin, 10313~Berlin, Germany\\
  and Zuse Institute Berlin, Takustr.~7, 14195~Berlin, Germany\\  \myemail{hoen@zib.de}\\
  \myorcid{0000-0003-1065-1651}}

\myaffil{Jacob von Holly-Ponientzietz}{%
  Zuse Institute Berlin, Takustr.~7, 14195~Berlin, Germany\\
  \myemail{von.holly-ponientzietz@zib.de}\\
  \myorcid{0009-0002-2601-3689}}

\myaffil{Rolf van der Hulst}{%
  University of Twente, Department of Applied Mathematics,  P.O.~Box~217, 7500~AE~Enschede, The Netherlands\\
  \myemail{r.p.vanderhulst@utwente.nl}\\
  \myorcid{0000-0002-5941-3016}}

\myaffil{Dominik Kamp}{%
  University of Bayreuth, Chair of Economathematics, Universitaetsstrasse 30, 95440 Bayreuth, Germany\\
  \myemail{dominik.kamp@uni-bayreuth.de}\\
  \myorcid{0009-0005-5577-9992}}

\myaffil{Thorsten Koch}{%
  Technische Universit\"at Berlin, Chair of Software and Algorithms for Discrete Optimization, Stra\ss{}e des 17. Juni 135, 10623 Berlin, Germany, and\\
  Zuse Institute Berlin, Department A$^2$IM, Takustr. 7, 14195~Berlin, Germany\\
  \myemail{koch@zib.de}\\
  \myorcid{0000-0002-1967-0077}}

\myaffil{Kevin Kofler}{%
  DAGOPT Optimization Technologies GmbH\\
  \myemail{kofler@dagopt.com}}

\myaffil{Jurgen Lentz}{%
  RWTH Aachen University, Lehrstuhl f\"ur Operations Research, Kackertstr.~7, 52072~Aachen, Germany\\
  \myemail{jurgen.lentz@rwth-aachen.de}\\
  \myorcid{0009-0000-0531-412X}}

\myaffil{Marco L\"ubbecke}{%
  RWTH Aachen University, Lehrstuhl f\"ur Operations Research, Kackertstr.~7, 52072 Aachen, Germany\\
  \myemail{marco.luebbecke@rwth-aachen.de}\\
  \myorcid{0000-0002-2635-0522}}

\myaffil{Stephen J. Maher}{%
  GAMS Software GmbH, Germany\\
  \myemail{smaher@gams.com}\\
  \myorcid{0000-0003-3773-6882}}

\myaffil{Paul Matti Meinhold}{%
  Zuse Institute Berlin, Takustr.~7, 14195~Berlin, Germany\\
  \myemail{meinhold@zib.de}\\
  \myorcid{0009-0003-5477-9152}}

\myaffil{Gioni Mexi}{%
  Zuse Institute Berlin, Takustr.~7, 14195~Berlin, Germany\\
  \myemail{mexi@zib.de}
  }

\myaffil{Til Mohr}{%
  RWTH Aachen University, Lehrstuhl f\"ur Operations Research, Kackertstr.~7, 52072 Aachen, Germany\\
  \myemail{til.mohr@rwth-aachen.de}\\
  \myorcid{0009-0001-9842-210X}}

\myaffil{Erik M\"uhmer}{%
  RWTH Aachen University, Lehrstuhl f\"ur Operations Research, Kackertstr.~7, 52072 Aachen, Germany\\
  \myemail{erik.muehmer@rwth-aachen.de}\\
  \myorcid{0000-0003-1114-3800}}

\myaffil{Krunal Kishor Patel}{%
  CERC, Polytechnique Montr\'eal, 2500 Chemin de Polytechnique, Montr\'eal,
  H3T 1J4, QC, Canada\\
  \myemail{krunal.patel@polymtl.ca}\\
  \myorcid{0000-0001-7414-5040}
}

\myaffil{Marc E.~Pfetsch}{%
  Technische Universität Darmstadt, Fachbereich Mathematik, Dolivostr.~15, 64293~Darmstadt, Germany\\
  \myemail{pfetsch@mathematik.tu-darmstadt.de}\\
  \myorcid{0000-0002-0947-7193}}

\myaffil{Sebastian Pokutta}{%
  Technische Universit\"at Berlin, Stra\ss{}e des 17. Juni 135, 10623 Berlin, Germany, and\\
  Zuse Institute Berlin, Takustr.~7, 14195~Berlin, Germany\\
  \myemail{pokutta@zib.de}\\
  \myorcid{0000-0001-7365-3000}}

\myaffil{Chantal Reinartz Groba}{%
  RWTH Aachen University, Lehrstuhl f\"ur Operations Research, Kackertstr.~7, 52072~Aachen, Germany\\
  \myemail{chantal.michelle.reinartz@rwth-aachen.de}\\
  \myorcid{0009-0001-1820-3864}}

\myaffil{Felipe Serrano}{%
  COPT GmbH, Berlin, Germany\\
  \myemail{serrano@copt.de}\\
  \myorcid{0000-0002-7892-3951}}

\myaffil{Yuji Shinano}{%
	Zuse Institute Berlin, Takustr.~7, 14195~Berlin, Germany\\
	\myemail{shinano@zib.de}\\
	\myorcid{0000-0002-2902-882X}}

\myaffil{Mark Turner}{%
  Zuse Institute Berlin, Takustr.~7, 14195~Berlin, Germany\\
  \myemail{turner@zib.de}\\
  \myorcid{0000-0001-7270-1496}}

\myaffil{Stefan Vigerske}{%
  GAMS Software GmbH, c/o Zuse Institute Berlin, Takustr.~7, 14195~Berlin, Germany\\
  \myemail{svigerske@gams.com}\\
  \myorcid{0009-0001-2262-0601}}

\myaffil{Matthias Walter}{%
  University of Twente, Department of Applied Mathematics,  P.O.~Box~217, 7500~AE~Enschede, The Netherlands\\
  \myemail{m.walter@utwente.nl}\\
  \myorcid{0000-0002-6615-5983}}

\myaffil{Dieter Weninger}{%
  Friedrich-Alexander-Universität Erlangen-Nürnberg, Cauerstr.~11, 91058 Erlangen, Germany\\
  \myemail{dieter.weninger@fau.de}\\
  \myorcid{0000-0002-1333-8591}}

\myaffil{Liding Xu}{%
  Zuse Institute Berlin, Department AIS$^2$T, Takustr.~7, 14195~Berlin, Germany\\
  \myemail{liding.xu@zib.de}\\
  \myorcid{0000-0002-0286-1109}}


\end{document}